\documentclass[final,1p,times,number]{elsarticle}
\journal{Journal for revision}
\usepackage[colorlinks]{hyperref} 

\usepackage{amssymb} 
\usepackage{amsmath}
\usepackage{amsthm}
\usepackage{array}
\usepackage{color}
\usepackage{algorithm}
\usepackage{algpseudocode}
\newtheorem{theorem}{Theorem}
\newtheorem{corollary}{Corollary}
\newtheorem{lemma}{Lemma}
\newtheorem{definition}{Definition}
\newtheorem{example}{Example}
\newtheorem{proposition}{Proposition}
\newtheorem{remark}{Remark}

\DeclareMathOperator*{\maxi}{max}

\begin{document}

\begin{frontmatter}

\title{Formal Reduction of Singularly-Perturbed Linear Differential Systems}

\author{Moulay A. Barkatou}
\address{XLIM UMR 7252 ; DMI\\
         University of Limoges; CNRS\\
        123, Avenue Albert Thomas\\
     87060 Limoges, France\\
         moulay.barkatou@unilim.fr}
         \author{Suzy S. Maddah\footnote{A part of this work was developed within the author's Ph.D. thesis at the University of Limoges, DMI.}}
\address{INRIA Saclay $\hat{I}$le-de-France\\
B$\hat{a}$timent Alan Turing\\
1 rue Honoré d'Estienne d'Orves\\
     91120 Palaiseau, France\\
         suzy.maddah@inria.fr}
\begin{abstract}
In this article, we discuss formal invariants of singularly-perturbed linear differential systems in neighborhood of turning points and give algorithms which allow their computation. The algorithms proposed are implemented in the computer algebra system \textsc{Maple}. 
\end{abstract}


\end{frontmatter}

Let $\partial = \frac{d}{dx}$ and consider the decisively simple-looking differential equation~ \cite[Introduction, Subsection 1.5 ]{key214}
$$ \varepsilon^2\; \partial^2 f^2 \;-\; (x^3 \;-\; \varepsilon) f\; =\; 0 $$
over the $(x,  \varepsilon)$-region $ D:\; |x| \leq \alpha_0\; , \; 0 < |\varepsilon| \leq \varepsilon_0$ for some positive constants $\alpha_0$ and $\varepsilon_0$. Setting $F= \begin{bmatrix} f \\ \varepsilon \partial f \end{bmatrix},$ this equation can be rewritten as the following first-order system:
\begin{equation} \label{introexmsys} \varepsilon \; \partial F\;=\; A(x, \varepsilon)\; F \;=\; \begin{bmatrix} 0 & 1 \\ x^3 - \varepsilon & 0 \end{bmatrix} F. \end{equation}
One can observe that the Jordan form of $A_0(x):= A(x,0)$ is not stable in any neighborhood of $x=0$. We refer by \textit{turning points} to such points where the Jordan form of $A_0(x)$ changes, i.e., either the multiplicity of the eigenvalues or the degrees of the elementary divisors change in a neighborhood of such points~\cite[p. 57]{key59}. For example, $x=0$ might be a \textit{turning point} for~\eqref{introexmsys}. Our goal is to compute a formal solution in a neighborhood of $x=0$ as $\varepsilon \rightarrow 0$. 
\begin{figure}
\centering
\includegraphics[width=60mm]{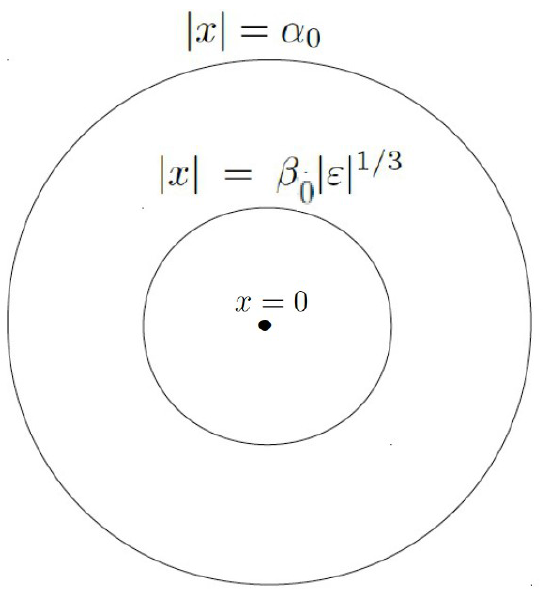}
\caption{}
\label{d1}
\end{figure}
\begin{itemize} 
\item In the region $\; \beta_0 |\varepsilon|^{1/3} \leq |x| \leq \alpha_0$ where $\beta_0$ is a positive constant, the transformation $F= T_1 G$ where $T_1=\begin{bmatrix}  1 & 0 \\ 0  & x^{3/2} \end{bmatrix}$ transforms system \eqref{introexmsys} into
$$(x^{-3} \varepsilon) x^{3/2} \; \partial G\;=\; \{ \begin{bmatrix} 0 & 1 \\ 1 & 0 \end{bmatrix} + (x^{-3} \varepsilon) \begin{bmatrix} 0 & 0 \\ -1 & \frac{-3}{2} x^{1/2} \end{bmatrix} \} G , $$
which can be rewritten, by a self-explanatory change of notation $x^{-3} \varepsilon = \xi$, as 
$$ \xi x^{3/2} \partial G \;=\; \{ \begin{bmatrix} 0 & 1 \\ 1 & 0 \end{bmatrix} + \xi \begin{bmatrix} 0 & 0 \\ -1 & \frac{-3}{2} x^{1/2} \end{bmatrix} \} G .$$
A transformation $ G = T_2 U$ where $T_2= \begin{bmatrix}  \frac{1}{2} & \frac{1}{2} \\ -\frac{1}{2}  & \frac{1}{2} \end{bmatrix} + O (\xi)$ would then result in the block-diagonalized system
$$ \xi x^{3/2} \partial U \;=\; \{ \begin{bmatrix} -1 & 0 \\ 0 & 1 \end{bmatrix} + O(\xi)  \} U.$$
Consequently, the system splits into two first-order linear differential scalar equations. We can then write down the following fundamental matrix of formal solutions for \eqref{introexmsys} in this region: $$F_{outer} = T_1\; T_2 \; exp (\begin{bmatrix} - \frac{2}{5}x^{5/2} \varepsilon^{-1} + O(x^{-1/2}) & 0 \\ 0 & \frac{2}{5}x^{5/2} \varepsilon^{-1} + O(x^{-1/2})\end{bmatrix}) .$$
\item In the region $|x| \leq \beta_0 |\varepsilon|^{1/3}$, we shall perform a so-called stretching transformation. Namely we set $\tau = x \varepsilon^{-1/3}$ and $\partial_{\tau} = \frac{d}{d\tau}$. Then, for all $\tau$ such as $|\tau| < \infty$, except possibly for neighborhoods of the roots of $\tau^3 - 1=0$\footnote{The turning points other than $\tau =0$, namely the roots of this equation do not explicitly correspond to the original turning point. They are referred to as \text{secondary turning points}~\cite{key214}.}, the transformation $F = L_1\; G$ where $$L_1 = \begin{bmatrix}  1 & 0 \\ 0  & \varepsilon^{1/2} \end{bmatrix}$$ reduces system~\eqref{introexmsys} to
\begin{equation} \label{introexmsysg} \varepsilon^{1/6} \partial_{\tau} \; G\;=\;  \begin{bmatrix} 0 & 1 \\ \tau^3 - 1 & 0 \end{bmatrix}  G .\end{equation} 
And the transformation $G = L_2 \;U$ where $$L_2 =  \left[ \begin {array}{cc} \frac{1}{2}+\frac{1}{8}\,{\tau}^{3} + O (\tau^{6})&\frac{1}{2
}+\frac{1}{8}\,{\tau}^{3} + O (\tau^{6})
\\ \noalign{\medskip}-\frac{i}{2}+\frac{i}{8}{\omega}^{3} + O (\tau^{6})& \frac{i}{2} -\frac{i}{8}{\tau}^{3} + O (\tau^{6})\end {array}
 \right] ,$$ reduces system~\eqref{introexmsysg} to 
$$\varepsilon^{1/6} \partial_{\tau} U = \{ \begin{bmatrix} -i + \frac{i}{2}\tau^3 + O (\tau^{6}) & 0 \\ 0 & i - \frac{i}{2}\tau^3 + O (\tau^{6}) \end{bmatrix} + O (\varepsilon^{1/6}) \} U .$$
A fundamental matrix of formal solution of system~\eqref{introexmsys} is then given by:
$$F_{inner} = L_1\; L_2 \; exp (\begin{bmatrix} \frac{-i \tau + (i/8)\tau^4 + O (\tau^{7}) }{ \varepsilon^{1/6} }& 0 \\ 0 & \frac{i \tau - (i/8)\tau^4 + O (\tau^{7}) }{ \varepsilon^{1/6} }\end{bmatrix}) .$$
\end{itemize}

We call $F_{outer}$ (resp. $F_{inner}$) an outer (resp. inner) solution as it is obtained in the outer (resp. inner) region around $x=0$. 
The corresponding differential systems are sometimes referred to as outer and inner differential systems as well. However, unlike the above particular example, intermediate regions might be encountered.
 
In this article, we study the construction of formal solutions (the formal reduction), of systems of the following form: 
\begin{equation}
\label{paramh}
\varepsilon^h \; \partial F\;=  A(x, \varepsilon) F, 
\end{equation} 
where the entries of the matrix $A(x, \varepsilon)$ are formal power series in $\varepsilon$, whose coefficients are formal power series in $x$; and $h$ is an integer which we call the \textit{$\varepsilon$-rank} of the system.

In~\cite[p. 74]{key221}, Iwano summarized as follows the problems needed to be resolved to obtain the complete knowledge about the asymptotic behavior of the solutions of a singularly-perturbed linear differential system:
\begin{description}
\item[$(\rm P1)$] \textit{Divide a domain $[D]$ in $(x, \varepsilon)$-space into a finite number of sub-domains so that the solution behaves quite differently as $\varepsilon$ tends to zero in each of these sub-domains;}
\item[$(\rm P2)$] \textit{Find out a complete set of asymptotic expressions of independent solutions in each of these sub-domains};
\item[$(\rm P3)$] \textit{Determine the so-called connection formula; i.e.\ a relation connecting two different complete sets of the asymptotic expressions obtained in $(2)$}.
\end{description}
In order to resolve $(\rm P1)$ and $(\rm P2)$, Iwano proceeded by associating a convex polygon to first-order systems in analogy with the scalar case, after imposing a precise triangular structure on $A(x, \varepsilon)$. $(\rm P3)$ remains generally unresolved and it is out of the scope of this article as well (see, e.g.\, ~\cite{key203,key201,key60} and references therein). The hypotheses on $A(x,\varepsilon)$ were eventually relaxed for special types of systems in a series of papers by Iwano, Sibuya, Wasow, Nakano, Nishimoto, and others (see references in~\cite{key59}). A famous and prevalent method among scientists is the \textsc{Wkb} method (see, e.g.\ \cite{key220}), which, roughly speaking, reduces the system at hand into one whose asymptotic behavior is known in the literature. However, as pointed out in~\cite{key216}, in view of the great variety and complexity of systems given by~\eqref{paramh}, it is rarely expected that a system can be reduced into an already investigated simpler form.

The methods proposed in the literature for the symbolic resolution of such systems are not yet fully algorithmic, exclude turning points (see~\cite{key61} and references therein), treat systems of dimension two only, rely on Arnold-Wasow form~\cite{key10}, and/or impose further restrictions on the structure of the input matrix. Moreover, there exists no computer algebra package dedicated to the symbolic resolution of neither system~\eqref{paramh} nor the scalar $n^{th}$-order scalar equation. The widely-used softwares \textsc{Maple} and \textsc{Mathematica} content themselves to the computation of outer formal solutions for scalar equations, and so does~\cite{key77}. 

In this article, we attempt to resolve $(\rm P2)$ algorithmically, that is give an algorithm to construct a fundamental matrix of formal solutions of an input system. We also attempt to give an insight into $(\rm P1)$ within some examples. 
  
As illustrated by the introductory example, the formal reduction of a singularly-perturbed system given by~\eqref{paramh} leads inevitably to the consideration of more general systems. These considerations will be discussed in Section~\ref{ringofcoefficients}. In Section~\ref{prelim}, we give necessary preliminaries including results on $n^{th}$-order scalar equations, the base cases, and  the well-known Splitting lemma which splits the system into subsystems of lower dimensions, whenever possible.  We discuss our contribution in the remaining sections, as we deviate from the classical reduction to give algorithms which act on the system directly without resorting to an equivalent scalar equation or to the Arnold form: First we refine in Sections~\ref{turnpt} and~\ref{moser} the algorithms which we proposed in~\cite{key102} for the resolution of turning points and reduction of the $\varepsilon$-rank $h$ to its minimal integer value. Then in Section~\ref{exporder}, we associate an $\varepsilon$-Newton polygon to a given system and show its role in retrieving formal invariants. And finally, we make use in Section~\ref{algokatzmain} of these notions introduced to propose an overall formal reduction algorithm. 

As we have already indicated, the literature on this problem is vast. In the hope of keeping this article concise but self-contained, we restrict our presentation to constructive results which contribute directly to our process of formal reduction. The omitted proofs are collected in the appendix along with illustrating examples.

In this paper, we content ourselves to the process of formal reduction. Any reference to the asymptotic interpretation of formal solutions will be dropped in the sequel. One may consult in this direction the method of matched asymptotic expansions (see, e.g.~\cite{key3623,key60} and references therein) or composite asymptotic expansions (see, e.g. ~\cite{key203,key201} and references therein). 

In the sequel, we adopt $\mathbb{C}$ as the base field for the simplicity of the presentation. However, any other commutative field $\mathcal{F}$ of characteristic zero ($\mathbb{Q} \subseteq \mathcal{F} \; \subseteq \; \mathbb{C}$) can be considered instead. The algorithms presented herein require naturally algebraic extensions of the base field. However, such extensions can be restricted as described in~\cite[Section 7.2]{key24}. 

Since we are leading a local investigation, we assume that the input system has at most one singular point, otherwise the region of study can be shrunk. Moreover, we place this singular point at the origin, otherwise a translation in the independent variable can be performed ($x \rightarrow 1/x$ for a singularity at $\infty$). 
\section*{Notations}
\begin{itemize}
\item $\mathbb{C}$ is the field of complex numbers, $\mathbb{Q}$ is the field of rational numbers, and $\mathbb{Q}^- = \{ q \in \mathbb{Q}\;|\; q \leq 0\}$.
\item We denote the dimension of the systems under consideration by $n$. We use $r$ for the algebraic rank of specified matrices (leading matrix coefficients) and $h$ for the $\varepsilon$-rank. 
\item We use upper case letters for algebraic structures, matrices, and vectors; the lower case letters $i, j, k, l$ for indices; and some lower case letters (greek and latin) locally in proofs and sections, such as $u$, $v$, $\sigma$, $\rho$, $\varrho$, etc...
\item We use $x, t, \tau$ as independent variables and $\partial$ to signify a derivation. Since $x$ dominates this article,  we drop it in the derivation. We set $\partial = \frac{d}{dx}$ and $\partial_\tau= \frac{d}{d\tau}$.
\item  We use $F$ for the unknown $n$-dimensional column vector (or $n \times n$-matrix); we use $G$, $U$, $W$, and $Z$ in this context as well. For scalar equations, we use $f$ and $g$.
\item $I_{d \times n}$ (resp. $I_{n}$) stands for the identity matrix of dimensions $d \times n$ (resp. $n \times n$); and $O_{d \times n}$ (resp. $O_{n}$) stands for the zero matrix of dimensions $d \times n$ (resp. $n \times n$). If confusion is unlikely to arise, we simply use $0$ to denote $O_{d \times n}$ or $O_{n}$.
\item We say that $M \in \mathcal{M}_n(\rm R)$ whenever the matrix $M$ is a square matrix of size $n$ whose entries lie in a ring $\rm R$.
\item $GL_n(\rm R)$ is the general linear group of degree $n$ over $\rm R$ (the set of $n \times n$ invertible matrices together with the operation of matrix multiplication). 
\item $\mathbb{C}[[x]]$ is the ring of formal power series in $x$ whose coefficients lie in $\mathbb{C}$; and $\mathbb{C}((x))$ is its fraction field. The usual valuation of an element $g(x)$ of $\mathbb{C}((x))$, i.e.\ the order of $g$ in $x$, is denoted by $val_{x} (g)$, with $val_{x} (0) = + \infty$. The ring $\mathbb{C}[[x]]$ (resp. the field $\mathbb{C}((x))$) is endowed with a differential ring (resp. field) structure by considering the derivation $\partial$. 
\item $\mathbb{C}[[\varepsilon]]$ is the ring of formal power series in $\varepsilon$ whose coefficients lie in $\mathbb{C}$; and $\mathbb{C}((\varepsilon))$ is its fraction field. $\mathbb{C}[[x]][[\varepsilon]]$ (resp. $\mathbb{C}[[x]]((\varepsilon))$) is the ring of formal (resp. meromorphic) power series in $\varepsilon$ whose coefficients lie in $\mathbb{C}[[x]]$. Similarly, $\mathbb{C}((x))[[\varepsilon]]$ (resp. $\mathbb{C}((x))((\varepsilon))$) is the ring (resp. field) of formal (resp. meromorphic) power series in $\varepsilon$ whose coefficients lie in $\mathbb{C}((x))$. 
\item Any element $f \in \mathbb{C}((x))((\varepsilon))$ can be written as $f=\sum_{k \in \mathbb{Z}}^\infty f_k(x) \varepsilon^k$, where $f_k(x) \in \mathbb{C}((x))$ for all $k \in \mathbb{Z}$. The derivation $\partial$ extends naturally to $\mathbb{C}((x))((\varepsilon))$ by the formula $\partial f= \sum_{k \in \mathbb{Z}}^\infty (\partial f_k(x)) \varepsilon^k$. Moreover, the map $val_{\varepsilon} :\; \mathbb{C}((x))((\varepsilon))\; \rightarrow \; \mathbb{Z} \cup \infty$ defines a valuation over $\mathbb{C}((x))((\varepsilon))$, satisfying the following properties for all $f(x, \varepsilon),\; g(x, \varepsilon)$ in $\mathbb{C}((x))((\varepsilon))$:
\begin{itemize}
\item $val_{\varepsilon} \;(f)\; = \infty$ if and only if $f =0$; 
\item $val_{\varepsilon} \;(f g)\; = val_{\varepsilon} \;(f) + val_{\varepsilon} \;(g) $;
\item $val_{\varepsilon} \;(f + g)\; \geq min\;(val_{\varepsilon} \;(f),\;val_{\varepsilon} \;(g)) $. 
\end{itemize}
\item We give the blocks of a matrix $M$ with upper indices, e.g.
\[
M \;=\;\begin{bmatrix}
    M^{11} & M^{12}\\
    M^{21} & M^{22}
  \end{bmatrix},
\]
and the size of the different blocks is dropped unless it is not clear from the context.
\end{itemize}

\section{Ring of coefficients}
\label{ringofcoefficients}
In this first part, we investigate a suitable ring of coefficients to treat first-order singularly-perturbed linear differential systems and $n^{th}$-order scalar singularly-perturbed differential equations, in a neighborhood of a turning point (see also~\cite[Chapter 2]{key59} and references therein). For $(U_0, V_0) \in \mathbb{R}^2$, we put
$$P(U_0, V_0) = \{ (U,V) \in \mathbb{R}^2 \;|\; U \geq U_0 \;\text{and}\; V \geq V_0 \}.$$
 Let $f  = \sum_{k \in \mathbb{Z}}^\infty f_k(x) \varepsilon^k \in \mathbb{C}((x))((\varepsilon))$ and let $P_{f}$ be the union of $P(k, val_x(f_k)), k \in \mathbb{Z}$. Then the Newton polygon of $f$, denoted by $\mathcal{N}_f$, is the boundary of the convex hull in $\mathbb{R}^2$ of the set $P_{f}$. We now consider, for $\sigma, p \in \mathbb{Q}$ with $\sigma \leq 0$, the half-plane $$H_{\sigma,p} = \{ (U,V) \in \mathbb{R}^2 \;|\; V \geq \sigma U + p \} . $$
One can verify that
 $$ \rm \mathcal{K} = \{ \;f \in \mathbb{C}((x))((\varepsilon)) \;|\; \mathcal{N}_f \subset H_{\sigma, p}\; \text{for some}\; \sigma \in \mathbb{Q}^-, p \in \mathbb{Q} \;\} $$
 is a differential field endowed with the derivation $\partial=d/dx$ and it is the field of fractions of the ring $\rm{\mathcal{R}} = \rm{\mathcal{K}}  \cap \mathbb{C}((x))[[\varepsilon]]$. Moreover, the elements of $\rm \mathcal{K}$ can be characterized geometrically as follows:
$$ f = \sum_{k \in \mathbb{Z}} f_k(x) \varepsilon^k \in \rm{\mathcal{K}}  \iff \exists \; \sigma \in \mathbb{Q}^-, p \in \mathbb{Q} \;|\; \forall \; k, val_x(f_k) \geq \sigma k + p .$$
In fact, given a non-zero element $f$ of $\rm{\mathcal{K}} $, we set $\nu_f := val_\varepsilon (f)$. Among the half-planes which contain $\mathcal{N}_f$, we consider the half-planes which are bounded below by a straight line passing through the point $(\nu_f, val_x(f_{\nu_f}))$, i.e. the half-planes $H_{\sigma, p_f}$ with $p_f := val_x(f_{\nu_f}) - \sigma \nu_f$. Then, among the $H_{\sigma, p_f}$'s, we pick the half-plane of maximal slope, which we denote by $\sigma_f$. Evidently, $H_{\sigma_f, p_f}$  is determined by a straight line which passes through the point $(\nu_f, val_x(f_{\nu_f}))$ and point(s) $(k, val_x(f_k))$ for at least one $k > \nu$. We denote by $(\mathcal{S}_f)$ this straight line whose equation is given by $V = \sigma_f U + p_f$, and we say that it is associated to $f$. For $f=0$, we set $\sigma_f=0$ and $p_f=0$.
\begin{example}
\begin{figure}[H]
\label{points}
\centering
\includegraphics[width=60mm]{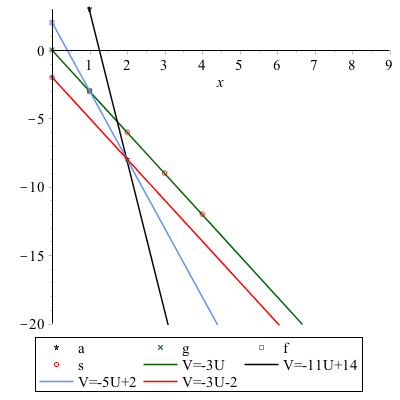}
\caption{Geometric interpretation of Example~\ref{sigmadraw}: $(\mathcal{S}_a)$, $(\mathcal{S}_g)$, $(\mathcal{S}_f)$, $(\mathcal{S}_s)$, are traced in black, green, blue, and red respectively.}
\end{figure}
\label{sigmadraw}
Given a nonzero element $f$ of $\;\mathbb{C}((x))((\varepsilon))$. If we can construct a straight line of finite slope which passes through $(\nu_f, val_x(f_{\nu_f}))$, and stays below the points $(k, val_x(f_{k}))$ for all $k > \nu$, then $f \in \rm{\mathcal{K}}$.  We illustrate the following examples in Figure~\ref{points}:
\begin{itemize}
\item Let $a= x^{3} \varepsilon + x^{-8} \varepsilon^2 = x^{3} \varepsilon (1 + x^{-11} \varepsilon)$ then $\sigma_a =-11$.
\item Let $g= \sum_{k=0}^\infty x^{-3k} \varepsilon^k$ then $\sigma_g =-3$.
\item Let $f= x^{2} + \sum_{k=1}^\infty x^{-3k} \varepsilon^k$  then $\sigma_f =-5$. 
\item Let $s= x^{-2} + \sum_{k=1}^\infty x^{-3k} \varepsilon^k$ then  $\sigma_s = -3$. 
\end{itemize}
\end{example} 

The following two lemmas give further insight into $\rm{\mathcal{K}}$:
\begin{lemma}
If $f(x, \varepsilon) \in \rm{\mathcal{K}} $ is non-zero, then there exists $e \in \mathbb{Z}^-$ such that $f(x, x^{e} \varepsilon) \in \mathbb{C}[[x]]((\varepsilon))$ .
\end{lemma}
\begin{proof}
Let $f = \sum_{k=\nu_f}^\infty f_k(x) \varepsilon^k$. We set $\mu_k= val_x(f_{k})$ and $f_k = x^{\mu_k} \tilde{f}_k$. Then for any $e \in \mathbb{Z}^-$ with $e \leq \sigma$ we have:
\begin{eqnarray*}f(x, \varepsilon) &=& \varepsilon^{\nu_f} (f_{\nu_f} +f_{\nu_f +1}\varepsilon + f_{\nu_f +2}\varepsilon^2 + f_{\nu_f +3}\varepsilon^3 +\dots ) \\ &=& \varepsilon^{\nu_f} x^{\mu_{\nu_f}}({\tilde{f}}_{\nu_f} + x^{\mu_{\nu_f+1} -\mu_{\nu_f}}{\tilde{f}}_{\nu_f +1}\varepsilon + x^{\mu_{\nu_f+2} -\mu_{\nu_f}} {\tilde{f}}_{\nu_f +2}\varepsilon^2 + x^{\mu_{\nu_f+3} -\mu_{\nu_f}}{\tilde{f}}_{\nu_f +3}\varepsilon^3 +\dots ) \\ &=&\varepsilon^{\nu_f} x^{\mu_{\nu_f}}({\tilde{f}}_{\nu_f} + x^{\mu_{\nu_f+1} -\mu_{\nu_f}-e}{\tilde{f}}_{\nu_f +1}(x^{e}   \varepsilon) + x^{\mu_{\nu_f+2} -\mu_{\nu_f} -2{e} } {\tilde{f}}_{\nu_f +2}(x^{e}  \varepsilon)^2 + x^{\mu_{\nu_f+3} -\mu_{\nu_f} -3{e} }{\tilde{f}}_{\nu_f +3}(x^{e}   \varepsilon)^3 +\dots ) .\end{eqnarray*}
But, $\mu_k\geq \sigma_f (k-\nu_f) + \mu_{\nu_f}$ for all $k \geq \nu_f$. Hence, $\mu_k\geq e (k-\nu_f) + \mu_{\nu_f}$ for all $k \geq \nu_f$. We can thus choose $e$ to be the largest integer which is less than or equal to $\sigma_f$.
\end{proof}
\begin{lemma}
The field of constants of the differential field $(\rm{\mathcal{K}}, \partial)$ is $\mathbb{C}((\varepsilon))$.
\end{lemma}
\begin{proof}
Let $f \in \rm{\mathcal{K}} $ and suppose that $\partial{f} =0$. Then, $\sum_{k=\nu_f}^{\infty} (\partial f_k(x)) \varepsilon^k =0$ which yields  $\partial f_k(x) =0$ for all $k \geq \nu_f$. Thus, for every $k \geq \nu_f$,  $f_k(x)= c_k$ for some $c_k \in \mathbb{C}$, which yields $f= \sum_{k=\nu_f}^{\infty} c_k \varepsilon^k$. 
\end{proof}
Henceforth, algorithmically, we can use the following handy representation of any element $f \in \rm{\mathcal{K}}$:
\begin{equation} \label{presentationxi} f(x, \varepsilon)=\sum_{k =\nu_f}^\infty f_k(x) \varepsilon^k = \sum_{k =\nu_f}^\infty x^{-k \sigma_f }  f_k(x) (x^{\sigma_f} \varepsilon)^k = x^{p_f} \sum_{k =\nu_f}^\infty x^{-k \sigma_f -p_f}  f_k(x) \xi_f^k, \quad \text{with} \quad \xi_f= x^{\sigma_f} \varepsilon .\end{equation}
Evidently, $val_x(x^{-\sigma_f k -p_f}  f_k(x)) \geq 0$ for all $k \geq \nu_f$. Moreover, under the notation $\xi_f$, we have:
$$ \partial (\sum_{k =\nu_f} f_k \xi_f^k) =  \sum_{k =\nu_f} (\partial (f_k) + \frac{\sigma_f k}{x} f_k) \xi_f^{k} .$$ 
\subsection{Input systems and input equations}
Motivated by the introductory example, and with the help of the above notations, we will investigate the following system rather than~\eqref{paramh}:
\begin{equation} \label{paramK} \partial F = A(x, \varepsilon) F \quad \text{where} \quad A \in \mathcal{M}_n(\rm{\mathcal{K}}) . \end{equation}
The valuation $val_\varepsilon$ (and the usual valuation $val_x$ of $\mathbb{C}((x))$) extends naturally to matrices $A = (a_{i,j}) \in \mathcal{M}_n(\rm{\mathcal{K}}) $ by $val_{\varepsilon}(A) = \operatorname{min}_{i,j} val_\varepsilon (a_{i,j})$. We can then define the Newton polygon of $A  \in \mathcal{M}_n(\rm{\mathcal{K}})$, $\nu_A \in \mathbb{Z}$, $p_A \in \mathbb{Q}$, and $\sigma_A \in \mathbb{Q}^-$, in an analogous manner to those of $f \in \rm{\mathcal{K}}$. Hence, we can express $A(x, \varepsilon)$ as follows:
\begin{equation} \label{psigmaa} 
A = \sum_{k=\nu_A}^{\infty} A_k \varepsilon_A^k, \; A_{\nu_A} \neq 0,\; p_A = val_x(A_{\nu_A}) - \sigma_A \nu_A,  \; \text{and} \;
val_x (A_k) \geq \sigma_A k + p_A  \; \forall \; k \geq \nu_A. 
\end{equation}
Hence by setting $\xi_A = x^{\sigma_A} \varepsilon$, we can write in analogy with \eqref{presentationxi}:
\begin{eqnarray*}  A(x,\varepsilon) &=& x^{p_A}\sum_{k = \nu_A}^{\infty} x^{-k \sigma_A - p_A} A_k(x) \; \xi_A^k \\ &=& \xi_A^{\nu_A} x^{p_A} \; \sum_{k =0}^{\infty} x^{-k \sigma_A - val_x(A_{\nu_A})} A_{\nu_A + k}(x) \;  \xi_A^k \quad \text{where} \quad {x^{-val_x(A_{\nu_A})}A_{\nu_A}}|_{x=0} \neq 0  . \end{eqnarray*}
Thus, we can rewrite system~\eqref{paramK} as follows:
$$ \xi_A^{-\nu_A} x^{-p_A} \partial F = (\sum_{k=0}^{\infty} x^{-k \sigma_A - val_x(A_{\nu_A})} A_{\nu_A + k}(x) \;  \xi_A^k) \; F , \quad \text{where} \quad  x^{-k \sigma_A -val_x(A_{\nu_A})} A_{\nu_A + k} \in \mathcal{M}_n(\mathbb{C}[[x]]) .$$
Similarly, if we consider a $n^{th}$-order differential  equation whose coefficients are elements of $\rm{\mathcal{K}}$: 
\begin{equation} \label{diff_scalar} a_{n}(x, \varepsilon)  \partial^n\; f \; +\;  a_{n-1}(x, \varepsilon) \; \partial^{n-1}\; f \; +\; \dots \;+\; a_{1}(x, \varepsilon) \; \partial f \;+\;  a_{0}(x, \varepsilon) \; f \;=\;0, \end{equation}
then setting $F = [f, \partial f, \partial^2 f, \dots, \partial^{n-1} f]^T$, we can express~\eqref{diff_scalar} equivalently as a system $\partial F = A(x,\varepsilon) F$ where $A(x,\varepsilon) \in \mathcal{M}_n(\rm{\mathcal{K}})$ is a companion matrix, and we define 
$\nu_a := \nu_A$, $p_a := p_A$, $\sigma_a := \sigma_A$, and $\xi_a := \xi_A$.
\subsubsection{Notations}
Thus, in this paper, we treat system~\eqref{paramK} in the following form: 
\begin{equation}
\label{param}
 [A_{\sigma_A}] \quad \quad \quad \xi_A^{h} \; x^{p_A} \; \partial F\; = \; A(x, \xi_A) F ,\quad  A(x, \xi_A) = \sum_{k=0}^\infty A_k(x) \xi_A^k \in \mathcal{M}_n(\rm{\mathcal{R}}) , \quad A_0(x=0) \neq 0, 
\end{equation} 
where  $\xi_A=x^{\sigma_A} \varepsilon$, $\sigma_A \in \mathbb{Q}^-$, $p_A \in \mathbb{Q}$, and $h \in \mathbb{Z}$. Under the definition of $\xi_A$, $A(x, \xi_A) = \sum_{k=0} A_k(x) \xi_A^k$ with $A_k(x) \in \mathcal{M}_n(\mathbb{C}[[x]])$ for all $k \geq 0$. We refer to $A_0(x)$ as the leading coefficient matrix and to $A_{0,0}:= A_0(x=0)$ as the leading constant matrix.  

We also consider $n^{th}$-order scalar linear differential equations of the form:
\begin{equation} \label{scalarparammm} [a_{\sigma_a}] \quad \quad \quad \partial^n f \; +\;  a_{n-1}(x, \xi_a) \partial^{n-1} f \; +\; \dots \;+\; a_{1}(x, \xi_a) \partial f \;+\;  a_{0}(x, \xi_a) f \;=\;0    ,\end{equation}
where $a_i(x, \xi_a) = \sum_{k\in \mathbb{Z}}^\infty a_{i,k}(x)\; \varepsilon^k \in \rm{\mathcal{K}}$ for all $ i \in \{0, \dots, n\}$, with $a_{i,k}(x) \in \mathbb{C}[[x]]$ for all $k \geq 0$.

We remark that a system $[A_{\sigma_A}]$ given by~\eqref{param} can be also expressed as a scalar equation of the form~\eqref{scalarparammm}. The theoretical possibility stems from the work of~\cite{key213} which discusses cyclic vectors. Moreover, an algorithm to compute a companion block diagonal form of $[A_{\sigma_A}]$ (and hence equivalent scalar form) can be obtained by generalizing the work in~\cite{key17} developed for unperturbed systems (see~\ref{appcompanion}). 

In the sequel, for the clarity of the presentation, the index $A$ (resp. $a$) will be dropped from $\nu_A$, $p_A$, $\sigma_A$, and $\xi_A$ (resp. $\nu_a$, $p_a$, $\sigma_a$, and $\xi_a$) whenever confusion is unlikely to arise.
\subsection{Equivalent systems}
\label{urgent}
Consider system $[A_{\sigma_A}]$ given by~\eqref{param}. Let $T \in GL_n(\rm \mathcal{K})$ then the transformation $F = T G$ (also called \textit{gauge transformation}) yields a system
\begin{equation} \label{gaugep}  [\tilde{A}_{\sigma_{\tilde{A}}}] \quad \quad \quad {\xi_{\tilde{A}}}^{\tilde{h}} x^{\tilde{p}}  \partial G\;=\; \tilde{A}(x, {\xi_{\tilde{A}}}) G, \quad \text{where} \quad \tilde{A} \in \mathcal{M}_n(\rm{\mathcal{R}}).\end{equation}
for some $\sigma_{\tilde{A}} \in \mathbb{Q}^-$, $\tilde{p} \in \mathbb{Q}$, and $\tilde{h} \in \mathbb{Z}$.  
We say that systems $[A_{\sigma_A}]$ and $T [A_{\sigma_A}]:= [\tilde{A}_{\sigma_{\tilde{A}}}] $ are \textit{equivalent}. Examples of such transformations and their applications are the transformations $T_1, T_2, L_1, L_2$ applied within the introductory example. 

In the sequel, we seek at each step a transformation which yields for an input system  $[A_{\sigma_A}]$, an equivalent system $[\tilde{A}_{\sigma_{\tilde{A}}}]$, so that either $\tilde{h} < h$ or $[\tilde{A}_{\sigma_{\tilde{A}}}]$ can be decoupled into systems of lower dimensions. Using a recursive process which employs this approach, we aim to break down any input system $[A_{\sigma_A}]$ into system(s) which can be treated with known methods (see basic cases of Section~\ref{basiccases}).

Evidently, $\sigma_{\tilde{A}}$ might differ from $\sigma_A$ and hence $\xi_A$ will be updated to $\xi_{\tilde{A}}$ after the application of $T$ (see the introductory example), which brings us to the next subsection.  

\subsection{Restraining index and $(\rm P1)$}
\label{resindex}

As we are lead naturally to the study of system  $[A_{\sigma_A}]$ with $\sigma_A \in \mathbb{Q}^-$ within the study of system~\eqref{param}, we will need to understand the inevitable growth of the order of the poles in $x$ within the reduction. For instance, in the introductory example, we start with an input system for which $\sigma_A=0$ and arrive at an decoupled system with $\sigma_{U} = -3$. 

As described in~\cite[Chapter 2]{key59} and as we have observed, these poles grow at worst linearly with $k$. From here initially stems our main motivation for our choice of the ring of coefficients and the representation in terms of $\xi_A$: the orders of the poles in $x$ which might be introduced in an input system within the process of reduction, are stored in $p$ and $\sigma$, which allows investigating them at any step within the reduction process. We can thus talk about a restraining index $\rho=-1/ \sigma$. This is the first step towards determining the inner and intermediate regions. In fact, the complicated behavior anticipated in the neighborhood of turning points can be investigated with the help of a sequence of positive rational numbers:
\begin{equation}
\label{radii}
[\rho] \quad \quad \quad \quad 0 =\;{\rho}_0 < \rho_1 < \rho_2 < \dots < \rho_m . 
\end{equation} 
With this sequence, the domain $|x| \leq x_0$ of $[D]$ can be divided into a finite number of sub-domains in each of which the solution of behaves quite differently (see~\cite[Intro, pp 2]{key62} and ~\cite[Chapter 1]{key221}). In case $[\rho]$ is known, one can apply adequate stretching transformations to the original input system, i.e. a change of the independent variable of the form $\tau = x \varepsilon^{-\rho_i},\; i \in \{1, \dots, m\}$. Then, using our proposed algorithm, we can construct inner and intermediate solutions. In this paper, we do not investigate $(\rm P1)$, i.e. we do not compute $[\rho]$. However, we motivate a plausible approach in some examples (see~\ref{appe}).

\section{Preliminaries}
\label{prelim}
  A detailed discussion of the scalar case ($n=1$) and regularly-perturbed systems ($h \leq 0$), is given in~\cite[Chapter 3, Section 1, p. 52 -56]{key61} for $\sigma_A=p_A=0$. In this section, we give a brief discussion which discards this restriction on $\sigma_A$ and $p_A$.  
\subsection{The base cases}
\label{basiccases}
\subsubsection{The case $ h \leq 0$}
\label{hzero}
If $h \leq 0$ then the solution of system $[A_{\sigma_A}]$ can be sought, up to any order $\mu$, upon presenting the solution as a power series in $\xi_A = x^{\sigma_A} \varepsilon$, i.e.\ 
$$ F \;=\; \sum_{k=0}^{\infty} \;F_k (x) \; \xi_A^{-h+k}.$$ 
The latter can then be inserted into   $[A_{\sigma_A}] \; \xi^h x^p \partial F = A F \;$ and the like powers of $\xi_A$ equated. This reduces the symbolic resolution to solving successively a set of inhomogeneous linear singular differential systems in $x$ solely:
\begin{eqnarray*}
 x^p\; \partial F_0\; &=&\;  ( A_0(x) \;+\; h \; \sigma_A\; x^{p-1} I) \; F_0 \\
 x^p\; \partial F_1 \;&=& \;( A_0(x)\; +\; (h-1) \;\sigma_A \; x^{p-1} I) \; F_1\; +\; A_1(x)\; F_0  \\
 \vdots \\
  x^p \;\; \partial F_{\mu} \;&=&\; ( A_0(x) \;+\; (h-\mu) \; \sigma_A \; x^{p-1} I) \; F_{\mu} + A_1(x) \; F_{\mu-1} \;+\; \dots \;+\;A_{\mu}(x)\; F_0  .\\
\end{eqnarray*}
For any $k \geq 0$, the solution of the homogeneous system singular in $x$ can be sought via the \textsc{Maple} package \textsc{Isolde}~\cite{key27} or the \textsc{Mathemagix} package \textsc{Lindalg}~\cite{key427} (which are both based on~\cite{key24}). Afterwards, the solution of the inhomogeneous system can be obtained by the method of variation of constants (see, e.g. \cite[Theorem 3, p. 11]{key6}). We remark that a transformation $F_{k} =  x^{-k \sigma_A} G $ reduces the inhomogeneous system in the dependent vector $F_{k}$ to:
$$x^p \; \partial G \;= ( A_0(x) \;+ \;h \; \sigma_A \; x^{p-1} I_n) \; G(x) \; , $$
which coincides with that of $F_0(x)$. 
 Thus we have:
$$F_{k}(x) \;= \;x^{-{k} \sigma_A} \; F_0(x) \;\int \; F_0^{-1}(x)\; x^{{k} \sigma_A} \; [ A_1(x)\; F_{{k}-1}(x)\; + \;\dots \;+ \;A_{k}(x)\; F_0(x)] \; dx .$$

\subsubsection{The scalar case}
\label{macutanpolygon}
The $n^{th}$-order singularly-perturbed scalar differential equations are treated in~\cite{key77} upon introducing an analog to the Newton polygon and polynomials (see e.g.~\cite{key105} and references therein for the unperturbed counterparts of these equations): the $\varepsilon$-polygon and $\varepsilon$-polynomials. In this subsection, we adapt this treatment to the more general equation $[a_{\sigma_a}]$ given by~\eqref{scalarparammm}:
$$[a_{\sigma_a}] \quad \quad \quad  \partial^n f \; +\;  a_{n-1}(x, \xi_a) \partial^{n-1} f \; +\; \dots \;+\; a_{1}(x, \xi_a) \partial f \;+\;  a_{0}(x, \xi_a) f \;=\;0.$$

We first give an analog to ~\cite[Definition 2.1]{key77} and then generalize~\cite[Proposition 4.1]{key77}, whose proof is easily adaptable from~\cite{key77} (see~\ref{appb}).
\begin{definition} Consider the scalar equation $[a_{\sigma_a}]$ given by~\eqref{scalarparammm}. For $(u_0, v_0) \in \mathbb{R}^2$, we put 
$$ P(U_0, V_0) = \{ (U,V) \in \mathbb{R}^2 | \; U \leq U_0 \; \text{and} \; V \geq V_0 \} . $$
Set $\nu_{i} = val_\varepsilon(a_i)$ for $i \in \{0, \dots, n\}$ and let $P_{a}$ be the union of $P(i, \nu_{i})$. Then the $\varepsilon$-polygon of $[a_{\sigma_a}]$, denoted by $\mathcal{N}_{\varepsilon}(a)$, is the intersection of $\mathbb{R}^+ \times \mathbb{R}$ with the convex hull in $\mathbb{R}^2$ of the set $P_{a}$.
We denote the slopes of the edges of $\mathcal{N}_{\varepsilon}(a)$ by $\{ e_1, \dots, e_\ell \}$. These slopes are non-negative rational numbers.  And, for every $ j \in \{1, \dots, \ell \}$, consider the algebraic equation given by
$$( E_j ) \quad  \;  \sum_{k=0}^{\ell}  {a}_{i_k, \;val_{\varepsilon}(a_{i_k})}(x) \; X^{(i_k - i_0)} = 0,$$
where $0 \leq i_0 < i_1 < \dots < i_\ell = n$ denote the integers $i$ for which $(i, \nu_i)$ lie on the edge of slope $e_j$ of the $\varepsilon$-polygon, and ${a}_{i, \nu_i}(x)= \xi_a^{-\nu_i}\; a_i(x, \xi_a) {|_{\xi_a=0}}$. We say that $(E_j)$ is the $\varepsilon$-polynomial associated to the slope $e_j$.
\end{definition}
\begin{proposition}
\label{newtoneqnxi}
Consider a nonzero $[a_{\sigma_a}]$ given by~\eqref{scalarparammm} and its $\varepsilon$-polygon $\mathcal{N}_{\varepsilon}(a)$ of slopes $\{ e_1, \dots, e_\ell \}$.  Let $$f(x, \xi_a) = exp (\int q(x, \xi_a) dx) \quad  \text{with} \quad q(x, \xi_a) \neq 0 \in \bigcup_{s \in \mathbb{N}^{*}} \overline{\mathbb{C}((x))} ((\xi_a^{1/s})) . $$ 
If $f(x,\xi_a)$ satisfies $[a_{\sigma_a}]$ formally then $$q(x, \xi_a) = \frac{1}{\xi_a^{e_j}} (X(x) + O(\xi_a))$$ for some $j \in \{1, \dots, \ell \}$, and $X(x)$ is one of the non zero roots of the $\varepsilon$-polynomial $(E_j)$ associated to $e_j$. 
\end{proposition}
 The full expansion of $q(x, \xi_a)$ can be obtained with successive substitutions of the form $f = exp (\int \frac{X(x)}{\xi_a^{e_j}} dx)\; g$~\cite[Proposition 4.2]{key77}. 
\begin{figure}[H]
\label{figures3}
\centering
\includegraphics[width=60mm]{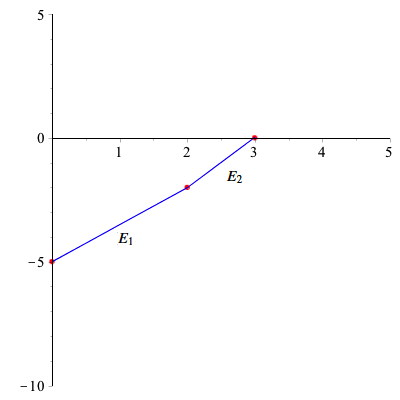}
\caption{The $\varepsilon$-polygon associated to \eqref{equiveq} in Example~\ref{exmwasow3}}
\end{figure}
\begin{example}
\label{exmwasow3}
Consider the following linear differential equation with $\sigma_a=0$:
\begin{equation} 
\label{equiveq}
[a_{\sigma_a}] \quad \quad \quad \quad \partial^3 f  \;-\; \frac{x}{\xi_a^2}\;\partial^2 f\; \;-\; \frac{1}{\xi_a^5} f \;=\; 0.
\end{equation}
We wish to study this equation according to its $\varepsilon$-polygon (see Figure~$3$). We have two slopes:
\begin{itemize}
\item The slope $e_1 = \frac{3}{2}$ for which $\; (E_1) \quad  x\;X^2 + 1 = 0$ whose nonzero solutions are $X = \pm \frac{i}{\sqrt{x}}$.  
\item The slope $e_2 = 2$ for which $\; (E_2) \quad X^3 - x X^2 = 0$ whose nonzero solution is $X = x$.
\end{itemize}
Thus, the leading terms (of the exponential part) of solutions of~\eqref{equiveq} are given by
$exp \;  (\int \frac{x}{\xi_a^2} dx) = exp \;  (\int \frac{x}{\varepsilon^2} dx)$ and $exp \; (\pm \int \frac{ i}{\sqrt{x} \xi_a^{3/2}} dx) = exp \; (\pm \int \frac{ i}{\sqrt{x} \varepsilon^{3/2}} dx)$.
\end{example}
One can also refer to~\cite[Section 5]{key77} for examples on the harmonic oscillator with damping and Orr-Sommerfeld equation. 

Based on the above, we give the following definition:
\begin{definition}
Consider the scalar equation $[a_{\sigma_a}]$ given by~\eqref{scalarparammm}. We call the largest slope of the $\varepsilon$-polygon $\mathcal{N}_{\varepsilon}(a)$, the $\varepsilon$-formal exponential order ($\varepsilon$-exponential order, in short) :
$$\omega_{\varepsilon} (a) = max \; \{ e_j , 1 \leq j \leq \ell \} .$$
\end{definition}
 Clearly, under the notations of this subsection, we have:
$$\omega_{\varepsilon} (a) \;=\; \maxi_{i=0}^{n-1} \;(0, - \frac{val_{\varepsilon}(a_i)}{n - i}).$$

\subsection{Consequences of the equivalence between a scalar equation and a system}
Due to their equivalence, the formal solutions of a first-order system $[A_{\sigma_A}]$ given by~\eqref{param} can be computed from an equivalent scalar equation $[a_{\sigma_a}]$ given by~\eqref{scalarparammm}, using the $\varepsilon$-polygon of $[a_{\sigma_a}]$. Moreover,~\eqref{radii} can be computed as well using the Iwano-Sibuya's polygon (which collects valuations with respect to both $x$ and $\varepsilon$~\cite{key62}). However such a treatment is unsatisfactory (see e.g.~\cite[Conclusion]{key77}) since it overlooks the information that we can derive from the system directly and demands an indirect treatment. Nevertheless, it plays a key role in the theoretical basis of the algorithm which we will develop in later parts. Most importantly, the equivalence between $[A_{\sigma_A}]$ and $[a_{\sigma_a}]$, and the invariance of the $\varepsilon$-exponential order under gauge transformations, allows the definition of the $\varepsilon$-exponential order of system $[A_{\sigma_A}]$ as follows:
\begin{definition}
Consider a system $[A_{\sigma_A}]$ given by~\eqref{param} and an equivalent  scalar equation $[a_{\sigma_a}]$ given by~\eqref{scalarparammm}. The $\varepsilon$-polynomials and $\varepsilon$-formal exponential order $\omega_{\varepsilon} (A)$ of $[A_{\sigma_A}]$ are those of $[a_{\sigma_a}]$. 
\end{definition}
It follows that, given a system $[A_{\sigma_A}]$, we have:
\begin{equation} \label{epsilonkappa} \omega_{\varepsilon} (A) = \max_{0 \leq i \leq n} \; (0,\; \frac{ - val_{\varepsilon}\; (a_i)}{n-i} ).\end{equation}
where $[a_{\sigma_a}]$ is a scalar equation equivalent to $[A_{\sigma_A}]$. 

In this paper, we give a direct treatment of the system which will lead eventually to the computation of $\omega_{\varepsilon} (A)$, the $\varepsilon$-polynomials, and consequently a basis of the space of formal solutions, without resorting to the equivalence of an input system to a scalar equation. Our treatment relies on the properties of the eigenvalues of the leading coefficient of the input system. We distinguish between three cases (distinct, unique, or zero eigenvalues), the first of which is classical and is recalled in the next subsection.

\subsection{$\varepsilon$-Block Diagonalization}
\label{blockparam}
Consider a system $[A_{\sigma_A}]$ given by~\eqref{param}. A classical tool in the perturbation theory is the so-called splitting  which separates off the existing distinct coalescence patterns. Whenever the leading constant matrix $A_{0,0} := A_0(x=0)$ admits at least two distinct eigenvalues, the system can be decoupled into subsystems of lower dimensions:
\begin{theorem}
\label{splitparam}
Consider system $[A_{\sigma_A}]$ given by~\eqref{param}
$$
 [A_{\sigma_A}] \quad \quad \quad \xi_A^{h} x^{p} \partial F\; = \; A(x, \xi_A) F ,\quad  A(x, \xi_A) \in \mathcal{M}_n(\rm{\mathcal{R}}) , \quad A_{0,0} \neq 0. $$ 
If $A_{0,0}  = \begin{bmatrix} A_{0,0}^{11} & O \\ O &  A_{0,0}^{22} \end{bmatrix}$ such that $A_{0,0}^{11}$ and $ A_{0,0}^{22}$ have no eigenvalues in common, then there exists a unique transformation $T(x, \xi_T) = \sum_{k=0}^\infty T_k(x) \xi_T^k \in GL_n(\rm{\mathcal{R}})$, given by 
$$T(x, \xi_T) = I + \sum_{k=0}^{\infty} \begin{bmatrix} O &  T^{12}_k(x)  \\ T^{21}_k(x)  & O \end{bmatrix} {\xi_T}^k ,$$
such that the transformation $F = T G$ gives $$
  [\tilde{A}_{\sigma_{\tilde{A}}}] \quad \quad \quad \xi_{\tilde{A}}^{h} x^p \; \partial G \;= \;\tilde{A}(x, \xi_{\tilde{A}}) \;G \;= \; \begin{bmatrix} \tilde{A}^{11} (x, \xi_{\tilde{A}}) & O \\ O &  \tilde{A}^{22}(x, \xi_{\tilde{A}}) \end{bmatrix} \; G \quad \text{where} \quad \tilde{A}_{0,0}= A_{0,0}. $$ Moreover $\sigma_T, \sigma_{\tilde{A}}\geq \sigma$ where $\sigma = \begin{cases} \sigma_A +\frac{p-1}{h} \quad \text{if} \; p<1 \\ \sigma_A \quad \text{otherwise} \end{cases}.$
\end{theorem}
Well-known proofs proceed in two steps: Block-diagonalizing $A_0(x)$, followed by successive block-diagonalization of coefficients (see, e.g., \cite[Theorem XII-4-1, p. 381]{key71}  and~\cite[Chapter 3, p. 56-59]{key61} for the particular case of $\sigma_A=p_A=0$, and \cite[Theorem 8.1, p. 70-81]{key59} for the general case). We give a constructive proof in~\ref{appb}. We remark that the specified form of $A_{0,0}$ in Theorem~\ref{splitparam} is non-restrictive. It suffices that $A_{0,0}$ has at least two distinct eigenvalues. Such a block-diagonal form can be then attained by Jordan constant transformation. In case the base field is not algebraically closed, the weaker transformation given in~\cite[Lemma A.1]{key26} can be applied instead, since it does not require any algebraic field extensions.

After splitting the input system, we proceed in our reduction for each of the decoupled subsystems in parallel. We can thus assume, without any loss of generality, that the leading constant matrix of the input system $[A_{\sigma_A}]$ has a unique eigenvalue $\gamma \in \mathbb{C}$. Thus, upon applying the so-called eigenvalue shifting, i.e.\ \begin{equation} \label{evshiftparam} F = E G =  G \; \exp (\int \xi_A^{-h} x^{-p} \gamma  dx), \end{equation} one can verify that the resulting system $E [A_{\sigma_A}]$ has a nilpotent leading constant matrix. Thus, it remains to discuss this case of nilpotency. This is the goal of the two following sections.

\section{Two-fold rank reduction}
\label{twofoldrankreduction}
Without loss of generality, we can now assume that system $[A_{\sigma_A}]$ given by~\eqref{param} is such that $A_{0,0}$ is nilpotent. At this stage of reduction, our approach diverges from the classical  indirect ones which require Arnold-Wasow form, a cyclic vector, or reduction to a companion form. Two cases arise:
\begin{itemize}
\item $A_{0,0}$ is nilpotent but $A_0(x)$ is not: In this case, the eigenvalues of $A_0(x)$ might coalesce in some neighborhood of $x=0$, which would cause the so-called turning point. In Subsection~\ref{turnpt}, we propose an algorithm to compute a transformation which reduces the input system to a system for which both $A_{0,0}$ and $A_0(x)$ are nilpotent. 
\item Both $A_{0,0}$ and $A_0(x)$ are nilpotent: In this case, we propose in Subsection~\ref{moser} an algorithm to compute a transformation which can reduce the $\varepsilon$-rank $h$ of system $[A_{\sigma_A}]$ to its minimal nonzero integer value. This minimality is an intrinsic property of the system and its solutions. Its computation paves the way for the retrieval of the $\varepsilon$-exponential order of the system, as we show later in Section~\ref{exporder}.
\end{itemize}
\subsection{Resolution of turning points}
\label{turnpt}
We illustrate the source of turning points in the following example:
\begin{example}\cite[p. 223]{key60}
\label{exmturnpt}
Let $a(x) , b(x)$ be holomorphic for $x$ in a region $\Omega  \in \mathbb{C}$ and consider
$$M(x)=\begin{bmatrix}  x & 1 & a(x) \\  0 & x & b(x) \\ 0 & 0 & 0   \end{bmatrix} .$$
The transition matrix $T(x)=\begin{bmatrix}  \alpha(x) & \beta(x) &  \gamma(x)(- a(x) x^{-1} + b(x) x^{-2})\\  0 & \alpha(x)  &  - \gamma(x) b(x) x^{-1} \\ 0 & 0 & \gamma(x) \end{bmatrix}$, with arbitrary scalar functions $\alpha (x) , \beta(x), \gamma(x)$ such that $\alpha(x) \neq 0$ and $\gamma(x) \neq 0$ for all $x \in \Omega$, yields for $x \neq 0$:
 $$ J(x)=T ^{-1} M(x) T(x)=\begin{bmatrix}  x & 1 &0 \\  0 & x & 0 \\ 0 & 0 & 0   \end{bmatrix} . $$ 
If $\alpha (x) , \beta(x), \gamma(x)$ are holomorphic, so is $T$, provided that zero is not a point of $\Omega$. If $0 \in \Omega$, then $T(x)$ has a pole at $x=0$ unless the following conditions are satisfied:
 $$ b(0)=0\;  \; \text{and} \;  \partial b (0) = a(0).$$
 Therefore, $M(x)$ is in general not holomorphically similar to $J(x)$ in regions that contain $x=0$. Moreover, if the second condition is not satisfied, two possibilities arise:
 \begin{itemize}
 \item if $b(0)=0$ then there exists a constant invertible matrix $S$ such that 
 $ S ^{-1} M(0) S=\begin{bmatrix}  0 & 1 &0 \\  0 &0 & 0 \\ 0 & 0 & 0   \end{bmatrix} . $ 
 Then $M(x)$ is pointwise similar to $J(x)$, even in regions containing $x=0$.
 \item if $b(0) \neq 0$ then the Jordan form of $M(0)$ is 
 $J(0) =\begin{bmatrix}  0 & 1 &0 \\  0 &0 & 1 \\ 0 & 0 & 0   \end{bmatrix} $
 so that $J(x)$ is not holomorphic at $x=0$. 
 \end{itemize}
\end{example}
 Hence, it might occur that the Jordan matrix is not holomorphic in some region or the Jordan matrix is itself holomorphic but, nevertheless, not holomorphically (although perhaps point-wise) similar to $M(x)$. Below is an example which exhibits a system whose leading coefficient matrix has such complications:
\begin{example}\cite[p. 57]{key59}
\label{exmwasow6}
Consider the following system with $\sigma_A = 0$, i.e. $\xi_A = \varepsilon$: 
$$ [A_{\sigma_A}] \quad \quad \xi_A^{2}  \partial F = \begin{bmatrix}  0 & 1& 0 \\ 0 & 0 & 1 \\  \xi_A & 0 & x \end{bmatrix} F \quad \text{where} \quad A_0(x)= \begin{bmatrix}  0 & 1& 0 \\ 0 & 0 & 1 \\0 & 0 & x \end{bmatrix} .$$
It follows from the form of $A_0(x)$ that the origin is possibly a turning point. We remark that this system is equivalent to~\eqref{equiveq}, as it is obtained from the latter by setting $F = [f, \xi_A^2 \partial f, \xi_A^4 \partial^2 f ]^T$.
\end{example}
 We recall the following proposition which we first gave in~\cite{key102}, and we refine its proof: 
\begin{proposition}~\cite[Section 3]{key102}
\label{nestedprop}
Consider system~\eqref{param} given by 
$$
[A_{\sigma_A}] \quad \quad \quad  \xi_A^{h} x^{p} \; \partial F\; = \;A(x, \xi_A) \;F\; =\;  \sum_{k=0}^{\infty} A_k(x) \xi_A^k  \; F
$$
with $\xi_A = x^{\sigma_A} \varepsilon$. Suppose that the leading constant matrix $A_{0,0}$ is nilpotent and $A_0(x)$ has at least one nonzero eigenvalue.  Then there exists an invertible polynomial transformation $T$ in a root of $x$ such that the transformation $F = T G$ and a re-adjustment of the independent variable $x$ results in a system  whose leading constant matrix has at least one nonzero constant eigenvalue. 
\end{proposition}
\begin{example}
Consider the matricial form of Weber's equation ($\sigma_A =0$, $p_A=0$) given by:
$$[A_{\sigma_A}] \quad \quad \quad \xi_A \frac{dF}{dx} = \begin{bmatrix} 0 & 1 \\ x^2 & 0\end{bmatrix} F,$$
with $\xi_A = x^{\sigma_A} \varepsilon$. Let $T(x) = \begin{bmatrix} 1& 0 \\ 0 & x \end{bmatrix}$ then $F = T G$ yields:
$$[B_{\sigma_A}] \quad \quad \quad \xi \frac{dG}{dx} = \begin{bmatrix} 0 & x \\ x & -\xi_A x^{-1}\end{bmatrix} G.$$
Upon factorizing $x$, we have:
$$ \xi_A x^{-1} \frac{dG}{dx} = \begin{bmatrix} 0 & 1 \\ 1 & -\xi_A x^{-2}\end{bmatrix} G.$$
We compute $\sigma_{\tilde{A}} = -2$ and thus set $\xi_{\tilde{A}} = x^{-2} \varepsilon$. The above system can be rewritten as:
$$[\tilde{A}_{\sigma_{\tilde{A}}}] \quad \quad \quad \quad \xi_{\tilde{A}} x \frac{dG}{dx} = \begin{bmatrix} 0 & 1 \\ 1 & -\xi_{\tilde{A}} \end{bmatrix} G .$$
Clearly, $\tilde{A}_{0,0}$ has $2$ constant eigenvalues and so the splitting lemma can be applied to decouple this system into two scalar equations. 
\end{example}
\begin{proof}(Proposition~\ref{nestedprop})
The eigenvalues of $A_0(x)$ admit a formal expansion in the fractional powers of $x$ in the neighborhood of $x=0$ (see, e.g.,~\cite{key56}).  We are interested only in their first nonzero terms. Let $\mu(x) = \sum_{j=0}^{\infty} \mu_j x^{j/s}$ be a nonzero eigenvalue of $A_0(x)$ with  $s \in \mathbb{N}^*$, $j \wedge s =1$, and whose leading exponent, i.e.\ smallest $j/s$ for which $\mu_j \neq 0$, is minimal among the other nonzero eigenvalues. Without loss of generality, we can assume that $s=1$, otherwise we set $x = t^{s}$. By~\cite{key41,key54}, there exists $T \in GL_n(\mathbb{C}(x))$ such that for $B_0(x) = T^{-1} A_0(x) T = B_{0,0} + B_{1,0} x + \dots$:
$$\theta_{B_0(x)}(\lambda) = {x}^{rank(B_{0,0})} \det (\lambda I + \frac{{B}_{0,0}}{x} + B_{1,0} )|_{x=0}$$
does not vanish identically in $\lambda$. Let $\nu >0$ denote the  valuation of $B_0(x)$ in $x$. By ~\cite[Proposition 1]{key54}, there are $n - \operatorname{deg} (\theta)$ eigenvalues of $B_0(x)$ whose leading exponents lie in $[\nu,\nu +1[$, and $ \operatorname{deg} (\theta)$ eigenvalues for which the leading exponent is equal to or greater than $\nu +1$. Then, applying the transformation $F = T(x) G$ to the system $[A_{\sigma_A}]$ yields
$$
[B_{\sigma_A}] \quad \quad \quad  \xi_A^{h} x^{p} \; \partial G\; = \;B(x, \xi_A) \;G\; =\;  \sum_{k=0}^{\infty} B_k(x) \xi_A^k  \; G
$$
where 
$B_{0}(x) := x^{\nu} (B_{0,0} + B_{1,0} x + \dots )$ and $B_{0,0}$ has $n - \operatorname{deg}(\theta)$ eigenvalues whose leading exponents lie in $[0,1[$. 
Let $\tilde{A}(x, \xi_A) = x^{-\nu} B(x,\xi_A)$ then the order of poles introduced in $x$ in $\tilde{A}(x, \xi_A)$ is at worst  $(\operatorname{span}_x (T) + \nu)$, where $\operatorname{span}_x(T)$ is the difference between the valuation and the degree of the polynomial transformation $T(x)$ in $x$. Hence, $\sigma_{\tilde{A}} \geq \sigma_A -\operatorname{span}_x (T) - \nu$  and $[B_{\sigma_A}]$ can be rewritten as
$$
 [\tilde{A}_{\sigma_{\tilde{A}}}] \quad \quad \quad  \xi_{\tilde{A}}^{h} x^{\tilde{p}} \; \partial G\; = \;\tilde{A}(x, {\xi}_{\tilde{A}}) \; G=  \sum_{k=0}^{\infty} \tilde{A}_k(x) {\xi}^k_{\tilde{A}}  \; G
$$
where ${\xi}_{\tilde{A}}= x^{\sigma_{\tilde{A}}} \varepsilon$, $\tilde{p} \in \mathbb{Z}$, $\tilde{A}_k(x) \in \mathcal{M}_n([[x]])$, $\tilde{A}_0(x) = x^{-\nu} B_0(x)$, and $\tilde{A}_{0,0} = B_{0,0}$. Then $\mu(x)$ is an eigenvalue of ${\tilde{A}}_0(x)$ with a minimal leading exponent and hence it is among those whose leading exponents lie in $[0,1[$. By our assumption $s=1$, and hence the leading exponent of $\mu(x)$ is zero and $\mu_0 \neq 0$. Since $\mu_0$ is a nonzero eigenvalue of $\tilde{A}_{0,0}$, it follows that the $\tilde{A}_{0,0}$ is non-nilpotent.
\end{proof}
We remark that the eigenvalues of $A_0(x)$ are the roots of the algebraic scalar equation $f(x, \lambda)= \det (A_0(x) - \lambda I_n) = 0$ and can be computed by Newton-Puiseux algorithm. The sought polynomial transformation $T$ can be computed via~\textsc{Isolde},~\textsc{miniIsolde},  or~\textsc{Lindalg}.
\begin{remark}
\label{defturnpt}
Proposition~\ref{nestedprop} leads to the following observation about the detection of a turning point:
Consider an input system $[A_{\sigma_A}]$ given by~\eqref{paramh} with $\sigma_A=0$, $p_A=0$, and restraining index $\rho_A=-1/\sigma_A$. If $[A_{\sigma_A}]$ has a turning point at $x=0$ then, by the end of the process of formal reduction of this system (i.e. whenever all decoupled subsystems are either of non-positive $\varepsilon$-rank or are scalar equations), at least one of its decoupled subsystems has a finite restraining index.
\end{remark}
\begin{example}
Consider the following system whose $\sigma_A=0$ and $\xi_A = \varepsilon$:
$$[A_{\sigma_A}] \quad \quad \quad \xi_A^2 \; \partial F \; =\;   \begin{bmatrix}  0 & 1& 0 \\ 0 & 0 & 1 \\  \xi_A & x & 0 \end{bmatrix} F \quad \text{where} \quad A_0(x)= \begin{bmatrix}  0 & 1& 0 \\ 0 & 0 & 1 \\0 & x & 0 \end{bmatrix} .$$
We first compute $s=2$. We set $x=t^2$ and compute 
$T = \begin{bmatrix} 0 & 0 & 1 \\ 0 & t & 0 \\ 0 & 0& t^2 \end{bmatrix}$.  Or equivalently, we consider $T = \begin{bmatrix} 0 & 0 & 1 \\ 0 & x^{1/2} & 0 \\ 0 & 0& x \end{bmatrix} .$ Then $F = T G$ yields:
 $$ \xi_A^{2} \; \partial G\; = x^{1/2} \{ \begin{bmatrix}  0 & 1& 0 \\ 1 & 0 & 0 \\  0 & 1 & 0 \end{bmatrix}  +  \begin{bmatrix}  0 & 0& 0 \\ 0 & 0 & 0 \\ 1 & 0 & 0 \end{bmatrix} x^{-3/2}  \xi_A  +  \begin{bmatrix}  0 & 0 & 0 \\ 0 & - 1 & 0 \\0 & 0 & -2 \end{bmatrix} x^{-1}  \xi_A^2 \}  \; G , $$
which we rewrite as:
 $$[\tilde{A}_{\sigma_{\tilde{A}}}] \quad \quad \quad \quad \xi_{\tilde{A}}^{2} x^{5/2} \; \partial G \; =  \{ \begin{bmatrix}  0 & 1& 0 \\ 1 & 0 & 0 \\  0 & 1 & 0 \end{bmatrix}  +  \begin{bmatrix}  0 & 0& 0 \\ 0 & 0 & 0 \\ 1 & 0 & 0 \end{bmatrix}\xi_{\tilde{A}}+  \begin{bmatrix}  0 & 0 & 0 \\ 0 & - x^2 & 0 \\0 & 0 & -2 x^2 \end{bmatrix} \xi_{\tilde{A}}^2\; \} \; G , $$
 where $\xi_{\tilde{A}}= x^{-3/2} \varepsilon$. Now that the leading term is a constant matrix with three distinct eigenvalues, we can proceed by applying the Splitting of Theorem \ref{splitparam}. 
\end{example}
In the above example, the system could be decoupled thoroughly into three scalar equations. However, we might encounter different scenarios as well, as illustrated in the following example.

\begin{example}\cite[p. 57]{key59}
\label{exmwasow}
We recall the system of Example \ref{exmwasow1} (which is equivalent to \eqref{equiveq}) with $\xi_A=\varepsilon$: $$[A_{\sigma_A}] \quad \quad \quad \xi_A^2  \; \partial F\; =  \; \begin{bmatrix}  0 & 1& 0 \\ 0 & 0 & 1 \\  \xi_A & 0 & x \end{bmatrix} F .$$ As mentioned before, we have a turning point at $x=0$. $A_{0,0}$ is nilpotent and the eigenvalues of  $A_0(x)$ are $0$ and $x$,  whence $s=1$. Let $T= \operatorname{diag} (1,x,x^2)$  then $F=TG$ yields:
$$\xi_A^2 \; \partial G \; =\;  x \{ \begin{bmatrix}  0 & 1& 0 \\ 0 & 0 & 1 \\ 0& 0 & 1\end{bmatrix}  +  \begin{bmatrix}  0 & 0& 0 \\ 0 & 0 & 0 \\ 1& 0 & 0\end{bmatrix} x^{-3} \xi_A +  \begin{bmatrix}  0 & 0& 0 \\ 0 & -1 & 0 \\ 0& 0 & -2 \end{bmatrix} x^{-2} \xi_A^{2} \; \} \; G .$$
Setting $\xi_{\tilde{A}}= x^{-3} \varepsilon $, the former can be rewritten equivalently as:
$$[\tilde{A}_{\sigma_{\tilde{A}}}] \quad \quad \quad \quad \xi_{\tilde{A}}^2 x^5 \; \partial G \;  =  \{ \begin{bmatrix}  0 & 1& 0 \\ 0 & 0 & 1 \\ 0& 0 & 1\end{bmatrix}  +  \begin{bmatrix}  0 & 0& 0 \\ 0 & 0 & 0 \\ 1& 0 & 0\end{bmatrix} \xi_{\tilde{A}} +  \begin{bmatrix}  0 & 0& 0 \\ 0 & - x^4& 0 \\ 0& 0 & -2 x^4 \end{bmatrix}  \xi_{\tilde{A}}^{2} \; \} \; G .$$
The leading constant matrix $\tilde{A}_{0,0}$ is no longer  nilpotent. Hence the system can be decoupled into two subsystems upon setting $G = T W$ where $$T=\begin{bmatrix}  1 & 0& 1 \\ 0 & 1 & 1 \\ 0& 0 & 1\end{bmatrix} + \begin{bmatrix}  -1& -1&0 \\ -1 & -1& -2 \\ -1& -1 & 0\end{bmatrix} \xi_{\tilde{A}} + O(\xi_{\tilde{A}}^{2}). $$  The resulting equivalent system then consists of the two decoupled lower dimension systems where $W= {[W^1, W^2]}^T$, $\sigma_{B} = \sigma_{C}=\sigma_{\tilde{A}}= -3$:
\begin{eqnarray*} & [B_{\sigma_{B}}]& \quad \quad \quad \quad \xi_{B}^2 x^5  \; \partial W^1 =  \{ \begin{bmatrix}  0 & 1\\ 0 & 0  \end{bmatrix}  +  \begin{bmatrix}  -1 & 0 \\ -1 & 0 \end{bmatrix} \xi_{B}  + \begin{bmatrix}  1 & -1 \\ 1 & -1 + x^4  \end{bmatrix} \xi_{B}^{2} + O(\xi_{B}^3) \; \}  \; W^1 .\\ 
& [C_{\sigma_{C}}]& \quad \quad \quad \quad  \xi_{C}^2 x^5 \;  \partial W^2 = \{ 1+ \xi_{C} + (1+ 2x^4) \xi_{C}^2 + O(\xi_{C}^3)\; \}\;  W^2 . 
\end{eqnarray*}
The exponential part (and hence formal solution) of the second subsystem is clearly $ \int \; \xi_{C}^{-2} x^{-5}  (1+ O(\xi_{C})) \; dx = \frac{1}{2} \varepsilon^{-2} x^2 (1 + O( \varepsilon^{-2} x^2))$. One can observe that the eigenvalue $x$ of the leading matrix coefficient  of the  input system $[A_{\sigma_A}]$ is recovered as expected.  This is in accordance with the exponential parts obtained for~\eqref{equiveq}. As for the first subsystem, $B_0(x)$ and $B_{0,0}$ are simultaneously nilpotent. 
Due to this dual nilpotency, in order to proceed in the formal reduction of $[B_{\sigma_{B}}]$, we will make use of the $\varepsilon$-rank reduction in the following subsection.
\end{example}
\subsection{$\varepsilon$-Rank reduction}
\label{moser}
We consider again system $[A_{\sigma_A}]$ given by~\eqref{param}. We assume without any loss of generality that $A_0(x)$ and $A_{0,0}$ are simultaneously nilpotent. In this section, we investigate the $\varepsilon$-rank reduction of the system, i.e.\ we seek to determine the minimal integer $\varepsilon$-rank among all systems equivalent to $[A_{\sigma_A}]$. If this minimal integer rank turns out to be non-positive then we continue the reduction by treating the system as in Section~\ref{hzero}. Otherwise, the minimal integer rank gives an upper bound for the $\varepsilon$-exponential order, which allows us to proceed to Section~\ref{exporder}.
 
In analogy to its unperturbed counterpart, we define the $\varepsilon$-Moser rank and $\varepsilon$-Moser invariant of system $[A_{\sigma_A}]$ to be the following rational numbers respectively:
\begin{eqnarray*} m_{\varepsilon}(A) &= &\; \text{max} \;  (0, h + \frac{\;rank\;(A_0(x))}{n}) \quad 
\text{and} \\ \mu_{\varepsilon} (A) &=& \; \text{min} \; \{ m_{\varepsilon} (T[A_{\sigma_A}]) \; \text{for all}\; T\; \text{ in}\; GL_n(\rm \mathcal{K})\}. \end{eqnarray*} 
\begin{definition}
System   $[A_{\sigma_A}]$ (the matrix $A(x, \xi_A)$ respectively) is called $\varepsilon$-reducible if $m_{\varepsilon}(A)~>~\mu_{\varepsilon} (A)$. Otherwise it is said to be $\varepsilon$-irreducible.
\end{definition}
If $m_{\varepsilon} (A) \leq 1$, then $h=0$. Hence, we restrict our attention to the case of $m_{\varepsilon} (A) > 1$. 

This definition is not to be mixed neither with the usual sense of Moser-irreducible unperturbed system $[A]$; nor with the sense of reduced system already employed in the literature, i.e.\ the system whose matrix is the leading coefficient matrix of $[A_{\sigma_A}]$: $ x^p \partial F = A_0(x) F$. In fact, the motivation behind Moser's work in~\cite{key19} was to determine the nature of the singularity (regular or irregular) of an unperturbed system $x^p \partial F = A(x) F,\; A(x) \in \mathcal{M}_n(\mathbb{C}[[x]])$. Consequently, the reduction of the so-called \textit{Poincar\'e rank} $p$ to its minimal integer value was investigated. However, given our perturbed system $[A_{\sigma_A}]$, it seems more plausible to reduce $h$, rather than $p$, to its minimal integer value, since a system with $h \leq 0$ can be treated as in Subsection~\ref{basiccases}.

Since $\mu_{\varepsilon} (A)$ cannot be computed from the outset, we aim in this section to generalize Moser's criterion to detect whether the $\varepsilon$-rank of a given perturbed system is minimal among all equivalent systems. We remark that Moser's criterion has been generalized as well to linear functional matrix equations in~\cite{key64}, and borrowed from the theory of differential systems in~\cite{key54} to investigate efficient algorithmic resolution of the perturbed algebraic eigenvalue-eigenvector problem. However, despite their utility and efficiency for such univariate systems, algorithms based on this criterion are not considered so far over bivariate fields.

In~\cite{key102}, we generalize the work of~\cite{key41,key19} for the case $\sigma_A=0$. Herein, we establish the following theorem without any restriction on $\sigma_A$ (see~\cite[Section 5.5, p. 101]{key1000} and~\ref{appc} for the proof). This yields an algorithm which takes an input system $[A_{\sigma_A}]$, and outputs a transformation $T \in GL_n(\rm \mathcal{K})$ and an equivalent system $[\tilde{A}_{\sigma_{\tilde{A}}}] := T[A_{\sigma_A}]$ which has a minimal $\varepsilon$-rank among all systems equivalent to $[A_{\sigma_A}]$.  
\begin{theorem}
\label{mosernilpotent1}
Consider system~\eqref{param} given by
$$
[A_{\sigma_A}] \quad \quad \quad \xi_A^{h} x^{p} \partial F = A(x, \xi_A) F=  \sum_{k=0}^{\infty} A_k(x) {\xi_A}^k  \; F ,
$$
and suppose that $h>1$ and $m_{\varepsilon}(A)>1$. The polynomial 
\begin{equation} \label{thetaparam} \theta_A(\lambda) := {\xi_A}^{rank\;(A_0(x))} \det (\lambda I + \frac{A_0(x)}{\xi_A} + A_1(x) )|_{\xi_A=0} \end{equation}
vanishes identically in $\lambda$, if and only if there exists a transformation $F = T G$, $T \in GL_n(\rm \mathcal{K})$, such that the equivalent system $[\tilde{A}_{\sigma_{\tilde{A}}}] := T [A_{\sigma_A}]$  has either a strictly lower $\varepsilon$-rank or a leading coefficient matrix  with a strictly lower algebraic rank. Moreover, $\sigma_{T}, \sigma_{\tilde{A}} \geq \sigma$ where  $$\sigma =  \begin{cases} \sigma_A + \frac{p-1}{h}  \quad \text{if} \quad p<1 \\  \sigma_A  \quad \text{otherwise} \end{cases} .$$
\end{theorem}

\begin{example}
\label{algoexm}
Consider the system $[A_{\sigma_A}] \quad \xi_A^4 \partial F = A(x, \xi_A) F $ with $\sigma_A= 0$ and $$A(x, \xi_A) = \begin{bmatrix} 2 x \xi_A^3  & 3 x^2 \xi^4 & 2x \xi_A^2 & (2x+1) \xi_A^5 \\ 0 &  \xi_A^4 & 0 & 0\\ 0 & 0 & \xi_A^2 & 0 \\ -2 x & 0 & 0 & 0 \end{bmatrix} . $$
Our algorithm computes the transformation $F=T G$, where 
$T= \begin{bmatrix} 0& \xi_T^3 & 0 & 0 \\ 0 &  0 & \xi_T & 0\\ \xi_T^3  & 0 &0 & 0 \\ 0 & 0 & 0 & 1 \end{bmatrix} ,$
with $\sigma_T = \sigma_A=0$ which results in an equivalent $\varepsilon$-irreducible system $[A_{\sigma_{\tilde{A}}}] \quad \xi_{\tilde{A}}^2 \partial G= \tilde{A}(x, \xi_{\tilde{A}}) G$ where
$$\tilde{A}(x, \xi_{\tilde{A}})= \begin{bmatrix} 1  & 0& 0& 0 \\ 2x &  2 x \xi_{\tilde{A}} & 3 x^2 & 2x+1\\ 0 & 0 & \xi_{\tilde{A}}^2 & 0 \\ 0 & -2 x \xi_{\tilde{A}} & 0 & 0\end{bmatrix} ,$$
and $\sigma_{\tilde{A}}=0$. One can observe that the $\varepsilon$-rank is diminished by two.
\end{example}
\begin{example}
\label{exmwasow1}
Consider the first subsystem resulting from the reduction of Example~\ref{exmwasow}:
$$[B_{\sigma_{B}}] \quad \quad \quad \quad \xi_{B}^2 x^5  \; \partial W^1 =  \{ \begin{bmatrix}  0 & 1\\ 0 & 0  \end{bmatrix}  +  \begin{bmatrix}  -1 & 0 \\ -1 & 0 \end{bmatrix} \xi_{B}  + \begin{bmatrix}  1 & -1 \\ 1 & -1 + x^4  \end{bmatrix} \xi_{B}^{2} + O(\xi_{B}^3) \; \}  \; W^1 .$$
It is easy to see that the leading coefficient and constant matrices coincide and that $\theta_B (\lambda) = 1$. Hence, the $\varepsilon$-rank which is equal to $2$ is minimal.
\end{example}
To proceed in the formal reduction of Example~\ref{exmwasow1}, one needs to introduce a ramification in $\varepsilon$ (resp. $\xi_B$). The computation of the necessary ramification will be discussed in Section~\ref{exporder}. To establish the results therein, the $\varepsilon$-Moser invariant of equivalent scalar equations is needed. We thus introduce it in Subsection~\ref{mosereqn}. Before proceeding however, we discuss the case $h=1$ which is not tackled by Theorem~\ref{mosernilpotent1}.
\subsubsection{The case $h=1$}
The limitation of Theorem~\ref{mosernilpotent1} to $h>1$ is of technical nature (see the proof of necessity). In this subsection, we treat the case $h=1$. We consider again the system given by~\eqref{param}:
$$
[A_{\sigma_A}] \quad \quad \quad \xi_A x^{p} \partial F = A(x, \xi_A) F =  \sum_{k=0}^{\infty} A_k(x) \xi_A^k \; F .  
$$
Let $h_{true}$ denote the minimal integer $\varepsilon$-rank which can be attained upon applying a transformation in $GL_n(\rm{\mathcal{K}})$ and $\omega_\varepsilon(A)s$ denote the $\varepsilon$-exponential order of $[A_{\sigma_A}]$. In this subsection, we give  a method to decide whether $h_{true} =1$ or $h_{true}=0$. We first apply the substitution $\varepsilon = \tilde{\varepsilon}^{n+1}$ to $[A_{\sigma_A}]$. We denote its $\tilde{\varepsilon}$-rank by $\tilde{h}$.  Since $\tilde{h}= n+1>1$, we can then apply the $\tilde{\varepsilon}$-rank reduction of Theorem~\ref{mosernilpotent1} and we denote by ${\tilde{h}}_{true}$ and $\omega_{\tilde{\varepsilon}}$ the true $\tilde{\varepsilon}$-rank of the resulting system and its $\tilde{\varepsilon}$-exponential order respectively. Under these notations, we have:
\begin{lemma}
\label{mosernilpotent3}
If ${\tilde{h}}_{true}=1$ then $h_{true}=0$.
\end{lemma}
\begin{proof}
Let $\omega_{\varepsilon} = \frac{\ell}{d}$ where $\ell$ and $d$ are coprime natural numbers. Then, on the one hand, $d, \ell \geq 0$, $d<n$ and $\omega_{\varepsilon} \leq h_{true}$. On the other hand, $\omega_{\tilde{\varepsilon}} = (n+1) \omega_{\varepsilon}$ and $\omega_{\tilde{\varepsilon}} \leq {\tilde{h}}_{true}$ (the minimal $\varepsilon$-rank bounds the $\varepsilon$-exponential order, see Corollary~\ref{qkappa}).  Hence, we have, 
$$\omega_{\tilde{\varepsilon}}= (n+1) \omega_{\varepsilon} = (n+1) \frac{\ell}{d} \geq (n+1) \frac{\ell}{n} \geq (1 + \frac{1}{n}) \ell .$$
Hence, if $\ell \neq 0$ then
$$  1 < (1 + \frac{1}{n}) \ell \leq \tilde{\omega}_{\tilde{\varepsilon}} \leq {\tilde{h}}_{true}.$$
\end{proof}
Lemma~\ref{mosernilpotent3} moves the problem from the case $h=1$ to the case $h>1$. It can be restated as follows: Within the formal reduction, whenever the case $h=1$ is encountered, a ramification $\varepsilon = \tilde{\varepsilon}^{n+1}$ is applied so that the reduction of Theorem~\ref{mosernilpotent1} can be applied. After this reduction, one of the following two cases arises:
\begin{itemize}
\item If the $\varepsilon$-rank of the resulting system is less than or equal to one then $h_{true}$ of the original system is zero, and so we stop our reduction (the $\varepsilon$-exponential part have been computed completely). If one wishes to continue reduction, then the classical Arnold-Wasow approach can be employed: First, put the system in Arnold's form~\cite{key10} and then use the Iwano-Sibuya's polygon to determine the proper polynomial transformation required to arrive at an equivalent system whose $\varepsilon$-rank is zero~\cite{key62}. The reduction then continues as explained in Subsection~\ref{hzero}. 
\item Otherwise, we proceed with the reduction: If the leading coefficient matrix of the ramified system has two distinct eigenvalues, then we treat turning points (if any) and apply splitting lemma. Otherwise, we proceed to Section~\ref{exporder}. 
\end{itemize}
We remark that in the implementation, we first try to find whether there exist constant vectors in the left null space of $G_{A}(\lambda)$ (see~\ref{appc}). If such vectors do exist, we use them to construct a transformation which might reduce the rank. If the $\varepsilon$-rank remains one after all the constant vectors are exhausted, we compute a candidate for the ramification in $\xi_A$ from the characteristic polynomial. If this candidate, after applying again the $\varepsilon$-rank reduction results in a system with non-nilpotent leading coefficient, we keep it. Otherwise, we use Arnold-Wasow approach.

\subsection{ $\varepsilon$-Moser invariant of a scalar equation}
\label{mosereqn}
Due to the equivalence between a scalar differential equation $[a_{\sigma_a}]$ and its companion system $[A_{\sigma_A}]$, it is natural to define and question the $\varepsilon$-Moser invariant of the former, in the hope to gain more insight into the problem. This is the goal of this subsection which is fulfilled by generalizing the analogous notion discussed for unperturbed scalar linear differential equations in~\cite[Part IV]{key19}. 

We consider again the singularly-perturbed linear differential equations given by~\eqref{scalarparammm}
$$ [a_{\sigma_a}] \quad \quad \quad \partial^n f \; +\;  a_{n-1}(x, \xi_a) \; \partial^{n-1} f \; +\; \dots \;+\; a_{1}(x, \xi_a) \; \partial f \;+\;  a_{0}(x, \xi_a) f \;=\;0    ,$$
where $a_{n}(x, \xi_a) = 1$. We prove the following proposition:
\begin{proposition}
\label{moserscalarparam}
Given the differential equation $[a_{\sigma_a}]$. Let $\tau, \nu$ be the smallest integers such that $val_\varepsilon \; (a_i (x, \xi_a)) \geq (i-n) (\tau -1) - \nu .$ In other words, let 
\begin{eqnarray}
\label{scalarparam11} \kappa & = &\; min\; \{ \varrho  \in \mathbb{N} | \; val_\varepsilon\;(a_i) + (n-i) \varrho \geq 0,\; 0 \leq i \leq n \} \\
\label{scalarparam22} \nu & = &\; max\; \{ (i-n) (\tau-1) - val_\varepsilon\;(a_i),\; 0 \leq i \leq n \}
\end{eqnarray}
Then the $\varepsilon$-Moser invariant of the system given by the associated companion matrix is $$\mu_{\varepsilon} (a) = \kappa +\frac{\nu}{n} .$$  
\end{proposition}
\begin{remark}[Geometric interpretation] Consider again the $\varepsilon$-polygon $\mathcal{N}_{\varepsilon}(a)$ of $[a_{\sigma_a}]$ constructed in Subsection~\ref{macutanpolygon}) in a $(U, V)$-plane. First, one can construct the straight line passing through the point $(n, 0)$, with the smallest integer slope $\kappa$, which stays below $\mathcal{N}_{\varepsilon}(a)$.  It is the straight line of equation $V = \kappa (U-n)$. Then, one finds among all the parallel lines of slope $(\kappa-1)$, the highest straight line with an integer $V$-intercept which stays below $\mathcal{N}_{\varepsilon}(a)$. The $V$-intercept of the latter is $-\nu - n (\kappa-1)$ where $\nu$ is an integer. In other words, the latter has the equation $V = (U-n) (\kappa-1) - \nu$.
\end{remark}
\begin{example}
\begin{figure}
\centering
\includegraphics[width=60mm]{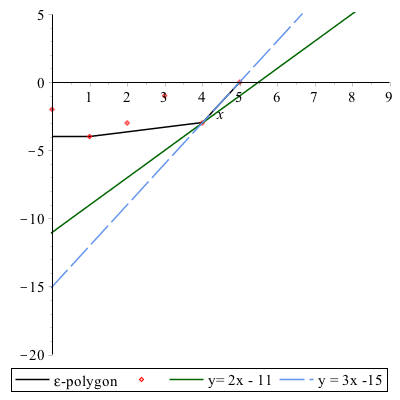}
\caption{Geometric Interpretation}
\label{figure3}
\end{figure}
Let $\sigma_a=-1$ and $\xi_a = x^{\sigma_a} \varepsilon$. We consider the scalar differential equation
$$[a_{\sigma_a}] \quad \partial^5 f + (x \xi_a^{-3} + x) \partial^4 f + 2x^2 \xi_a^{-1} \partial^3 f + (\xi_a^{-3} + 1) \partial^2 f + (-3 \xi_a^{-4} + x^2\xi_a^{-2}) \partial f - \xi_a^{-2} f = 0.$$
First, we plot the points $(i, val_\varepsilon(a_i)) , 0 \leq i \leq n=5$ and the $\varepsilon$-polygon (see Figure~\ref{figure3}). Next, we plot the straight line passing through the point $(n, 0) = (5,0)$ with smallest integer slope such that it stays below all these point. This yields $V = 3U - 15 = 3 (U-5)$  and $\tau = 3$. Then, we plot the straight line of slope $\tau -1= 2$ which stays below these points: $V= 2 U - 11 = 2(U-5) - 1$. This yields $\nu = 1$. Finally, we consider $\gamma = (\gamma_0, \dots, \gamma_4) = (-11, -9, -7, -5, -3)$.  We observe that $i_0 = n - \nu = 4$ and so the equivalent system is given by:
$$[A_{\sigma_A}] \quad \quad \quad x \xi_A^{3} \partial W = \begin{bmatrix} 11 \xi_A^3 & x \xi_A & 0 & 0 & 0 \\ 0 & 9 \xi_A^3 & x \xi_A & 0 & 0 \\ 0 & 0 & 7 \xi_A^3 & x \xi_A & 0 \\ 0 &0 & 0 & 5 \xi_A^3 & x \xi_A \\ x \xi_A^9 & 3 x \xi_A^5 - x^3 \xi_A^7 & -x \xi_A^4  - x \xi_A^7 & -2x^3\xi_A^4 & -x^2- (x^2 - 3) \xi_A^3 \end{bmatrix} W ,$$
with $\sigma_A = \sigma_a = -1$. One can verify that the system $[A_{\sigma_A}]$ is $\varepsilon$-irreducible since  $\theta_A(\lambda) = x \lambda^4.$
\end{example}
\begin{proof}(Proposition \ref{moserscalarparam}) For clarlity, we set $\xi := \xi_a$. We define  $$\gamma_i = max\; \{ \kappa (i-n) , (\kappa-1) (i-n) - \nu \} ,  \; 0 \leq i \leq n-1 .$$
Then, by construction,  $\gamma_i \leq val_{\varepsilon} \; (a_i)$ and equality is attained for at least one $i \geq n - \nu$ (otherwise $\nu$ and $\kappa$ can be minimized). Let $i_0$ represent the smallest integer such that $\gamma_{i_0} = \kappa (i_0 -n)$, then $i_0 = n -\nu$. Geometrically, the $\gamma_i$'s represent a broken line dominated  by the $(i, val_{\varepsilon} \; (a_i))$. They will aid in the  construction of an $\varepsilon$-irreducible system which is equivalent to~\eqref{scalarparammm}. In fact, let
$$w_{i+1} =\xi^{ \gamma_i}\; \partial^i f , \quad 0 \leq i \leq n-1.$$
then we have 
$$ \begin{cases} 
\partial w_{i+1} =\xi^{-\kappa} \; [ \begin{cases}\xi w_{i+2} \quad  \quad \quad 0 \leq i < i_0 \\ w_{i+2} \quad \quad \quad  \quad \quad \quad i_0 \leq i \leq n-2 \end{cases} \;+\;\xi^{\kappa} \frac{\sigma_a \gamma_i}{x}\; w_{i+1}  ] , \\[5pt]
 w_{n} =\xi^{\gamma_{n-1}} \; \partial^{n-1}f =\xi^{- \kappa}  \; \partial^{n-1}f   , \\[5pt]
\partial w_{n} =  \frac{\sigma_a \gamma_{n-1}}{x}\; w_{n} +\xi^{- \kappa}  \partial^n f =\xi^{-\kappa} [\xi^{\kappa} \; \frac{\sigma_a \gamma_{n-1}}{x}\; w_{n} + \sum_{i=0}^{n-1} \; \alpha_{i}(x, \varepsilon) w_{i+1} ] . 
\end{cases}$$
where  ${\alpha}_i (x, \xi) = - a_i (x, \varepsilon)\xi^{-\gamma_{i}} =  \sum_{k \in \mathbb{Z}}^{\infty} a_{i,k}(x)\xi^{k-\gamma_i}.$ Let $W= (w_1, \dots w_n)^T$ and $A(x, \xi) =$
$$ \begin{bmatrix} \sigma_a \gamma_{0}\xi^{\kappa} & x\xi & \quad & \quad & \quad & \quad & \quad \\ \quad & \quad & \ddots & \quad & \quad & \quad & \quad \\ \quad & \quad & \quad & \quad & \quad & \quad & \quad \\ \quad & \quad & \sigma_a \gamma_{i_0-1}\xi^{\kappa} &x\xi & \quad & \quad & \quad \\ \quad & \quad & \quad &  \sigma_a \gamma_{i_0}\xi^{\kappa}  & x& \quad & \quad \\
\quad & \quad & \quad & \quad & \sigma_a \gamma_{i_0+1}\xi^{\kappa}  & \quad & \quad \\  \quad & \quad & \quad & \quad & \quad & \quad & \quad  \\  \quad & \quad & \quad & \quad & \quad & \quad & x \\ x {\alpha}_0 & x {\alpha}_1 & \dots & x {\alpha}_{i_0} & x {\alpha}_{i_0 + 1} & \dots & x {\alpha}_{n-1} +  \sigma_a \gamma_{n-1}\xi^{\kappa} \end{bmatrix} . $$
Then one can verify that,
$$[A_{\sigma_a}] \quad \quad \quad  x \xi^{\kappa}  \partial W = A(x, \xi) W.$$
It remains to prove that this system is $\varepsilon$-irreducible. We first remark  that $A(x, 0)$ has rank $\nu = n - i_0$ due to the linear independence of its last $\nu$ rows. In fact, it's clear that it has $\nu -1$ linearly independent rows (the rows with the $x$'s). Moreover, ${\alpha}_{i} (x, \xi) = - a_{i} (x, \xi)\xi^{-\gamma_{i}} = - a_{i} (x, \xi)\xi^{-val_\varepsilon (a_{i})}$, and so ${\alpha}_{i} (x, 0) \neq 0$ for at least one $i \in \{ 0, \dots, i_0\}$. Thus the last $\nu$ rows are linearly independent. Setting $\xi = x^{\sigma_a} \varepsilon$ we have:
\begin{eqnarray*} &\theta_A(\lambda) &= \xi^{\nu} \det ( \lambda I + \frac{ A_0(x)}{\xi} + A_1 (x)) |_{\xi =0 } \\ &= & \xi^{\nu} \; \begin{vmatrix}  \lambda + \delta (\sigma_a \gamma_{0}) & x & \quad & \quad & \quad & \quad & \quad \\ \quad & \quad & \ddots & \quad & \quad & \quad & \quad \\ \quad & \quad & \quad & \quad & \quad & \quad & \quad \\ \quad & \quad &\lambda  + \delta (\sigma_a \gamma_{i_0-1})  &x & \quad & \quad & \quad \\ \quad & \quad & \quad &  \lambda  + \delta (\sigma_a \gamma_{i_0})   & \frac{x}{\xi} & \quad & \quad \\
\quad & \quad & \quad & \quad & \lambda  + \delta (\sigma_a \gamma_{i_0 +1})  & \quad & \quad \\  \quad & \quad & \quad & \quad & \quad & \quad & \quad  \\  \quad & \quad & \quad & \quad & \quad & \quad &  \frac{x}{\xi} \\ \frac{x {\alpha}_0}{\xi} & \frac{x {\alpha}_1}{\xi}  & \dots & \frac{x {\alpha}_{i_0}}{\xi}  & \frac{x {\alpha}_{i_0 + 1}}{\xi}  & \dots & \lambda + \delta (\sigma_a \gamma_{n-1}) + \frac{x {\alpha}_{n-1}}{\xi} \end{vmatrix}  \; |_{\xi =0 } \\
&=&   \xi^{\nu -n + i_0} \; \begin{vmatrix}  \lambda  + \delta (\sigma_a \gamma_{0}) & x & \quad & \quad \\ \quad & \quad & \ddots & \quad\\ \quad & \quad & \quad & \quad \\ \quad & \quad &\lambda  + \delta (\sigma_a \gamma_{i_0-1})  &x \\ x {\alpha}_0& x {\alpha}_1 & \dots & x {\alpha}_{i_0} \end{vmatrix} , \end{eqnarray*}
where $\delta=1$ if $\kappa =1$ and is zero otherwise.
Hence, since ${\alpha}_{i} (x, 0) \neq 0$ for at least one $i \in \{ 0, \dots, i_0\}$,   $\theta_A(\lambda)$ does not vanish identically in $\lambda$ and its highest possible degree is $i_0$. It follows from Theorem~\ref{mosernilpotent1} that the system is $\varepsilon$-irreducible and $\mu_{\varepsilon} (a) = \mu_{\varepsilon} (A)= \kappa + \frac{\nu}{n}$. 
\end{proof}

\begin{corollary}
\label{corocyclic}
Given an $\varepsilon$-irreducible system 
$\; [A_{\sigma_A}] \quad 
\xi_A^{h} x^{p} \partial F = A(x, \xi_A) F  
$. Let $[a_{\sigma_a}]$ be the scalar equation corresponding to a  companion system equivalent to $[A_{\sigma_A}]$. Then for $\kappa$ and $\nu$ defined in \eqref{scalarparam11} and \eqref{scalarparam22} we have: $\kappa=h$ and $\nu= rank (A_0)$. 
\end{corollary}
\begin{proof}
 By the $\varepsilon$-irreducibility of $[A_{\sigma_A}]$ and its equivalence to the companion system, we have:
   $$\kappa+ \frac{\nu}{n} = \mu_{\varepsilon}(a) = \mu_{\varepsilon}(A) = m_{\varepsilon}(A) = h+\frac{r}{n}$$
   Consequently, $\kappa=h$ and $\nu= rank (A_0)$. 
\end{proof}
\begin{corollary}
\label{qkappa}
Under the notations of Corollary \ref{corocyclic}, we have: 
$$\kappa -1 + \frac{\nu}{n} \leq \omega_{\varepsilon}(A) \leq \kappa  \quad \text{and} \quad h -1 + \frac{rank(A_0)}{n} \leq \omega_{\varepsilon}(A) \leq h .$$
\end{corollary}
\begin{proof}
Follows from \eqref{scalarparam11}, \eqref{scalarparam22}, and \eqref{epsilonkappa}. 
\end{proof}
\section{Computing the $\varepsilon$-exponential order}
\label{exporder}
We consider again system $[A_{\sigma_A}]$ which is given by~\eqref{param}: 
\begin{equation}
\label{katzmain4}
 [A_{\sigma_A}] \quad \quad \quad  x^p \xi_A^h\; \partial F\; = \; A(x, \xi_A) F = \sum_{k=0}^{\infty} A_k(x) \xi_A^k \;F .\end{equation}
In this section, we assume without loss of generality that $h>0$, $[A_{\sigma_A}]$ is $\varepsilon$-irreducible, and $A_0(x), A_{0,0}$ are both nilpotent. The goal of this section is to compute the $\varepsilon$-exponential order $\omega_{\varepsilon}(A)$ and $\varepsilon$-polynomials of the system, which give indispensable information within formal reduction. In particular, the former determines the ramification in $\varepsilon$ which can lead to a system whose leading matrix is non-nilpotent, so that the process of formal reduction can be resumed. This leads to the recursive Algorithm~\ref{algorithmkatzmain3}. Eventually, we can construct fundamental matrices of formal solutions in any given subdomain. We illustrate this algorithm by examples and motivate items for further investigation in~\ref{appe}. Let
\begin{equation}
\label{det11}
 \det (\lambda I - \frac{A(x, \varepsilon)}{x^p \xi_A^{h}}) = \lambda^n + \alpha_{n-1}(x, \xi_A) \lambda_{n-1} + \dots + \alpha_0 (x, \xi_A)  .\end{equation}
such that $\alpha_n =1$ and $\alpha_i(x, \xi_A) = \sum_{j= val_{\varepsilon}(\alpha_i)}^{\infty} \alpha_{i,j}(x) \xi_A^j \in {\rm \mathcal{K}}$ for $i \in \{0, \dots, n\}$. We define the $\varepsilon$-polygon $\mathcal{N}_{\varepsilon} (A)$ of $[A_{\sigma_A}]$ by taking $P_{\varepsilon}(A)$ of Section~\ref{macutanpolygon} to be the union of $P(i, \;val_{\varepsilon}\; (\alpha_i (x, \xi_A))$ for $i \in \{0, \dots, n\}$. We can then prove the following:
\begin{theorem}
\label{katzmain}
Consider the $\varepsilon$-irreducible system $[A_{\sigma_A}]$ given by~\eqref{katzmain4} with $h>0$ and~\eqref{det11}. If  $h > n - rank (A_0(x))$ then the $\varepsilon$-exponential order of $[A_{\sigma_A}]$ is given by 
$$\omega_{\varepsilon} (A)= \maxi_{0 \leq i < n}\; (\frac{- val_{\varepsilon}(\alpha_i)}{n - i}),$$
and its corresponding $\varepsilon$-polynomial is given by the algebraic equation
$$E_{\varepsilon}(X) = \sum_{k=0}^{\ell} x^{\sigma_A \cdot val_{\varepsilon} (\alpha_{i_k}) } {\alpha}_{i_k, val_{\varepsilon} (\alpha_{i_k})}\; X^{(i_k - i_0)}$$ 
where $0 \leq i_0 < i_1 < \dots < i_\ell = n$ denote the integers $i$ for which $\omega_{\varepsilon}(n-i)= - val_{\varepsilon}(\alpha_i)$ (i.e.\ lie on the edge of slope $\omega_{\varepsilon}$ of the $\varepsilon$-polygon $\mathcal{N}_{\varepsilon} (A)$ of $[A_{\sigma_A}]$); and ${\alpha}_{i, val_{\varepsilon} (\alpha_{i})}(x) = \xi_A^{-val_{\varepsilon}\; (\alpha_{i})}\; \alpha_{i} (x, \xi_A) |_{\xi_A = 0}$.
\end{theorem}
 The proof can be established by adapting to the parametrized setting the proofs of~\cite[Lemma 3, Lemma 4, Proposition 1, Theorem 1]{key24} (see~\ref{appd}). Not only does this theorem compute these invariants of the system, but it also allows a further reduction of the system as follows: Suppose that $\omega_{\varepsilon}(A)=\frac{\ell}{d}$ where $\ell,d$ are relatively prime positive integers. One can then set $\tilde{\varepsilon} = \varepsilon^{1/d}$ (or equivalently\footnote{In fact, we want to set $ \tilde{\varepsilon}^d = \varepsilon$. And so we choose an integer $s$ such that $\tilde{\xi_A}^d = \xi_A x^s$ which is equivalent to $x^{\sigma d} \tilde{\varepsilon}^d = x^{\sigma +s} \varepsilon$.} $\tilde{\xi_A} = \xi_A^{1/d} x^{\sigma (d-1)/d}$) in $[A_{\sigma_A}]$ and perform the $\tilde{\varepsilon}$-rank reduction (the minimal $\varepsilon$-rank is $\ell$). This will lead to an equivalent system whose $\varepsilon$-rank is $\ell$ and whose leading coefficient matrix  has at least $d$ distinct eigenvalues. The system can thus be decoupled into subsystems of lower dimensions. By repeating these procedures for each of the resulting subsystems, we can decouple the initial system into subsystem(s) of lower dimension(s) or zero $\varepsilon$-rank. This leads to the recursive algorithm of Section~\ref{algokatzmain}. \\ 
We remark that the condition $h > n - rank (A_0(x))$ in Theorem~\ref{katzmain} is non-restrictive. We can always arrive to a system satisfying this condition by a generalization of~\cite[Lemma 5]{key24}, which is based on applying the affinity $(U,V) \rightarrow (U, d\;V)$ for some integer $d$ to the $\varepsilon$-polygon. In fact, we can prove the following lemma:
\begin{lemma}
\label{katzlemma}
Consider an $\varepsilon$-irreducible system $[A_{\sigma_A}]$ given by~\eqref{katzmain4} with $h>0$ and $r=\; rank\; (A_0(x))$. Let $d$ be an integer such that $d \geq \frac{n}{h-1 +\frac{r}{n} }$ and let $$
\tilde{\xi_A}^{\tilde{h}} x^{\tilde{p}} \partial G= \tilde{A}(x, \tilde{\xi_A}) G .  
$$
be the $\tilde{\varepsilon}$-irreducible differential system obtained by the ramification $\varepsilon = {\tilde{\varepsilon}}^d$ (or equivalently $\tilde{\xi_A} = \xi_A^{1/d} x^{\sigma_A (d-1)/d}$ and $\sigma_{\tilde{A}} = \sigma_A (d-1)$) and performing $\tilde{\varepsilon}$-rank reduction. Then $\mu_{\tilde{\varepsilon}} (\tilde{A}) = \tilde{h} + \frac{\tilde{r}}{n}$, where $\tilde{r} = \;rank\;(\tilde{A}_0(x))$ and $\tilde{h} + \tilde{r} > n$. \end{lemma}
\begin{proof}
By Corollary~\ref{qkappa} we have, $h-1+\frac{r}{n} \leq \omega_\varepsilon(A) \leq h$. And due to $\varepsilon = \tilde{\varepsilon}^d$, we have $\omega_{\varepsilon} (A)= \frac{ \omega_{\tilde{\varepsilon}} (\tilde{A})}{d}$. Hence, $\omega_{\tilde{\varepsilon}} (\tilde{A}) \geq d (h-1+\frac{r}{n}) \geq n $. But $\tilde{h} \geq \omega_{\tilde{\varepsilon}} (\tilde{A})$ and $\tilde{r} \geq 1$, which yields $\tilde{h} + \tilde{r} \geq  \omega_{\tilde{\varepsilon}} (\tilde{A}) + 1\geq d (h-1+\frac{r}{n})  + 1  \geq n+1 >n$.
\end{proof}

\section{Formal reduction algorithm}
\label{algokatzmain}
\begin{algorithm}
\caption{ $\textsc{Exp\_param} (h, p, \sigma_A, A(x, \xi_A))$ : Computes the $\varepsilon$-exponential part of system $[A_{\sigma_A}]$ (performs formal reduction)}.
\label{algorithmkatzmain3}
\textbf{Input:} $h, p, \sigma_A, \;A(x,\xi_A)$ of system~\eqref{param} ($p=\sigma_A=0$ for system \eqref{paramh})  \\
\textbf{Output:} $Q(x^{1/s}, \varepsilon^{-1/d} )$
\begin{algorithmic}
\State  $Q \gets  \operatorname{diag} (0, \dots 0)$;
\While {$h >0$ and $n \neq 1$} \do
\State \If {$A_{0,0}$ has at least two distinct eigenvalues}
\State Apply the $\varepsilon$-block diagonalization of Section~\ref{blockparam};
\State $\textsc{Exp\_param} (h, p, \sigma_A, \tilde{A}_{11}(x,\varepsilon))$; Update $Q$ and $\sigma_A$;
\State $\textsc{Exp\_param} (h, p, \sigma_A, \tilde{A}_{22}(x,\varepsilon))$; Update $Q$ and $\sigma_A$;
\ElsIf {$A_{0,0}$ has one non-zero eigenvalue}
\State Update $Q$ from the eigenvalues of $A_{0,0}$;
\State $A(x,\xi_A) \gets $ perform eigenvalue shifting~\eqref{evshiftparam}; ($A_{0,0}$ is now nilpotent); 
\State $\textsc{Exp\_param} (h, p, \sigma_A, A(x,\xi_A))$; Update $Q$;
\ElsIf {$A_0(x)$ is not nilpotent} 
\State $A(x,\xi_A), p, \sigma_A, \gets $ apply turning point resolution of Subsection~\ref{turnpt}; ($A_{0,0}$ is now non-nilpotent and $\sigma_A$ is updated); 
\State $\textsc{Exp\_param} (h, p, \sigma_A, A(x,\xi_A))$; Update $Q$;
\Else
\State $A(x,\xi_A), h \gets$ $\varepsilon$-rank reduction of Section~\ref{moser}; 
\If {$h >0$ (i.e. $h_{true} >0$)}
\If {$A_{0,0}$ has at least two distinct eigenvalues}
\State Apply the $\varepsilon$-block diagonalization of Section~\ref{blockparam};
\State $\textsc{Exp\_param} (h, p, \sigma_A, \tilde{A}_{11}(x,\varepsilon))$; Update $Q$ and $\sigma_A$;
\State $\textsc{Exp\_param} (h, p, \sigma_A, \tilde{A}_{22}(x,\varepsilon))$; Update $Q$ and $\sigma_A$;
\ElsIf {$A_{0,0}$ has one non-zero eigenvalue}
\State Update $Q$ from the eigenvalues of $A_{0,0}$;
\State $A(x,\xi_A) \gets $ perform eigenvalue shifting~\eqref{evshiftparam}; ($A_{0,0}$ is now nilpotent); 
\State $\textsc{Exp\_param} (h, p, \sigma_A, A(x,\xi_A))$; Update $Q$;
\ElsIf {$A_0(x)$ is not nilpotent} 
\State $A(x,\xi_A), p, \sigma_A, \gets $ apply turning point resolution of Subsection~\ref{turnpt}; ($A_{0,0}$ is now non-nilpotent and $\sigma_A$ is updated); 
\State Update $Q$ from the eigenvalues of $A_{0,0}$;
\State $A(x,\xi_A) \gets $ perform Eigenvalue shifting; ($A_{0,0}$ is now nilpotent);
\State $\textsc{Exp\_param} (h, p, \sigma_A, A(x,\xi_A))$; Update $Q$;
\Else 
\State Use Theorem~\ref{katzmain} of Section~\ref{katzmain}; 
\State $\omega_{\varepsilon} = \frac{\ell}{d}$; $\varepsilon \gets \varepsilon^d$; 
\State $h, A(x,\xi_A) \gets $ Apply $\varepsilon$-rank reduction of Section~\ref{moser} ($h \gets \ell$); 
\State $\textsc{Exp\_param} (h, p, \sigma_A, A(x,\xi_A))$; Update $Q$;
\EndIf
\EndIf
\EndIf
\EndWhile \\
\Return{(Q)}.
\end{algorithmic}
\end{algorithm}
With the algorithms of Splitting lemma, turning point resolution, $\varepsilon$-rank reduction, and as can be verified by Algorithm~\ref{algorithmkatzmain3} below, we have re-established constructively the well-known  general form for a fundamental matrix of formal outer solutions for an input system $[A_{\varepsilon}]$ given by~\eqref{paramh}:
 \begin{equation} \label{solkatzmain3}
F =  (\sum_{k=0}^{\infty} \Phi_k(x^{1/s}) \; (x^\sigma \varepsilon)^{k/d}) \exp( \int \mathcal{Q} (x^{1/s}, \varepsilon^{-1/d})), \end{equation}
where $s, d$ are positive integers; $\sigma$ is a nonpositive rational number; $\mathcal{Q}$ is a diagonal matrix whose entries are polynomials in $\varepsilon^{-1/d}$ with coefficients in $\mathbb{C}((x^{1/s}))$. We refer to $Q:= \int \mathcal{Q}$ as the $\varepsilon$-\textit{exponential part} (logarithms in a root of $x$ might arise as a result of integration); and the entries of the $\Phi_k(x^{1/s})$'s are root-meromorphic in $x$ (see~\cite[Introduction]{key210} or~\cite{key202}). 
\begin{remark}
Under the notations and statements of Theorem~\ref{katzmain}, the leading term of $Q$ is given by
$$\frac{1}{\varepsilon^{\omega_{\varepsilon}(A)}}\; \int \operatorname{diag} (X_1(x), \dots, X_{deg \; (E_{\varepsilon})}(x), 0, \dots, 0) \; dx,$$
where the ${X_i}$'s denote the roots of the $\varepsilon$-polynomial $E_{\varepsilon}(X)$. 
A similar statement can be stated in terms of the eigenvalues of the leading matrix coefficient, $x^p$, and $\xi_A^h$.
\end{remark}
 We sum up our main results in the formal reduction algorithm, Algorithm~\ref{algorithmkatzmain3}, which computes the $\varepsilon$-exponential part and consequently a fundamental matrix of formal solutions~\eqref{solkatzmain3} in a given subdomain. We recall that: 
\begin{itemize}
\item If $n=1$ then we proceed by integration up to the first $h$ terms;
\item If $h \leq 0$ then we follow Subsection~\ref{hzero} (using \textsc{ISOLDE}, \textsc{miniIsolde}, or \textsc{Lindalg}).
\end{itemize}

\section{Conclusion and further investigations}
\label{conclusion}
In this article, we give an algorithm, implemented in \textsc{Maple}\footnote{The package is available at: $http://www.specfun.inria.fr/~smaddah/Research\_.html$}, which computes a fundamental matrix of outer formal solutions for singularly-perturbed linear differential systems in a neighborhood of a turning point. The subprocedures discussed are stand-alone algorithms (Splitting, rank reduction, and resolution of turning points). The formal reduction algorithm is based on the generalization of the algorithm given in~\cite{key24} for the unperturbed counterparts of such systems.

Our results are presented in the formal setting. However, the growth of the order of poles in $x$ is tracked within the formal reduction. This gives information on the restraining index of the system, and consequently, an adequate stretching transformation can be chosen and performed. This furnishes the first step in resolving Iwano's first problem and computing, for an input system, $[\rho]$ which is given by~\eqref{radii}:
\begin{equation*}
[\rho] \quad \quad \quad \quad 0 =\;{\rho}_0 < \rho_1 < \rho_2 < \dots < \rho_m . 
\end{equation*} 

In~\ref{appe}, we try to motivate the employment of our proposed formal reduction algorithm to resolve the same problem for a general system: Suppose that we start the reduction with the input system $[A_{\sigma_A}]$, $\sigma_A=0$, given by~\eqref{param} as an input. Then, upon applying Algorithm~\ref{algorithmkatzmain3}, we can determine the outer solutions and the final restraining index $\sigma_{final}$ which allows full reduction. If the restraining index is nonzero then we have a turning point and $\rho_m = -1/\sigma_{final}$. We can then perform the stretching $\tau = x \varepsilon^{-\rho_1}$ in the input system $[A_{\sigma_A}]$, and apply to it  Algorithm~\ref{algorithmkatzmain3}. This determines $\rho_2$. We show an example where the same process can be repeated with $\rho_2$ to determine $\rho_3$, and iteratively we reach a system whose final restraining index is zero. In future work, we hope to investigate this approach and its correctness. 

Another field of investigation is the relation between the algebraic eigenvalues of the matrix $A(x,\varepsilon)$ of such systems and the exponential part of the solution (see Example~\ref{introexm2} and ~\cite[Chapter 2]{key1000}). One can then benefit from the existing work on fractional power series expansions of solutions of bivariate algebraic equations, to compute the restraining index (see, e.g. \cite{key501,key502} and references therein). On the other hand, it would be interesting to investigate a differential-like reduction for a two-parameter perturbation of a JCF (see ~\cite[Chapter 2]{key1000} and references therein for the one-parameter case).  In fact, one can observe that the main role in the reduction process is reserved to the similarity term of $T[A_{\sigma_A}]$, i.e. $T^{-1} A_{\sigma_A} T$ rather than $T^{-1} \partial T$. Hence, for a non-differential operator where $T[A_{\sigma_A}]= T^{-1} A_{\sigma_A} T$, the discussion is not expected to deviate substantially from the discussion presented here for a differential one. 

The generalization of other notions and efficient algorithms can be investigated as well (see, e.g.\ simple system \cite{key25,key40} and ~\cite[Appendix A]{key1000}). 

In Examples \ref{exmbender1} and \ref{exmbender2} we treat input systems with $\sigma_A <0$. The algorithms generalized herein have been generalized to treat systems with an essential singularity~\cite{key4954} and difference systems as well(see, e.g.\ \cite{key101,key64}).  This motivates investigating the adaptation of our proposed algorithms at least to the difference setting.

There remain as well the questions which fall under the complexity, e.g.\ studying the complexity of this algorithm; the bounds on number of coefficients needed in computations; developing efficient algorithms for computing a cyclic vector; and comparing it to the results of our direct algorithm. 

And finally, to give a full answer on the behavior of the solutions, the related problems of connection, matching, and secondary turning points, are to be studied as well (see, e.g.\ \cite{key203,key201} and references therein). 

\section{Acknowledgements}
We would like to thank Dr. Reinhard Schaefke for pointing out inaccuracies in the previous version. Our discussions with him highly contributed to improving the quality of this manuscript.

\bibliographystyle{model1b-num-names} 
\bibliography{<your-bib-database>}

\appendix

\section{Computing a companion block diagonal form}
\label{appcompanion}
It is well-known that a scalar equation $[a_{\sigma_a}]$ given by~\eqref{scalarparammm} can be expressed as a first-order system by setting $$F = [f, \partial f, \partial^2 f, \dots, \partial^{n-1} f]^T.$$ However, the opposite direction of this transition is nontrivial although possible. The theoretical possibility stems from the work of~\cite{key213} with the so called cyclic vectors. Consider again system~\eqref{param} given by:
$$[A_{\sigma_A}] \quad \quad \quad \xi_A^{h} \; x^{p_A} \; \partial F\; = \; A(x, \xi_A) F ,\quad \xi_A = x^{\sigma_A} \varepsilon \; \text{and} \; A_{0,0} \neq 0.
$$
We also recall that with the definition of $\xi_A$, $A(x, \xi_A) = \sum_{k=0} A_k(x) \xi_A^k$ with $A_k(x) \in \mathcal{M}_n(\mathbb{C}[[x]])$ for all $k>0$.
In this appendix, we compute a companion block diagonal form for $[A_{\sigma_A}]$, in a finite number of steps, in a neighborhood of a turning point. This algorithm is an adaptation of that of~\cite{key17}, developed for the unperturbed counterpart of $[A_{\sigma_A}]$. It relies on a sequence of polynomial (shearing) transformations in $x, \xi_A$ and elementary row/column operations. 
\begin{itemize}
\item Shearing transformations: \begin{equation} \label{shearingopp} T = \operatorname{diag} (x^{\alpha_1} \xi_A^{\beta_1}, x^{\alpha_2} \xi_A^{\beta_2}, \dots, x^{\alpha_n} \xi_A^{\beta_n} ) \end{equation} where the $\alpha$'s and $\beta$'s are respectively rational numbers and integers. We remark that shearing transformations may alter $\sigma_A$ because of the $\alpha$'s.
\item Elementary transformations: We consider transformations of the form \begin{equation} \label{elemopp}T = I_n + P_{i,j}(a) \end{equation} where $P_{i,j}$ is a zero matrix except for entry in the $i^{th}$ row and $j^{th}$ column for $i,j \in \{ 1, \dots, n\}$ and  $a \in \rm{\mathcal{R}}$ such that $a(x=0, \xi_A=0)\neq 0$.  Obviously, $T \in GL_n(\rm{\mathcal{R}})$. The equivalent system $[\tilde{A}_{\sigma_{\tilde{A}}}] = T [A_{\sigma_A}]$ is such that $\tilde{h}=h$, $\tilde{p}=p$, and $\sigma_{\tilde{A}} \geq \operatorname{min} (\frac{p}{h} + \sigma_A, \sigma_A)$. Moreover, it can be easily verified that the effect of this transformaiton on the system $[A_{\sigma_A}]$ is as follows:
\begin{description}
\item[$M_{i}(a)$:] Multiplies the  $i^{th}$ column and $i^{th}$ row by $a$ and $1/a$ respectively. It also adds ($-\xi_A^{h} x^{p} \frac{\partial a}{a}$) to the diagonal entry in the $(i,i)^{th}$ position.
\item[$C_{i,j}(a)$:] Adds to the $j^{th}$ column the $i^{th}$ column multiplied by $a$, adds to the $i^{th}$ row the $j^{th}$ row multiplied by $-a$,  and adds ($-\xi_A^{h} x^{p} \partial a$) to the entry in the $(i,j)^{th}$ position.
\item[$R_{i,j}(-a)$:] Adds to the $i^{th}$ row the $j^{th}$ row multiplied by $a$, adds to the $j^{th}$ column the $i^{th}$ column multiplied by $-a$,  and adds ($\xi_A^{h} x^{p} \partial a$) to the entry in the $(i,j)^{th}$ position.
\end{description}
\end{itemize}
A careful choice of these transformations will allow us to establish the following result constructively:
\begin{theorem}
\label{compp}
Consider the system $[A_{\sigma_A}]$ given by~\eqref{param}. Then there exists a transformation $T$ which is a product of transformations of the forms~\eqref{shearingopp} and ~\eqref{elemopp} in a root of $x$, such that the equivalent system $\; [\tilde{A}_{\sigma_{\tilde{A}}}] \; \xi_{\tilde{A}}^{\tilde{h}} x^{\tilde{p}} \partial F = \tilde{A} F$  has the block diagonal form:
$$ \tilde{A}(x, \xi_{\tilde{A}}) = \operatorname{diag} ( \tilde{A}^{11}(x, \xi_{\tilde{A}}), \dots,  \tilde{A}^{\ell \ell}(x, \xi_{\tilde{A}}))$$
where the $ \tilde{A}_{11}, \dots,  \tilde{A}^{\ell \ell}$ are companion matrices.
\end{theorem}
The proof of this theorem relies on two lemmas which we establish herein. Before proceeding, we can always assume, without loss of generality, that $A_{0,0} = A(x=0,\xi_A=0)$ is in Frobenius canonical form, i.e.,
\begin{equation} \label{frobp} A_{0,0} = \operatorname{diag} (A_{0,0}^{11}, \dots, A_{0,0}^{\ell \ell}) \end{equation} where the block submatrices
$A_{0,0}^{11}, \dots, A_{0,0}^{\ell \ell}$ are constant matrices in companion form, of dimensions $n_1 \geq n_2 \geq \dots \geq n_\ell \geq 1$ ($n_1 + n_2 + \dots + n_\ell = n$). 
\begin{lemma}
\label{arnoldpa}
Consider the system $[A_{\xi_A}]$ given by~\eqref{param} whose leading constant coefficient $A_{0,0}$ is in Frobenius canonical form. Then, by a sequence of elementary operations of the form~\eqref{elemopp}, the equivalent system $[\tilde{A}_{\sigma_{\tilde{A}}}]\; \xi_{\tilde{A}}^{h} x^{p} \partial F = \tilde{A} F$ has the following block decomposition:
\begin{equation} \label{arnoldpa1} \tilde{A}(x, \xi_{\tilde{A}}) = \begin{bmatrix} S(x, \xi_{\tilde{A}}) & U(x, \xi_{\tilde{A}}) \\ V(x, \xi_{\tilde{A}}) & W(x, \xi_{\tilde{A}}) \end{bmatrix} , \end{equation}
where 
\begin{itemize}
\item $S(x, \xi_{\tilde{A}})$ is a $n_1$-square matrix  in companion form:
$$S(x, \xi_{\tilde{A}}) = \begin{bmatrix} 0 & 1& \quad & \quad & 0 \\ 0 & 0 & \ddots & \quad \\ & \quad & \quad & \quad & \quad \\ & \quad & \quad & \quad  & 1\\ s_1 & s_2 & \dots & \quad &s_{n_1} \end{bmatrix} ;$$ 
\item $U(x, \xi_{\tilde{A}})$, $V(x, \xi_{\tilde{A}})$ are of dimensions $n_1\times (n-n_1)$ and $(n-n_1) \times n_1$ respectively, $U_{0,0} = O_{n_1\times (n-n_1)}$, $V_{0,0}=O_{(n-n_1) \times n_1}$ , and of the forms:
$$U(x, \xi_{\tilde{A}}) = \begin{bmatrix} 0 & \hdots  & 0 \\  \vdots & \quad  & \vdots \\ 0 & \hdots  &0\\ u_1  & \dots  &u_{n-n_1} \end{bmatrix} \quad \quad \text{and} \quad \quad V(x, \xi_{\tilde{A}}) = \begin{bmatrix} v_1 & 0 & \hdots & 0 \\  \vdots & \vdots & \quad & \vdots \\  v_{n-n_1} & 0 & \dots  &0 \end{bmatrix};$$
\item $W(x, \xi_{\tilde{A}})$ is a square matrix of dimension $n-n_1$;
\item Furthermore, $\tilde{A}_{0,0}=A_{0,0}$.
\end{itemize}
\end{lemma}
\begin{proof}
We first assume $n_1 >1$. Otherwise, $A(x, \xi_A) = [a_{ij}]_{1 \leq i, j \leq n}$ is already in form~\eqref{arnoldpa1}. Since $A_{00}$ has the form~\eqref{frobp}, $a_{i, i+1}(x=0, \xi=0) =1$. Hence,  $1/a_{i, i+1}(x, \xi_A) \in \rm{\mathcal{R}}$ and $1/a_{i, i+1}(x=0, \xi_A=0) =1$. 
\begin{description}
\item[Step 1 ] For each row number $i$ from 1 to $n_1-1$:
\begin{description}
\item[(1.1)] We use $1/a_{i, i+1}(x, \xi_A)$ as a pivot and apply $M_{i+1}(1/a_{i, i+1}(x, \xi_A))$ to set the entries in the $(i, i+1)$ positions to $1$. We remark that the term  ($-\xi_A^{h} x^{p} (a) \partial(1/a)$) is added to the diagonal entry in the $(i+1,i+1)^{th}$ position. 
\item[(1.2)] We now make use of these $1$'s to set the entries in the $(i,j)$ positions for $i\in \{1, \dots, n_1 -1 \}$, $j \in \{1, \dots, n\}$, $j \neq i+1$ to zeroes: For each row number $j \neq i+1$ from 1 to $n$, we apply $C_{i+1, j}(-a_{ij})$. We remark that the term  ($+\xi_A^{h} x^{p} \partial a$) is added to the entry in the $(i+1,j)$ position. Clearly, the rows $1$ to $i-1$ are not altered by these elementary operations.
\end{description}
\item[Step 2 ] We have attained the anticipated form for the first $n_1$ rows, i.e. we have shown how to construct an equivalent system whose first $n_1$ rows have the described properties of $\begin{bmatrix} S(x, \xi_A) & U(x, \xi_A) \end{bmatrix}$. If $n_1 < n$ then we can proceed to work on the last $n-n_1$ rows. 
\begin{description}
\item[(2.1)] Again, we make use of these $1$'s created in Step $[(1.1)]$ to set the entries in the $(i,j)$ positions for $i\in \{n_1 +1, \dots, n \}$ and $j \in \{2, \dots, n_1\}$ to zeroes: For each $i, j$ in these ranges, we apply $R_{i, j-1}(-a_{ij})$ successively for each $i$ from $n_1 +1$ to $n$ with an inner loop of $j$ from $n_1$ to $2$. \marginpar{\color{blue} the inner loop shld start from $n_1$ because the term coming from the derivation settles in the $(i,j-1)$ position and so will be deleted afterwards. I should specify how $\sigma_{\tilde{A}}$ changes w.r.t. $\sigma_A$.} Evidently,  each such operation does not alter the form created by Step $1$ for the first $n_1$ rows,  and by its preceeding operations within Step $2$.
\end{description}
\end{description}
Finally, since each of the operations performed is of the form~\eqref{elemopp}, for $x=0$ and $\xi_A=0$, such an operation reduces to $I_n$. Hence, the form of $A_{0,0}$ is preserved.
\end{proof}
\begin{remark}
One can observe that by a repetitve application of Lemma~\ref{arnoldpa} on the diagonal blocks of dimensions $n_1 \geq n_2 \geq \dots \geq n_{\ell}$, we arrive at the aformentioned Arnold form.
\end{remark}
\begin{lemma}
\label{arnoldpb}
Consider the system $[A_{\sigma_A}]$ given by~\eqref{param}. Suppose that its leading constant coefficient $A_{0,0}$ is in Frobenius canonical form and $A(x, \xi_A)$ is in form~\eqref{arnoldpa1}. Then, by a sequence of constant transformations, elementary operations of the form~\eqref{elemopp}, and shearings of the form~\eqref{shearingopp}, the equivalent system $\; [\tilde{A}_{\sigma_{\tilde{A}}}] \; \xi_{\tilde{A}}^{h} x^{\tilde{p}} \partial F = \tilde{A} F$ is such that has the following block diagonal form:
\begin{equation} \label{arnoldpb1} \tilde{A}(x, \xi_{\tilde{A}}) = \begin{bmatrix} \tilde{S}(x, \xi_{\tilde{A}}) & O \\ O & \tilde{W}(x, \xi_{\tilde{A}}) \end{bmatrix} , \end{equation}
where $\tilde{S}(x, \xi_{\tilde{A}})$ is a $r$-square companion matrix where $r \geq n_1$ and $ \tilde{W}(x, \xi_{\tilde{A}})$ is a $n-r$ square matrix.
\end{lemma}

\begin{proof}
If $U(x, \xi_A)$ and $V(x,\xi_A)$ of~\eqref{arnoldpa1} are smiltaneously zero matrices or $n_1 = n$ then $A(x, \xi_A)$ is already in the form~\eqref{arnoldpb1}. Otherwise we suppose that $U(x, \xi_A)$ is nonzero and limit our discussion to this case. If $V(x, \xi_A)$ is nonzero as well, then we can proceed analogously. Let $u_j = \sum_{k=0}^{\infty} u_{j,k}(x) \xi_A^{k}$ for $j \in \{1, \dots, n-n_1 \}$ then we define
$$\beta = \operatorname{min}_{1 \leq j \leq n-n_1} (\;val_{\varepsilon}\;(u_j)) \quad \text{and} \quad \alpha = \operatorname{min}_{1 \leq j \leq n-n_1} (\;val_{x}\;(u_{j, \beta})).$$
We remark that $\beta \geq 0$, and if $\beta =0$ then $\alpha >0$ since $U_{0,0}$ is a zero matrix. Let $$T = \operatorname{diag} (x^{\alpha} \xi_A^{\beta} I_{n_1}, I_{n -n_1}) \quad \text{and} \quad  [B_{\sigma_B}] = T[A_{\sigma_A}].$$ Then we have:
$$ B(x, \xi_B) = \begin{bmatrix} S(x, \xi_A) - x^{p+\alpha -1} \xi_A^{\beta + h} (\alpha + \beta \sigma_A) I_{n_1} & x^{-\alpha} \xi_A^{-\beta} U(x, \xi_A) \\ x^{\alpha} \xi_A^{\beta} V(x, \xi_A) & W(x, \xi_A) \end{bmatrix} .$$
$$B_0(x) = \begin{bmatrix} S_0(x) & x^{-\alpha} \xi_A^{-\beta} U_0(x) \\ O & W_0(x)\end{bmatrix} $$
and
$$B_{0,0} = \begin{bmatrix} S_{0,0}& U_{0,0} \\ O & W_{0,0} \end{bmatrix} $$
where $S_{0,0}$ is a companion matrix of dimension $n_1$, $U_{0,0}$ is a companion matrix whose last row has at least one nonzero element, and $W_{0,0}$ is a block diagonal square matrix whose blocks are in companion form. Let $e_i$ denote the $i^{th}$ row of the identity matrix $I_n$ and $B_{0,0}^i$ denote the $i^{th}$ power of $B_{0,0}$. W then have:
\begin{eqnarray*}
e_1 B_{0,0}^i &=& \sum_{j = 1}^i c_{i,j} e_j + e_{i+1} \quad \text{for} \quad i=1, \dots, n_1 -1 \\
e_1 B_{0,0}^{n_1} &=& \sum_{j = 1}^{n_1} c_{n_1,j} e_j + \sum_{j = 1}^{n-n_1} u{j,\beta} e_{{n_1}+j}
\end{eqnarray*}
where $c_{i,j} \in \mathbb{C}$. Since there exists at least one $j \in \{1, \dots, n-n_1\}$ such that $u_{j,\beta} \neq 0$, then the vectors $e_1, e_1 B_{0,0}, e_1B_{0,0}^2, e_1 B_{0,0}^{n_1}$ are linearly independent. Hence, denoting by $r$ the dimension of the first block of the Frobenius form of $B_{0,0}$, i.e. the degree of the minimal polynomial of $B_{0,0}$, we have $r>n_1$. Hence, by putting the leading constant matrix in Frobenius form, applying Lemma~\eqref{arnoldpa}, and the procedure described herein, we either arrive at the desired form in a finite number of steps.
\end{proof}

\begin{proof}(Theorem~\ref{compp})
The proof can be attained by a recursive application of  Lemmas~\ref{arnoldpa} and~\ref{arnoldpb}. 
\end{proof}
\begin{remark}
We remark that this method can give a cyclic vector and that $[\tilde{A}_{\sigma_{\tilde{A}}}]$ is not uniquely determined by $[A_{\sigma_A}]$. 
\end{remark}

\section{Proofs for Section~\ref{blockparam}}
\label{appb}
\begin{proof}(Proposition~\ref{newtoneqnxi})
Let $f(x, \xi) = exp (\int q(x, \xi) dx)$ with $q(x, \xi) \in \bigcup_{s \in \mathbb{N}^{*}} \overline{\mathbb{C}((x))} ((\varepsilon^{1/s}))$. 
We remark first that 
$$\forall i \in \mathbb{N} , \quad \partial^i f = P_i(q) f$$
where 
$$P_i(q) = q^i + \dots + \partial^{i-1}_x q$$ is a polynomial in $(q, \partial q, \dots, \partial^{i-1} q)$.
Then, $f(x, \xi)$ verifies~\eqref{scalarparammm} if and only if $q(x, \xi)$ is a solution of 
$$ \sum_{i=0}^n a_{i}(x, \xi) P_i(q) .$$
We denote by $e$ the polar order of $q$ in $\xi$  ($\varepsilon$-order), and write
$$q(x, \xi) = \xi^e (s(x) + O(\xi)), q \not\equiv 0 .$$
Then, $\forall i \in \{ 0, \dots, n\}, val_{\varepsilon}(P_i)= -i e$ and in fact $P_i(q(x, \xi))  = \xi^{-id} (s(x)^i + O(\xi))$.

Thus,$$ \sum_{i=0}^n a_i(x, \xi) P_i(q) = \sum_{i=0}^n \xi^{\nu_i - i d} W_i(q), $$ where $val_{\varepsilon} (W_i) \geq 0$. Let $v(e) = min_{0 \leq i \leq n} (\nu_i - i e)$ and $I(e)= \{ i \in \{ 0, \dots, n\} \;|\; val_{\varepsilon}(a_i - i e) = v(e) \}$. Hence,
$$ \sum_{i=0}^n a_{i}(x, \xi) P_i(q) = \xi^{v(e)} ( \sum_{i \in I(e)} a_{i, \nu_i} (x) s^i(x) + O(\xi) ).$$
Thus, $q(x, \xi) = \xi^e (s(x) + O(\xi)), q \not\equiv 0 $ satisfies $\sum_{i=0}^n a_{i}(x, \xi) P_i(q) \equiv 0$ if and only if $s(x)$ is a non zero solution of the equation
$$ \sum_{i \in I(e)} a_{i, \nu_i} (x) X^i (x) =0 .$$
This equation has non trivial solutions if and only if $|I(e)| > 1$, i.e. if and only if $e \in \{ e_1, \dots, e_\ell \}$. And then, $X(x)$ is one of the nonzero roots of the associated determinant equation.
\end{proof}
\begin{proof} (Theorem~\ref{splitparam}) 
We proceed in two steps:\\
\textbf{Step 1}: We first block-diagonalize the leading matrix coefficient: Suppose that there exists $T_0(x) \in Gl_n(\mathbb{C}[[x]])$ s.t. the substitution $F = T_0(x)G$ in $[A_{\sigma_A}]$ yields an equivalent system $[\tilde{A}_{\sigma_{A}}] \quad \xi_A^h x^p \partial  G = \tilde{A}(x,\xi_A) G$ where $\tilde{A}(x,\xi_A) = \sum_{k=0}^\infty \tilde{A}_k(x) \xi_A^k$, $\tilde{A}_k(x) \in \mathcal{M}_n(\mathbb{C}[[x]])$ for all $k \geq 0$, and the leading coefficient matrix $\tilde{A}_0(x)$ is non-zero and has the following block-form in accordance with $A_{0,0}$:
$$\tilde{A}_0(x) = \begin{bmatrix} \tilde{A}_0^{11}(x) & O \\ O &  \tilde{A}_{0}^{22}(x) \end{bmatrix}.$$
Then we have:
$$A \;T_0(x)\; -\; T_0(x)\; \tilde{A} = \xi_A^h x^p \partial T_0(x)$$
which yields
\begin{equation} 
\label{splitparam1k} 
A_k(x) \;T_0(x)\; -\; T_0(x)\; \tilde{A}_k(x) = \delta_{kh} \partial T_0(x), \quad \text{where} \; \delta_{kh}=1 \; \text{if}\; k=h \; \text{and}\; \delta_{kh}=0 \; \text{otherwise}. 
\end{equation}
In particular (since $h>0$):
\begin{equation}
\label{splitgauge}
A_0(x) \;T_0(x)\; -\; T_0(x)\; \tilde{A}_0(x) \;=\;0. 
\end{equation}
We further assume that $T_0(x) \in \mathcal{M}_n(\mathbb{C}[[x]])$ is in the following form
$$T_0(x) = \begin{bmatrix} I & T_0^{12}(x) \\ T_0^{21}(x) & I \end{bmatrix}.$$
Then, for $1 \leq \varrho \neq \varsigma \leq 2$, we have:
\begin{equation}
\label{splitparam1}
\begin{cases}
A_0^{\varrho \varrho}(x) \;-\; \tilde{A}_0^{ \varrho \varrho}(x) \;+\; A_0^{ \varrho \varsigma }(x)\; T_0^{\varsigma   \varrho}(x)\;= \;O \\
A_0^{\varsigma  \varrho}(x) \;+\; A_0^{\varsigma \varsigma}(x)\; T_0^{\varsigma  \varrho}(x) \;- \;T_0^{\varsigma \varrho}(x) \; \tilde{A}_0^{ \varrho  \varrho}(x) \;=\; O
\end{cases}
\end{equation}
Inserting the series expansions $A_0(x) = \sum_{i=0}^{\infty} A_{i,0} x^i$ and $\tilde{A}_0(x) = \sum_{i=0}^{\infty} \tilde{A}_{i,0} x^i$ in \eqref{splitparam1}, and equating the power-like coefficients, we get for $i =0$: $$ \begin{cases} A^{\varrho \varrho}_{0,0}\;=\;\tilde{A}^{ \varrho \varrho}_{0,0}\\ A^{\varsigma \varsigma}_{0,0}\; T^{\varsigma  \varrho}_{0,0} - T^{\varsigma \varrho}_{0,0} \; \tilde{A}^{ \varrho  \varrho}_{0,0} = O  \end{cases}$$ which are satisfied by setting $ \tilde{A}^{\varsigma \varsigma}_{0,0} = A^{\varsigma \varsigma}_{00}$ and $T^{\varsigma \varsigma}_{0,0} = I, \; T^{\varsigma \varrho}_{0,0} = O .$
And for $i \geq 1$ we get:
\begin{eqnarray}  \label{splitparam4}    A_{0,0}^{\varsigma \varsigma } T_{i,0}^{\varsigma \varrho} - T_{i,0}^{ \varsigma \varrho} A_{0,0}^{\varrho \varrho} &=& - A_{i,0}^{\varsigma \varrho} - \sum_{j=1}^{i-1} ( A^{\varsigma \varsigma}_{i-j,0} T^{\varsigma \varrho}_{j,0}  -  T^{\varsigma \varrho}_{j,0} \tilde{A}^{\varrho \varrho}_{i-j,0} )\\
  \label{splitparam3} 
{\tilde{A}}^{\varrho \varrho}_{i,0} &= & A_{i,0}^{\varrho \varrho} + \sum_{j=1}^{i-1} A^{\varrho \varsigma}_{i-j,0} T^{\varsigma \varrho}_{j,0} .
\end{eqnarray}
It's clear that \eqref{splitparam4} is a set of Sylvester matrix equations that possess a unique solution due to the assumption on the disjoint spectra of $A_{0,0}^{11}$ and $A_{0,0}^{22}$ (see, e.g., \cite[Appendix A.1, p. 212-213]{key6}). Remarking that the right hand side depends solely on the ${A}_{j,0}$, $T_{j,0}$, $\tilde{A}_{j,0}$ with $j<i$, equations \eqref{splitparam4} and \eqref{splitparam3} 
are successively soluble. Hence, such $T_0(x) \in Gl_n(\mathbb{C}[[x]])$ can be constructed. \\
Clearly, $p_T=\nu_T=0$ and $\sigma_T \geq \sigma_A$. Moreover, it follows from~\eqref{splitparam1k} that $\nu_{\tilde{A}}=0$, $p_{\tilde{A}}=0$, and $\sigma_{\tilde{A}} \geq \sigma_A$.\\
\textbf{Step 2}: Due to the first step, we can assume without loss of generality that the leading coefficient  matrix $A_0(x)$ of the system $[A_{\sigma_A}]$  has a block-diagonal form in accordance with that of $A_{0,0}$.
It suffices to seek $T$ and $\tilde{A}$ in the form $T = \sum_{k=0}^\infty T_k(x) \xi$ and $\tilde{A} = \sum_{k=0}^\infty \tilde{A}_k(x) \xi$ respectively where $T_k(x), \tilde{A}_k(x) \in GL_n(\mathbb{C}[[x]])$ for all $k>0$.
We set $T_0(x)=I_n$, $\tilde{A}_0(x)= A_0(x)$, and we rewrite $$A = \sum_{k=0}^\infty A_k(x) \xi_A = \sum_{k=0}^\infty B_k(x) \xi, \quad \text{where} \; B_k(x) = x^{(\sigma_A-\sigma) k} A_k(x).$$
Since $\sigma_A \geq \sigma$, it follows that $B_k(x) \in \mathcal{M}_n(\mathbb{C}[[x]])$ for all $k>0$. With the expansions in $\xi$, we get the following set of ansatz for $k \geq 0$:
 \begin{equation} \label{gaugeparam1} \sum_{i=0}^k ( B_{k-i}(x) T_i(x) - T_i(x) {\tilde{A}}_{k-i}(x) ) = x^p (\partial + \frac{\sigma (k-h)}{x}) T_{k-h} ,\quad \text{where} \quad T_{k-h} = 0 \;\text{for}\; k<h. \end{equation}
It then follows from~\eqref{gaugeparam1} that  for $k>1$ we have:
\begin{eqnarray*}  B_{0}(x) T_k(x) - T_k(x) {\tilde{A}}_{0}(x) &=&\sum_{i=1}^{k-1} B_{k-i}(x) T_i(x) - T_i(x) {\tilde{A}}_{k-i}(x) \\ &+& x^p (\partial + \frac{\sigma (k-h)}{x}) T_{k-h}(x) - {\tilde{A}}_{k}(x) + B_{k}(x) . \end{eqnarray*}
 By inserting the claimed block forms of $T_k(x)$ and ${\tilde{A}}_{k}(x)$ for all $k>1$ in the above equation, with $1 \leq \varrho \neq \varsigma \leq 2$, we obtain:
\begin{eqnarray*}
\tilde{A}_k^{ \varrho \varrho}(x) &=& B_k^{ \varrho \varrho }(x) \\
B_0^{\varsigma \varsigma}(x)\; T_k^{\varsigma  \varrho}(x) \;- \;T_k^{\varsigma \varrho}(x) \; \tilde{A}_0^{ \varrho  \varrho}(x) &=& \sum_{i=1}^{k-1} B_{k-i}^{\varsigma \varsigma}(x)\; T_i^{\varsigma  \varrho}(x) \;- \;T_{k-i}^{\varsigma \varrho}(x) \; \tilde{A}_i^{ \varrho  \varrho}(x) \\ &+& x^p (\partial + \frac{\sigma (k-h)}{x}) T_{k-h}^{\varsigma  \varrho}(x) + B_k^{\varsigma  \varrho}(x) .
\end{eqnarray*}
Hence, for every $k>1$, the computation of $T_k(x)$ requires only lower order terms of $T, B,$ and $ \tilde{A}$. Hence, by setting $T_k(x=0) =0$ and expanding each of the equations for a fixed $k>1$ w.r.t. $x$, we arrive at constant Sylvester matrix equations, which can be resolved successively due to the disjoint spectrum of $A_{0,0}^{11}$ and $A_{0,0}^{22}$, to  compute $T_k(x)$ up to any desired order in $x$. 
\end{proof}
\section{Proofs for Section~\ref{moser}: $\varepsilon$-rank reduction }
\label{appc}
In this section, we establish our constructive proof of the necessity of Theorem~\ref{mosernilpotent1}. For sufficiency, one can consult~\cite[Section 5.5, p. 101]{key1000}. We prove the following:

\begin{theorem}
\label{mosernilpotent2}
Consider the system~\eqref{param} given by:
$$
[A_{\sigma_A}] \quad \quad \quad x^{p} \xi_A^{h} \partial F \; = A(x, \xi_A) F =  \sum_{k=0}^{\infty} A_k(x) {\xi_A}^k  \; F
$$
with $h>1$ and $m_{\varepsilon}(A)>1$ and set  $rank\;(A_0(x))=r$ . Suppose that the polynomial 
$$\theta_A(\lambda) := {\xi_A}^{r}  \det (\lambda I + \frac{A_0(x)}{\xi} + A_1(x) )|_{\xi_A=0}$$ 
vanishes identically in $\lambda$. Then there exists a transformation $F = R(x, \xi_R) G$ such that the equivalent system
$$
[\tilde{A}_{\sigma_{\tilde{A}}}] \quad \quad \quad {\xi_{\tilde{A}}}^{\tilde{h}} x^{\tilde{p}} \partial G \; =\; \tilde{A} (x, \xi_{\tilde{A}}) \; G = \; \sum_{k=0}^{\infty} \tilde{A}_k(x) \xi_{\tilde{A}}^k  \; G
$$
satisfies either $\tilde{h} < h$ or $\tilde{h} = h$ and $rank (\tilde{A}_0(x)) < r$.
Moreover, such a $R(x, \xi_R)$ can be constructed by \cite[Algorithm 1]{key102} where it is always chosen to be a product of unimodular transformations in $GL_n(\mathbb{C}[[x]])$ and polynomial transformations in $\xi_A$. 
\end{theorem}
The proof of the theorem will be given after a set of intermediate results. 
\begin{lemma}
\label{gauss1}
Given the system 
 $$
[A_{\sigma_A}] \quad \quad \quad x^p \xi_A^h \; \partial F=  A(x, \xi_A) F = \sum_{k=0}^{\infty} A_k(x) \xi_A^k \; F ,  $$ with $r= rank\;(A_0(x))$. There exists a unimodular transformation $U(x)$ in $GL_n(\mathbb{C}[[x]])$ such that the leading coefficient matrix $\tilde{A}_0(x)$ of the equivalent system~\eqref{gaugep} given by 
$$
[\tilde{A}_{\sigma_{\tilde{A}}}] \quad \quad \quad \xi_{\tilde{A}}^h x^p \partial G = \tilde{A}(x, \xi_{\tilde{A}}) G ,\quad \sigma_{\tilde{A}} \geq \operatorname{min} (\frac{p}{h}+\sigma_A, \sigma_A) $$
has the form
$$\tilde{A}_0(x)  = \begin{bmatrix} \tilde{A}_0^{11} (x) & O \\ \tilde{A}_0^{21} (x) & O \end{bmatrix} $$
where  $\tilde{A}_0^{11} (x) $ is a square matrix of dimension $r$ and $\begin{bmatrix} \tilde{A}_0^{11}(x) \\ \tilde{A}_0^{21}(x) \end{bmatrix}$ is a $n \times r$ matrix of full column rank $r$. 
\end{lemma}
\begin{remark}
\label{gauss2}
In practice, $U(x)$  can be obtained by performing Gaussian elimination on the columns of $A_0(x)$ taking as  pivots  the elements of minimum valuation (order) in $x$.
\end{remark}
\begin{proof}
By Remark \ref{noeffect1}, $\tilde{A}_0(x) = U^{-1}(x) A_0(x) U(x)$. Hence it suffices to search a similarity transformation $U(x)$. Since $\mathbb{C}[[x]]$ is a principal ideal domain (the ideals of $\mathbb{C}[[x]]$ are of the form $x^k \mathbb{C}[[x]]$), it is well known that one can construct unimodular transformations $Q(x), U(x)$ lying in $GL_n(\mathbb{C}[[x]])$ such that the matrix $Q(x) A_0(x) U(x)$ has the Smith normal form
$$Q(x)\;A_0(x)\;U(x) = \operatorname{diag} (x^{\beta_1}, \dots, x^{\beta_r}, 0, \dots, 0) )$$
where $\det (U(0)) \ne 0$, $\det (Q(0)) \ne 0$ , and $\beta_1, \dots, \beta_r$ in $\mathbb{Z}$ with $0 \leq \beta_1 \leq \beta_2 \leq \dots \leq \beta_r$.\\
It follows that we can compute a unimodular matrix $U(x)$ in $GL_n(\mathbb{C}[[x]])$ so that its last $n-r$ columns form a $\mathbb{C}[[x]]$-basis of ker ($A_0(x)$).  \\
As for finding a lower bound of $\sigma_{\tilde{A}}$, one can see that since $U(x=0)$ is invertible,  we have:
$$val_x (\tilde{A}_k) \geq \begin{cases} val_x (\tilde{A}_k) \quad \text{if} \; k \neq h \\
val_x (\tilde{A}_h) + p \quad \text{otherwise}. \end{cases}$$
Thus, if $p<0$ then $\sigma_{\tilde{A}} \geq \sigma_A +\frac{p}{h}$.
\end{proof}
\begin{remark}
\label{noeffect1} Consider system $[A_{\sigma_A}]$ given by~\eqref{param}.
Let $T(x)= ( I_n + Q(x)) \in GL_n(\mathbb{C}[[x]])$ such that
$$ Q(x) = [q_{ij}]_{1 \leq i, j \leq n} \; , \quad s.t.\; \begin{cases} q_{nj} \in \mathbb{C}[[x]] \;,\quad \text{for} \; r+1 \leq j < n \\ q_{ij} = 0 \quad \text{elsewhere} .\end{cases}$$ 
One can then verify that the resulting equivalent system \eqref{gaugep} is 
$$[\tilde{A}_{\sigma_{\tilde{A}}}] \quad \quad \quad \xi_{\tilde{A}}^{h} x^{p} \partial G = \tilde{A}(x, \xi_{\tilde{A}}) G, \quad \sigma_{\tilde{A}} \geq \sigma_A$$
\begin{equation} \label{gauss5} \text{where} \quad \begin{cases}  \tilde{A}_0(x) =  ( I_n + Q(x))^{-1} \;A_0(x)\;  ( I_n + Q(x)) \\
 \tilde{A}_1(x)  =  ( I_n + Q(x))^{-1} \; A_1(x) \; ( I_n + Q(x))  - \gamma x^p \partial Q(x) ; \\
 \tilde{A}_2(x)  = \dots \end{cases} \end{equation}
where $\gamma =1$ if $h=1$ and $\gamma =0$ otherwise. 
But, due to \eqref{thetaparam}, we can limit our interest to $A_0(x)$, $A_1(x)$, $ \tilde{A}_0(x)$, and $ \tilde{A}_1(x)$ solely. Hence, if $h>1$ then it suffices to investigate $T^{-1} \; A(x, \xi_{\sigma_A}) \; T$.  We remark that if $h=1$ then the term  $T^{-1} \partial T(x)$ should be taken into account. This is the reason the case $h=1$ is treated separately. This distinction is purely technical as explained in the following remark and unlike for the Poincar\'e rank $p$ of an unperturbed system, it does not lead to a classification of regularity.
\end{remark}
 By Lemma~\eqref{gauss1}, we can suppose without loss of generality that  $A_0(x)$ is of the form \begin{equation} \label{gaussform} A_0(x)  = \begin{bmatrix} A_0^{11}(x)  & O \\ A_0^{21}(x) & O \end{bmatrix} \end{equation} with $r$ independent columns and $(n-r)$ zero columns. We partition $A_1(x)$ in accordance with $A_0(x)$, i.e. $ A_1(x)= \begin{bmatrix} A_1^{11} (x)& A_1^{12}(x) \\ A_1^{21}(x) & A_1^{22} (x) \end{bmatrix},$  and consider
\begin{equation} \label{glambdaform} G_{A}(\lambda)= \begin{bmatrix} A_0^{11}(x) & A_1^{12} (x)\\ A_0^{21}(x) & A_1^{22}(x) + \lambda  I_{n-r}\end{bmatrix} .\end{equation}
This consideration of $G_{A}(\lambda)$ gives an $\varepsilon$-reduction criterion equivalent to $\theta_{A}(\lambda)$. In fact, let $D(\xi_A)= \operatorname{diag} (\xi_A I_{r}, I_{n-r})$ where $r= rank\;(A_0(x))$. Then we can write $\xi_A^{-1} A(x,\xi_A) = N D^{-1}$ where $N:= N(x, \xi_A) \in \mathcal{M}_n(\rm{\mathcal{R}})$. Set $D_0 = D(0)$ and $N_0= N(x,0)$. Then we have 
\begin{eqnarray*} \det (G_{A}(\lambda)) & =& \det (N_0 + \lambda D_0)  =  \det (N + \lambda D)|_{\xi_A=0} \\
& = & ( \det (\frac{A(x,\xi_A)}{\xi_A} + \lambda I_n) \det (D))|_{\xi_A=0}  \\ &=&  (\det (\frac{A_0(x)}{\xi_A} + A_1(x) + \lambda I_n) \; \xi_A^{r} )|_{\xi_A=0} = \theta_A (\lambda).  \end{eqnarray*}  
Thus,
$\det (G_{A}(\lambda)) \equiv 0$ vanishes identically in $\lambda$ if and only if $ \theta_A(\lambda)$ does. 
We illustrate our progress with the following simple example. 
\begin{example}
\label{firststep}
~\cite[Example 2]{key102} Consider $\quad [A_{\sigma_A}]\quad \xi_A^h \partial F = A(x, \xi) F $ with $\sigma_A=0$, $h >1$, and $$A(x, \xi_A)=  \begin{bmatrix} \xi_A & -x^3 \xi_A & (1+x) \xi_A & 0 \\ x^2 & x \xi_A & 0 & -2x \xi_A \\ -x & 0 & 0 & 2 \xi_A \\ 0 & 2 & 0 & \xi_A^2\end{bmatrix} . $$ Clearly, $A_0(x)$ is nilpotent of rank $2$ and \begin{equation} \label{firststepg} G_{A}(\lambda)= \begin{bmatrix}  0 & 0 & x+1  & 0 \\ x^2 & 0 & 0 & -2 x\\ -x & 0 &  \lambda  & 2 \\ 0 & 2 & 0 &  \lambda  \end{bmatrix} .\end{equation}
\end{example}

 We then have the following proposition:
\begin{proposition}
\label{gauss3}
Given the system
$$
[A_{\sigma_A}] \quad \quad \quad \xi_A^h \; x^p \; \partial F=  A(x, \xi_A) F = \sum_{k=0}^{\infty} A_k(x) \xi_A^k \; F ,  $$ 
with $r = \;rank\;(A_0(x))$. If $m_{\varepsilon}(A) >1$ and $\det (G_A (\lambda)) \equiv 0$ then there exists a finite product of triangular matrices $T(x) = P ( I_n + Q(x)) $ where $P$ is a permutation and $\det T(x) = \pm 1$, such that for the equivalent system $$
[\tilde{A}_{\sigma_{\tilde{A}}}] \quad \quad \quad\xi_{\tilde{A}}^{h} x^{p} \partial G \; =\; \tilde{A} (x, \xi_{\tilde{A}}) \; G = \; \sum_{k=0}^{\infty} \tilde{A}_k(x) \xi_{\tilde{A}}^k  \; G
$$
we have:
\begin{equation} \label{particularform3} G_{\tilde{A}} (\lambda) = \begin{bmatrix} A_0^{11}(x) & U_1(x) & U_2(x) \\ V_1(x) & W_1(x) + \lambda I_{n- r -\varrho} & W_2(x) \\ M_1(x) & M_2(x) & W_3(x) + \lambda I_\varrho \end{bmatrix} ,\end{equation}
where $0 \leq \varrho \leq n-r,\; W_1,\; W_3$ are square matrices of order $(n-r-\varrho)$ and $\varrho$ respectively, and
\begin{equation} \label{particularconditionb} rank\;(\begin{bmatrix} A_0^{11}(x) & U_1(x) \\ M_1(x) & M_2(x) \end{bmatrix}) = rank\;(\begin{bmatrix} A_0^{11}(x) & U_1(x) \end{bmatrix}),\end{equation}
\begin{equation}  \label{particularconditiona} rank\;(\begin{bmatrix} A_0^{11}(x) & U_1(x) \end{bmatrix}) < r . \end{equation}
Moreover, $\sigma_{\tilde{A}} \geq \operatorname{min} (\frac{p}{h}+\sigma_A, \sigma_A)$.
\end{proposition}
We shall need the following Remark in the proof of the Proposition. 
\begin{remark}\cite[Remark 6]{key102}
\label{gq}
Suppose that $G_{A}(\lambda)$ has the form \eqref{particularform3} and there exists a transformation $T(x) \in GL_n(\mathbb{C}((x)))$ such that $G_{T[A_{\sigma_A}]}(\lambda)$ has the form
$$ G_{T[A_{\sigma_A}]}(\lambda)  = \begin{bmatrix} A_0^{11} & U_1 & U_2 \\ V_1 & W_1 + \lambda I_{n- r -\varrho} & W_2 \\ O & O & \tilde{W}_3 + \lambda I_\varrho \end{bmatrix} ,$$
where $0 \leq \varrho \leq n-r\; $ and $\tilde{W}_3$ is upper triangular with zero diagonal. Then, $$\det (G_{T[A_{\sigma_A}]}(\lambda)) = {\lambda}^\varrho \det ( \begin{bmatrix} A_0^{11} & U_1 \\ V_1 & W_1 + \lambda I_{n- r -\varrho} \end{bmatrix}).$$
If $\det (G_{T[A_{\sigma_A}]}(\lambda)) \equiv 0$ then we have 
 $\det (G_{T[A_{\sigma_A}]}(\lambda)) \equiv 0 $ as well (rank of leading coefficient matrix  is unchanged). Hence,
\begin{equation}
\label{singular}
\det ( \begin{bmatrix} A_0^{11} & U_1  \\ V_1 & W_1 + \lambda I_{n- r -\varrho} \end{bmatrix}) \equiv 0 .
\end{equation}
For a fixed $\varrho \in \{ 0, \dots, n-r \}$ we shall denote by $G_0^{(\varrho)}$ the matrix $ \begin{bmatrix} A_0^{11} & U_1  \\ V_1 & W_1 \end{bmatrix}$ of~ \eqref{particularform3}.
\end{remark}
\begin{proof} (Proposition ~\ref{gauss3})
Since $\det (G_A(\lambda))  \equiv 0$ then in particular, the matrix $G_A(\lambda=0)$ is singular. Let $E_1$ (respectively $E_2$) be the vector space spanned by the first $r$ (resp. last $n-r$) rows of $G_A(\lambda=0)$. We have
$$dim (E_1+E_2) = rank\;(G_A(\lambda=0)) < n .$$
If  $dim (E_1) < r$ then setting $\varrho =0$ suffices to fulfill our claim. Otherwise, since  $$dim(E_1 + E_2)  = dim(E_1) + dim(E_2) - dim(E_1 \cap E_2) < n , $$   it follows that 
 either $dim(E_2) < n-r$ or $dim(E_1 \cap E_2) > 0$. In both cases, there exists at least a  row vector $\varpi^{(1)}(x)=(\varpi^{(1)}_1 (x), \dots,\varpi^{(1)}_n(x))$ with entries in $\mathbb{C}((x))$ in the left null space of $G_A(\lambda=0)$, such that $\varpi^{(1)}_j(x) \neq 0$ for some $r+1 \leq j \leq n$. We can assume without loss of generality that $\varpi^{(1)}(x)$ has its entries in $\mathbb{C}[[x]]$. Indeed, this assumption can be guaranteed by a construction as in Remark~\ref{gauss2}. Let the constant matrix $P^{(1)}$ denote the product of permutation matrices which exchange the rows of $A$, so that $val_x (\varpi^{(1)}_n (x)) <  val_x (\varpi^{(1)}_j (x))\;, r+1 \leq j \leq n-1$, where $val_x$ denotes the $x$-adic valuation (order in $x$). Let 
$$Q^{(1)}(x) = [q^{(1)}_{ij}(x)]_{1 \leq i, j \leq n} \; , \quad s.t.\; \begin{cases} q^{(1)}_{nj}(x) = - \frac{\varpi^{(1)}_j(x)}{\varpi^{(1)}_n (x)}\;,\quad \text{for} \; r+1 \leq j < n  \\ q^{(1)}_{ij} = 0 \quad \text{elsewhere} \end{cases}$$ 
Thus, $P^{(1)}( I_n + Q^{(1)}(x))$ is unimodular in $GL_n(\mathbb{C}[[x]])$. \\
Set $F= F^0$, $A= A^{(0)}$, $\tilde{A}= A^{(\varrho)}$ and let $A^{(1)}$ be the matrix of the equivalent system $\xi_{A^{(1)}}^{h} x^p \partial F^{(1)} = A^{(1)}(x, \xi_{A^{(1)}}) F^{(1)}$ obtained by the transformation 
$$F^{(0)} = P^{(1)} ( I_n + Q^{(1)}(x))\; F^{(1)}.$$
Thus, by Remark \eqref{noeffect1}, $G_{A^{(1)}}(\lambda)$ has the form \eqref{particularform3} with \eqref{particularconditionb} and $\varrho=1$. \\
By Remark \ref{gq}, the matrix $G_{A^{(1)}}(\lambda = 0)$ is singular and the condition \eqref{particularconditiona} does not occur, then one can find, by the same argument as above a permutation matrix and a nozero vector $\varpi^{(2)}(x)$ in the left null space of $G^{(1)}(\lambda = 0)$.  Let 
$$Q^{(2)}(x) = [q^{(2)}_{ij}(x)]_{1 \leq i, j \leq n} \; , \quad s.t.\; \begin{cases} q^{(2)}_{n-1, j}  (x)= - \frac{\varpi^{(2)}_j (x)}{\varpi^{(2)}_{n-1}(x)}\;,\quad \text{for} \; r+1 \leq j < n-1\\ q^{(2)}_{ij} = 0 \quad \text{elsewhere} \end{cases}$$ 
The matrix $G_{A^{(2)}}(\lambda)$ is then of the form \eqref{particularform3} with \eqref{particularconditionb} and $\varrho=2$.\\
Consider the finite sequence of equivalent systems obtained by the transformation 
$$F^{(s-1)} = P^{(s)} ( I_n + Q^{(s)}(x))\; F^{(s)}$$ where $1 \leq s \leq \varrho$ and
$$Q^{(s)}(x) = [q^{(s)}_{ij}(x)]_{1 \leq i, j \leq n} \; , \quad s.t.\; \begin{cases} q^{(s)}_{n-s+1, j}  (x)= - \frac{\varpi^{(s)}_j (x)}{\varpi^{(s)}_{n-s+1}(x)}\;,\quad \text{for} \; r+1 \leq j < n-s+1\\ q^{(s)}_{ij} = 0 \quad \text{elsewhere} . \end{cases} $$  
Then this process yields an equivalent matrix $\tilde{A}(x, \xi_{\tilde{A}}) := A^{(\varrho)}(x, \xi_{A^{(\varrho)}})$ with \eqref{particularconditionb} for which either \eqref{particularconditiona} occurs or $\varrho=n-r$. But in the latter case one has, again by Remark~\ref{gq}, that $\det (A_0^{11}(x) )= 0$, and so \eqref{particularconditiona} occurs. 
\end{proof}
\begin{example}[Continue Example $\ref{firststep}$]
\label{exm7}
A simple calculation shows that $\det (G_{A}(\lambda)) \equiv 0$ hence $A$ is $\varepsilon$-reducible.  From \eqref{firststepg}, for $\lambda =0$, we have the singular matrix
$$G_{A}(\lambda=0) =   \begin{bmatrix}  0 & 0 & x+1  & 0 \\ x^2 & 0 & 0 & -2 x\\ -x & 0 &  0  & 2 \\ 0 & 2 & 0 & 0 \end{bmatrix}.\; \text{Let} \; T= \begin{bmatrix}  1 & 0 & 0  & 0 \\ 0 & 1 & 0 & 0\\ 0 & 0 & 0 & 1 \\ 0 & 0 & 1 &  0\end{bmatrix}$$ 
then the transformation $F = T G$ yields the equivalent system $\xi\xi_{\tilde{A}}^h \partial G = \tilde{A}(x,\xi_{\tilde{A}}) G$ where $\sigma_{\tilde{A}}= \sigma_A$ and
$$\tilde{A}(x,\xi_{\tilde{A}}) =   \begin{bmatrix} \xi_{\tilde{A}}  & -x^3 \xi_{\tilde{A}} & 0 & (1+x) \xi_{\tilde{A}}\\ x^2 & x \xi_{\tilde{A}} & -2x \xi_{\tilde{A}} & 0 \\ 0 & 2 & \xi_{\tilde{A}}^2 & 0 \\ -x & 0 & 2 \xi_{\tilde{A}} & 0 \end{bmatrix} . $$ 
$G_{\tilde{A}}(\lambda)$ has the form \eqref{particularform3} with $\varrho= 1$ and $r=2$. In fact, $$G_{\tilde{A}}(\lambda) = \begin{bmatrix} 0  & 0 & 0 & (1+x)  \\ x^2 & 0 & 0 & 0 \\ 0 & 2 & \lambda & 0 \\ -x & 0 & 2 &  \lambda \end{bmatrix} . $$
\end{example}
\begin{lemma}
\label{shearing}
Given the system
$$
[A_{\sigma_A}] \quad \quad \quad \xi_A^h \; x^p \; \partial F=  A(x, \xi_A) F = \sum_{k=0}^{\infty} A_k(x) \xi_A^k \; F .$$ 
Set $r = rank (A_0(x))$ and suppose that $m_{\varepsilon}(A) >1$ and $G_{A}(\lambda)$ has the form \eqref{particularform3} with conditions \eqref{particularconditionb} and \eqref{particularconditiona} satisfied. Consider the shearing  transformation $S(\xi_A) = \operatorname{diag} (\xi_A I_r , I_{n-r-\varrho}, \xi_A I_\varrho)$ if $\varrho \neq 0$ and $S(\xi_A)= \operatorname{diag} (\xi_A I_r , I_{n-r})$ otherwise. Then $F = S(\xi_A) G$ yields the equivalent system
$$
[\tilde{A}_{\sigma_{\tilde{A}}}] \quad \quad \quad \xi_{\tilde{A}}^h \; x^p \; \partial G=  \tilde{A}(x, \xi_{\tilde{A}}) G = \sum_{k=0}^{\infty} \tilde{A}_k(x) \xi_{\tilde{A}}^k \; G ,$$ 
for which $rank (\tilde{A}_0(x)) < r$ and $\sigma_{\tilde{A}} \geq \begin{cases} \sigma_A  \quad \text{if} \quad p>1 \\  \sigma_A + \frac{p-1}{h}  \quad \text{otherwise}. \end{cases} $
\end{lemma}
\begin{proof}
We partition $A(x, \xi_A)$ as follows ( we drop $(x, \xi_A)$ for clarity)
$$A(x, \xi_A)= \begin{bmatrix} A^{11}  & A^{12} & A^{13} \\ A^{21}  & A^{22} & A^{23} \\  A^{31}  & A^{32} & A^{33} \end{bmatrix}$$
where $A^{11} , A^{22} , A^{33} $ are of dimensions $r, n-r-\varrho,$ and  $\varrho$ respectively. 
It is easy to verify then  that \begin{eqnarray*} \tilde{A}(x, \xi_A) &=& S^{-1} A S - S^{-1} \partial S \\ &=& S^{-1} A S - \xi_A^h x^{p-1} \operatorname{diag} ( \sigma I_r , O_{n-r-\varrho}, \sigma I_\varrho)\\ &=& \begin{bmatrix} A^{11} - \xi_A^h x^{p-1} \sigma I_r &  \xi_A^{-1} A^{12} & A^{13} \\   \xi_A A^{21}  & A^{22} &   \xi_A A^{23} \\  A^{31}  &   \xi_A^{-1} A^{32} & A^{33} - \xi_A^h x^{p-1} \sigma_A I_\varrho \end{bmatrix}.\end{eqnarray*}
Hence, the new leading coefficient matrix is 
$$\tilde{A}_0(x) = \begin{bmatrix} A_0^{11} & U_1 & O \\ O& O &O \\ M_1 & M_2 & O \end{bmatrix}$$  and $rank (\tilde{A}_0(x)) = rank( A_0^{11} \quad  U_1) < r.$
Moreover, due to the entries of the form $\xi_A^h x^{p-1} \sigma_A$ in $\tilde{A}$, $\sigma_{A}$ should be adjusted as claimed if $p<1$. 
\end{proof}
\begin{algorithm}
\label{algoparamgeneral}
\caption{$\varepsilon$-Rank Reduction of System  $[A_{\sigma_A}]$ for $h>1$ }
\label{algorithmparam}
\textbf{Input:} $h,\; p, \; \sigma_A,\;A(x,\xi_A)$ of   $[A_{\sigma_{A}}]$ \\
\textbf{Output:} $R \in GL_n(\rm \mathcal{K})$ and an equivalent system $[\tilde{A}_{\sigma_{\tilde{A}}}]$ which is $\varepsilon$-irreducible.
\begin{algorithmic}
\State $R  \gets  I_{n}$; 
\State $h  \gets  \varepsilon$-rank of $A$;
\State $U(x)  \gets $ Lemma~\ref{gauss1} so that $U^{-1} A_0(x) U$ has form \eqref{gaussform}; 
\State $R  \gets  R U$;
\State $A \gets  U^{-1} A U - \xi_A^h x^p U^{-1} \partial U$; Update $\sigma_A$;
\State $d= \det (G_{\lambda}(A))$;
\While {$d=0$ and $h>0$} \do \\
\If{$h>1$}
\State $T(x), \varrho  \gets $ Proposition \ref{gauss3} ; 
\State $A \gets  T^{-1} A T - \xi_A^h x^p T^{-1} \partial T$; Update $\sigma_A$;
\State $S(\xi_A)  \gets $ Lemma \ref{shearing}; 
\State $A \gets  S^{-1} A S - \xi_A^h x^p S^{-1} \partial S$; Update $\sigma_A$;
\State $R  \gets  R T S$;
\EndIf
\State $U(x)  \gets $  Lemma~\ref{gauss1}; Update $\sigma_A$;
\State $R  \gets  R U$;
\State $A \gets  U^{-1} A U - \xi_A^h x^p U^{-1} \partial U$;
\State $d = \det (G_{\lambda}(A))$;
\State $h  \gets  \varepsilon$-rank of $A$;
\EndWhile\\
\Return{(R, A,h)}.
\end{algorithmic}
\end{algorithm}
\begin{example}[Continue Example~\ref{exm7}] Let $S (\xi_{\tilde{A}}) = \operatorname{diag} (\xi_{\tilde{A}} , \xi_{\tilde{A}}, 1, \xi_{\tilde{A}})$ then $G = S(\xi_{\tilde{A}}) U$ yields $\xi_{\tilde{\tilde{A}}}^h \partial U = \tilde{\tilde{A}}(x, \xi) U$ where $\sigma_{\tilde{\tilde{A}}}= \sigma_{\tilde{A}}$ and
$$\tilde{\tilde{A}}(x, \xi_{\tilde{\tilde{A}}}) = \begin{bmatrix} \xi_{\tilde{\tilde{A}}} & -x^3 \xi_{\tilde{\tilde{A}}}  & 0 & (1+x) \varepsilon \\ x^2 &x \xi_{\tilde{\tilde{A}}} & -2x & 0 \\ 0 & 2  \xi_{\tilde{\tilde{A}}} & \xi_{\tilde{\tilde{A}}}^2 & 0 \\ -x & 0  & 2 & 0\end{bmatrix} . $$ It is clear that the leading term $\tilde{\tilde{A}}_0(x)$ has rank $1<  2 = rank (A_0(x))$. 
\end{example}

\begin{proof} (Theorem \ref{mosernilpotent2}) By Lemma~\ref{gauss1}, we can assume that $A_0(x)$ is in the form \eqref{gaussform}. Then,  $G_{A}(\lambda)$ is constructed as in \eqref{glambdaform}. Since $\det (G_{A}(\lambda)) \equiv 0$, it suffices to take the change of basis $F = R G = T S G$, where $T$ and $S$ are as in Propositions \ref{gauss3} and Lemma \ref{shearing} respectively. 
\end{proof}
\begin{remark} The $\varepsilon$-reducibility of $[A_{\sigma_A}]$ implies that the rank of the leading coefficient matrix can be reduced (and consequently the $\varepsilon$-Moser rank) without necessarily reducing the $\varepsilon$-rank $h$ of the system. If the $\varepsilon$-reduction criterion holds for a sufficient number of equivalent systems then a repetitive application of such a transformation results in an equivalent system whose leading coefficient matrix has a zero rank, hence $h$ can be reduced at least by one (e.g. Example~\ref{algoexm}). 
\end{remark}

\section{Proofs for Section~\ref{exporder}}
\label{appd}
We give the proof of Theorem~\ref{katzmain} after establishing a series of useful lemmas. The following proofs are an adaptation to the parametrized setting of the proofs in~\cite[Lemma 3, Lemma 4, Proposition 1, Theorem 1]{key24} respectively. For clarity within these intermediate proofs, we will express systems  in the equivalent notation $[A_{\sigma_A}] \quad \partial F = A(x, \varepsilon) F$ where $A(x, \varepsilon) \in \mathcal{M}_n(\rm \mathcal{K})$, rather than $[A_{\sigma_A}] \quad \xi_A^h x^p \partial F = \tilde{A}(x, \xi_A) F$. Let
\begin{equation}
\label{det1}
 \det (\lambda I - A(x, \varepsilon)) = \lambda^n + \alpha_{n-1}(x, \varepsilon) \lambda_{n-1} + \dots + \alpha_0 (x, \varepsilon)  .\end{equation}
such that $\alpha_n =1$ and $\alpha_i(x, \varepsilon) = \sum_{j= val_{\varepsilon}(\alpha_i)}^{\infty} \alpha_{i,j}(x) \varepsilon^j \in \rm \mathcal{K}$ for $i \in \{0, \dots, n\}$. We define the $\varepsilon$-polygon $\mathcal{N}_{\varepsilon} (A)$ of $[A_{\sigma_A}]$ as in Section~\ref{macutanpolygon}, by taking $P_{\varepsilon}(A)$  to be the union of $P(i, \;val_{\varepsilon}\; (\alpha_i (x, \varepsilon))$ for $i \in \{ 1, \dots, n\}$. We thus prove the following theorem, of which Theorem~\ref{katzmain} is a straightforward corollary.
\begin{theorem}
\label{katzmain11}
Consider the $\varepsilon$-irreducible system $[A_{\sigma_A}] \quad \partial F = A(x, \varepsilon) F$ where $A(x, \varepsilon) \in \mathcal{M}_n(\rm \mathcal{K})$, $h>0$, and~\eqref{det1}. If  $h > n - rank (A_0(x))$ then the $\varepsilon$-formal exponential order is given by 
$$\omega_{\varepsilon} (A)= \maxi_{0 \leq i < n}\; (\frac{- val_{\varepsilon}(\alpha_i)}{n - i}).$$
Additionally, the corresponding $\varepsilon$-polynomial is given by the algebraic equation
$$E_{\varepsilon}(X) = \sum_{k=0}^{\ell} {\alpha}_{i_k, val_{\varepsilon} (\alpha_{i_k})}\; X^{(i_k - i_0)}$$ 
where $0 \leq i_0 < i_1 < \dots < i_\ell = n$ denote the integers $i$ for which $\omega_{\varepsilon}(n-i)= - val_{\varepsilon}(\alpha_i)$ (i.e.\ lie on the edge of slope $\omega_{\varepsilon}$ of the $\varepsilon$-polygon $\mathcal{N}_{\varepsilon} (A)$ of $[A_{\sigma_A}]$); and ${\alpha}_{i, val_{\varepsilon} (\alpha_{i})}(x) = \varepsilon^{-val_{\varepsilon}\; (\alpha_{i})}\; \alpha_{i} (x, \varepsilon) |_{\varepsilon = 0}$.
\end{theorem}

\begin{lemma}
\label{katzlemma1}
Let $A(x, \varepsilon)$, $W(x, \varepsilon)$ be matrices in $\mathcal{M}_n(\rm \mathcal{K})$ and $\mathcal{M}_n(\rm{\mathcal{R}})$ respectively.  Put $h= max (0, -val_{\varepsilon}\;(A(x, \varepsilon)))$ and let
$$ \det (\lambda I - A(x, \varepsilon)) -  \det (\lambda I - A(x, \varepsilon) + W(x, \varepsilon)) = \alpha_{n-1}(x, \varepsilon) \lambda^{n-1} + \alpha_{n-2}(x, \varepsilon) \lambda^{n-2} + \dots + \alpha_0(x, \varepsilon),$$
where $\alpha_{n-1}(x, \varepsilon), \alpha_{n-2}(x, \varepsilon), \dots, \alpha_0(x, \varepsilon)$ lie in $\rm \mathcal{K}$. Then 
$$val_{\varepsilon} (\alpha_{n-i}) \geq (1-i) h \quad 1 \leq i \leq n .$$ 
\end{lemma}
\begin{proof}
For any $i \in \{1, \dots, n\}$, it follows from Cramer's rule that $$\alpha_{n-i} = \sum_{l=0}^{\gamma}( w_l \prod_{s=1}^{i-1} a_{l,s}) $$ 
where for $0\leq l \leq \gamma,  1\leq s\leq i-1, $ $w_l$ are entries in $W(x, \varepsilon)$ and $a_{l,s}$ are entries in $W(x, \varepsilon)$ or $A(x, \varepsilon)$. Consequently, 
$$val_{\varepsilon} (\alpha_{n-i}) \geq (i-1) val_{\varepsilon} (A(x, \varepsilon)) \geq h . $$
\end{proof}
\begin{lemma}
\label{katzlemma2}
Let $A(x, \varepsilon)$, $B(x, \varepsilon)$ be two matrices in $\mathcal{M}_n(\rm \mathcal{K})$ such that $val_{\varepsilon} (A(x, \varepsilon)) \leq val_{\varepsilon} (B(x, \varepsilon))$. Consider the two systems $\; [A_{\sigma_A}] \; \partial F = A(x, \varepsilon) F$ and $\; [B_{\sigma_B}] \; \partial G = B(x, \varepsilon) G$. Suppose that there exists $T \in GL_n(\rm \mathcal{K})$ such that $ [B_{\sigma_B}] = T[A_{\sigma_A}]$. Put $h= max (0, -val_{\varepsilon}\;(A(x, \varepsilon)))$ and let
$$\begin{cases}  \det (\lambda I - A(x, \varepsilon)) =  \lambda^{n} + \alpha_{n-1} \lambda^{n-1} + \alpha_{n-2} \lambda^{n-2} + \dots + \alpha_0 \\
 \det (\lambda I - B(x, \varepsilon)) =  \lambda^{n} + \beta_{n-1} \lambda^{n-1} + \beta_{n-2} \lambda^{n-2} + \dots + \beta_0 .\end{cases}
$$
Then, we have:
$$val_{\varepsilon} (\alpha_{n-i}(x, \varepsilon) - \beta_{n-i}(x, \varepsilon)) \geq (1-i) h , \quad 1 \leq i \leq n .$$
\end{lemma}
\begin{proof}
By~\cite[Lemma 1]{key90} we can write $T(x, \varepsilon) = P(x, \varepsilon)\; \varepsilon^{\gamma}\; Q(x, \varepsilon)$ where $P, Q \in GL_n (\rm{\mathcal{R}})$ and $\varepsilon^{\gamma} = \operatorname{diag} (\varepsilon^{\gamma_1}, \dots, \varepsilon^{\gamma_n} )$ for some integers $(\gamma_1 \leq \gamma_2 \leq \dots \leq \gamma_n)$.
Consider $[\tilde{A}_{\sigma_{\tilde{A}}}] := P [A_{\sigma_A}]$ and $[\tilde{B}_{\sigma_{\tilde{B}}}] := Q^{-1} [B_{\sigma_B}]$. Then we have 
$$\tilde{B} = \varepsilon^{- \gamma} \tilde{A} \varepsilon^{\gamma}, \quad val_{\varepsilon} (\tilde{A}) = val_{\varepsilon} (A), \quad \text{and} \; val_{\varepsilon} (\tilde{B}) = val_{\varepsilon} (B). $$
It follows that,
$$\begin{cases}
  \det (\lambda I - \tilde{A}) =  \det (\lambda I - P^{-1} A P + P^{-1}\partial P) =  \det (\lambda I - A + (\partial P) P^{-1}) \\
  \det (\lambda I - \tilde{B}) =  \det (\lambda I - Q B Q^{-1}  + Q  \partial Q ^{-1}) = \det (\lambda I - B + (\partial Q ^{-1}) Q) \\
  \det (\lambda I - \tilde{B}) =  \det (\lambda I - \varepsilon^{- \gamma} \tilde{A} \varepsilon^{\gamma}) =  \det (\lambda I - \tilde{A} ) .
\end{cases}$$
Since $P, Q, P^{-1}, Q^{-1}, $ are units of $\mathcal{M}_n(\rm{\mathcal{R}})$, it follows that $(\partial P) P^{-1}$ and $(\partial Q ^{-1}) Q$ inherit this property as well. 
The rest of the proof follows as a consequence of Lemma~\ref{katzlemma1}. 
\end{proof}

\begin{proposition}
\label{katzprop1}
Consider the system $[A_{\sigma_A}] \quad \partial F = A(x, \varepsilon) F$  where $A \in \mathcal{M}_n(\rm \mathcal{K})$ and let:
$$ \det (\lambda I - A(x, \varepsilon)) =  \lambda^{n} + \alpha_{n-1}(x, \varepsilon) \lambda^{n-1} + \alpha_{n-2}(x, \varepsilon) \lambda^{n-2} + \dots + \alpha_0(x, \varepsilon).$$
Let $[C_{\sigma_C}] \quad \partial G = C(x,\varepsilon) G$ where \begin{equation} \label{companion} C(x, \varepsilon)\; = \; Companion \; (c_i (x, \varepsilon))_{0 \leq i \leq n-1} = \begin{bmatrix} 0 & 1 & 0 & \dots & 0 \\ 0 & 0& 1 & \dots  & 0 \\ \vdots&\vdots&\vdots&\vdots&\vdots \\  0 & 0& 0 & \dots  & 1 \\ c_0 & c_1 & c_2 & \dots & c_{n-1} \end{bmatrix}, \end{equation} be a companion system which is equivalent to $[A_{\sigma_A}]$ over $GL_n(\rm \mathcal{K})$.
 Then we have,
$$val_{\varepsilon}\; (\alpha_{i} - c_{i}) \geq  (1- (n-i)) h \quad 0 \leq i \leq n-1 .$$
\end{proposition}

\begin{proof}
Let $[c_{\sigma_c}]$ with $c_n = 1$ be the scalar equation representing $[C_{\sigma_C}]$ ($\sigma_c = \sigma_C$). Consider $\kappa$ of $[c_{\sigma_c}]$ as defined in Section~\ref{mosereqn}. We define $ \beta_i = (n-i) \kappa,$ for $i \in \{ 0 , \dots, n-1\}$ and $[D_{\sigma_D}]= \varepsilon^{\beta}[C_{\sigma_C}]$. Then $D(x, \varepsilon) = \varepsilon^{-\beta} C(x, \varepsilon) \varepsilon^{\beta}$ where $\varepsilon^\beta = \operatorname{diag} (\varepsilon^{\beta_0}, \dots, \varepsilon^{\beta_{n-1}})$. It follows that $D(x, \varepsilon) = \varepsilon^{-\kappa}\; Companion\; (\varepsilon^{\beta_i} c_i(x, \varepsilon))_{0 \leq i \leq n-1}$. We have, $val_{\varepsilon} (D(x, \varepsilon)) \geq -\kappa$ since $val_{\varepsilon} (\varepsilon^{\beta_i} c_i) = (n-i) \kappa +  val (c_i) \geq 0$.  
By the equivalence between $[A_{\sigma_A}]$ and $[C_{\sigma_C}]$ (resp. $[D_{\sigma_D}]$) we have $h \geq \kappa$ (resp. $val_{\varepsilon} (D(x, \varepsilon)) \geq -\kappa \geq -h$). Let  
$$ \det (\lambda I - D(x, \varepsilon)) = \lambda^n + d_{n-1}(x, \varepsilon) \lambda^{n-1} + \dots + d_0(x, \varepsilon) .$$
Hence, by Lemma~\ref{katzlemma2}, we have
$$val_{\varepsilon}\; (\alpha_i - d_i) \geq (1-(n-i)) h , \; 0 \leq i \leq n . $$
Moreover, by Lemma~\ref{katzlemma1} 
$$val_{\varepsilon} (d_i - c_i) \geq (1-(n-i))\; max\; (0, -\; val_{\varepsilon}\; (D(x, \varepsilon))) \geq (1-(n-i)) \kappa \geq (1-(n-i)) h . $$
It follows that $$val_{\varepsilon} (\alpha_i - c_i) \geq \; min\; (val_{\varepsilon}\; (\alpha_i - d_i), \; val_{\varepsilon}\; ( d_i - c_i)) \geq (1- (n-i)) h . $$
\end{proof}
\begin{proof}(Theorem~\ref{katzmain11})
Let $[C_{\sigma_C}] \quad \partial G = C(x,\varepsilon) G$ be as in Proposition~\ref{katzprop1}. Due to their equivalence, $[A_{\sigma_A}]$ and $[C_{\sigma_C}]$ have the same $\varepsilon$-exponential order and $\varepsilon$-polynomial. Hence, we have
$$\omega_{\varepsilon}:=\omega_{\varepsilon} (A) = \omega_{\varepsilon} (C) = max_{0 \leq i < n}\; (\frac{- val_{\varepsilon}(c_i)}{n - i}) \quad \text{and} \quad E_{\varepsilon}(X) = \sum_{k=0}^{\ell} c_{i_k, val_{\varepsilon} (c_{i_k})} X^{(i_k - i_0)} $$ 
where $0 \leq i_0 < i_1 < \dots < i_\ell = n$ denote the integer $i$ for which $\omega_{\varepsilon}(n-i)= - val_{\varepsilon}(c_i)$. 
By Corollary~\ref{qkappa} and Proposition~\ref{katzprop1} one has 
\begin{eqnarray*}
val_{\varepsilon}\; (\alpha_i - c_i) & \geq & (i-n+1) h = (i-n)  \omega_{\varepsilon} + (i-n)(h- \omega_{\varepsilon}) + h \\
&\geq&  \omega_{\varepsilon} (i-n) + (-n) (1- \frac{r}{n}) + h \geq  \omega_{\varepsilon} (i-n) + r+ h-n .
\end{eqnarray*}
It follows from $h + r - n>0$ that $ val_{\varepsilon}\; (\alpha_i - c_i) > (i-n) \omega_{\varepsilon},\; 0 \leq i \leq n$. 
\begin{itemize}
\item If $i \in \{  i_0, \dots i_\ell\}$ then $\omega_{\varepsilon} (i-n) = val_{\varepsilon}\; (c_i)$ and $val_{\varepsilon}\;(\alpha_i - c_i) > \; val_{\varepsilon}\; (c_i)$. Hence, $val_{\varepsilon}\; (\alpha_i) = \; val_{\varepsilon}\; (c_i)$ and $\alpha_{i, \; val_{\varepsilon}\; (\alpha_i)} = c_{i, \; val\; (c_i)}$. 
\item Else $val_{\varepsilon}\;(\alpha_i - c_i) > \omega_{\varepsilon} (i-n)$ and $val_{\varepsilon}\; (c_i) > \omega_{\varepsilon} (i-n)$. Thus, 
$$  val_{\varepsilon}\; (\alpha_i) \geq \; min\; (val_{\varepsilon}\; (\alpha_i - c_i), \; val_{\varepsilon}\; (c_i)) > \omega_{\varepsilon}(i-n) .$$
\end{itemize}
This completes the proof. 
\end{proof}
We also illustrate Lemma~\ref{katzlemma} with the following example:
\begin{example}[Lemma~\ref{katzlemma}]
Consider the system $[A_{\sigma_A}] \quad  \xi_A^2 \partial F = A(x,\xi_A) F$ with $\sigma_A=0$ and
$$A(x,\xi_A)= \left[ \begin {array}{ccccccc} 0&0&{\xi_A}^{2}&0&\xi_A&x&0
\\ \noalign{\medskip}{\xi_A}^{2}&{\xi_A}^{3}&x\xi_A&0&0&0&0
\\ \noalign{\medskip}0&0&0&x&0&0&4\\ \noalign{\medskip}3\,\xi_A\,{x
}^{2}&0&0&0&{\xi_A}^{2}&x\xi_A&0\\ \noalign{\medskip}0&x\xi_A
&0&0&0&{\xi_A}^{2}&0\\ \noalign{\medskip}0&0&x\xi_A&{\xi_A}^{
2}&{\xi_A}^{2}&0&0\\ \noalign{\medskip}0&{\xi_A}^{3}&0&0&0&0&x
\xi_A\end {array} \right].$$
$A(x, \xi_A)$ is $\varepsilon$-irreducible, and we have $n=7$, $h=2$, and $r:= rank(A_0(x))=2$. Thus, the condition $h > n-r$ of Theorem~\ref{katzmain} is not verified. Let us first consider a random ramification in $\xi_A$ regardless of Lemma~\ref{katzlemma}. For instance, let us consider $\xi_A = \tilde{\xi}_A^3$ and apply $\tilde{\varepsilon}$-rank reduction. This yields system  $[C_{\sigma_C}] \quad \tilde{\xi}_C^5 \partial G = C(x,\tilde{\xi}_C) G$ with $\sigma_C=0$ and
$$C(x,\tilde{\xi}_C)=  \left[ \begin {array}{ccccccc} 0&0&0&x\tilde{\xi}_C&{\tilde{\xi}}_C^{4}&{\tilde{\xi}}_C^{4}&0\\ \noalign{\medskip}1/4\,{\tilde{\xi}}_C^{2}{x}^{2}&x {\tilde{\xi}}_C^{2}&{\tilde{\xi}}_C^{7}&0&1/4\,\ {\tilde{\xi}}_C \, \left( {\tilde{\xi}}_C^{3}-{
x}^{2} \right) &1/4\,{\tilde{\xi}}_C^{4}x&3/4\,{\tilde{\xi}}_C\,{x}^{3}
\\ \noalign{\medskip}0&0&{\tilde{\xi}}_C^{8}&x{\tilde{\xi}}_C^{2}&0&0&{\tilde{\xi}}_C^{5}\\ \noalign{\medskip}0&4&0&0&0&0&0\\ \noalign{\medskip}x{\tilde{\xi}}_C
^{3}&0&0&0&0&{\tilde{\xi}}_C^{5}&3\,{\tilde{\xi}}_C^{2}{x}^{2}
\\ \noalign{\medskip}{\tilde{\xi}}_C^{6}&0&x{\tilde{\xi}}_C^{2}&0&0&0&0
\\ \noalign{\medskip}x&0&0&{\tilde{\xi}}_C^{5}&0&{\tilde{\xi}}_C^{2}&0
\end {array} \right] , $$
for which, $\tilde{h}=5$ and $\tilde{r}=2$. Thus, the condition $\tilde{h}>n-\tilde{r}$ is still not necessarily verified after a random ramification. To guarantee that we will arrive at a system verifying this condition, we make use of Lemma~\ref{katzlemma} and choose an integer $d$ such that $d \geq 7/(1 + 2/7)= 49/9$. Let $d=7 > 49/9$, perform $\xi_A = \tilde{\xi}_A^7$ in $[A_{\sigma_A}]$, and then $\tilde{\varepsilon}$-rank reduction. This yields the system 
$[B_{\sigma_B}] \quad \tilde{\xi}_B^{11} \partial G = B(x,\tilde{\xi}_B) G$
with $\sigma_B=0$ and
$$ B(x,\tilde{\xi}_B) = \left[ \begin {array}{ccccccc} {\tilde{\xi}}_B^{4}x&1/4\,{\tilde{\xi}}_B^{4}{x}
^{2}&{\tilde{\xi}}_B^{15}&0&1/4\,\tilde{\xi}_B\, \left( {\tilde{\xi}}_B^{7}-{x}^{2}
 \right) &1/4\,{\tilde{\xi}}_B^{8}x&3/4\,\tilde{\xi}_B\,{x}^{3}
\\ \noalign{\medskip}0&0&0&\tilde{\xi}_B\,x&{\tilde{\xi}}_B^{8}&{\tilde{\xi}}_B^{8}&0
\\ \noalign{\medskip}0&0&{\tilde{\xi}}_B^{18}&{\tilde{\xi}}_B^{4}x&0&0&{\tilde{\xi}
}_B^{11}\\ \noalign{\medskip}4&0&0&0&0&0&0\\ \noalign{\medskip}0&{
\tilde{\xi}}_B^{7}x&0&0&0&{\tilde{\xi}}_B^{11}&3\,{\tilde{\xi}}_B^{4}{x}^{2}
\\ \noalign{\medskip}0&{\tilde{\xi}}_B^{14}&{\tilde{\xi}}_B^{4}x&0&0&0&0
\\ \noalign{\medskip}0&x&0&{\tilde{\xi}}_B^{11}&0&{\tilde{\xi}}_B^{4}&0
\end {array} \right] .
$$
For system $[B_{\sigma_B}]$, we have: $\tilde{h} = 11$, $\tilde{r}=2$, and so the condition $\tilde{h} = 11 > n-\tilde{r} = 7-2$ holds, and so, Theorem~\ref{katzmain} can now be applied. From $\det(\lambda I - \frac{B(x,\tilde{\xi}_B)}{{\tilde{\xi}}_B^{11}})$, we can compute $\omega_{varepsilon} (M) = \frac{42}{4}$. And so, we introduce the ramification $\tilde{\xi}_B = \tilde{\tilde{\xi}}_B^4$ and apply $\tilde{\tilde{\varepsilon}}$-rank reduction which gives a $\tilde{\tilde{\varepsilon}}$-system of order $42$ and whose leading coefficient has $4$ nonzero eigenvalues.
 In fact, we compute
$[\tilde{M}_{\tilde{\tilde{\xi}}_M}] \quad \tilde{\tilde{\xi}}_M^{42} \partial W = \tilde{M}(x,\tilde{\tilde{\xi}}_M) W$
with $\sigma_M=0$ and
$$ \tilde{M}(x,\tilde{\tilde{\xi}}_M) =  \left[ \begin {array}{ccccccc} {\tilde{\tilde{\xi}}_M}^{14}x&1/4\,{\tilde{\tilde{\xi}}_M}^{14}{
x}^{2}&0&0&1/4\,{\tilde{\tilde{\xi}}_M}^{28}-1/4\,{x}^{2}&1/4\,{\tilde{\tilde{\xi}}_M}^{28}x&3/4
\,{x}^{3}\\ \noalign{\medskip}0&0&0&x&{\tilde{\tilde{\xi}}_M}^{28}&{\tilde{\tilde{\xi}}_M}^{28}&0
\\ \noalign{\medskip}0&0&0&{\tilde{\tilde{\xi}}_M}^{14}x&0&0&{\tilde{\tilde{\xi}}_M}^{42}
\\ \noalign{\medskip}4&0&0&0&0&0&0\\ \noalign{\medskip}0&{\tilde{\tilde{\xi}}_M}^{
28}x&0&0&0&{\tilde{\tilde{\xi}}_M}^{42}&3\,{\tilde{\tilde{\xi}}_M}^{14}{x}^{2}
\\ \noalign{\medskip}0&0&{\tilde{\tilde{\xi}}_M}^{14}x&0&0&0&0\\ \noalign{\medskip}0
&x&0&{\tilde{\tilde{\xi}}_M}^{42}&0&{\tilde{\tilde{\xi}}_M}^{14}&0\end {array} \right] .
$$
We remark that one can observe that this procedure was not necessary for such system since $gcd(14, 42,28) =14$. Thus a simplification would leave us with a system of order $3$. Thus, in the implementation we first try to ignore the condition $h > n-r$ and apply Theorem~\ref{katzmain} directly to the system $[A_{\sigma_A}]$. For this particular example we obtain $\omega_{\varepsilon}(A) = 3/2$. Then by performing $\xi_A = \tilde{\xi}_A^2$ and the $\tilde{\varepsilon}$-rank reduction, we arrive at a $\tilde{\varepsilon}$-irreducible system of order $3$ whose leading coefficient has $4$ nonzero eigenvalues. In fact, we compute $[\tilde{A}_{\sigma_{\tilde{A}}}] \quad  {\tilde{\xi}}_{\tilde{A}}^{3} \partial G = \tilde{A}(x,{\tilde{\xi}}) G$
where $\sigma_{\tilde{A}} = 0$ and
$$ \tilde{A}(x,{\tilde{\xi}}_{\tilde{A}}) =  \left[ \begin {array}{ccccccc} 0&0&0&x&{\tilde{\xi}}_{\tilde{A}}^{2}&{\tilde{\xi}}_{\tilde{A}}^{2}&0
\\ \noalign{\medskip}1/4\,{\tilde{\xi}}_{\tilde{A}} \,{x}^{2}&{\tilde{\xi}}_{\tilde{A}}\,x&{\tilde{\xi}}_{\tilde{A}}^{4}
&0&-1/4\,{x}^{2}+1/4\,{\tilde{\xi}}_{\tilde{A}}^{2}&1/4\,x{\tilde{\xi}}_{\tilde{A}}^{2}&3/4\,{x}^{3}
\\ \noalign{\medskip}0&0&{\tilde{\xi}}_{\tilde{A}}^{5}&{\tilde{\xi}}_{\tilde{A}}\,x&0&0&{\tilde{\xi}}_{\tilde{A}}^{3}
\\ \noalign{\medskip}0&4&0&0&0&0&0\\ \noalign{\medskip}x{\tilde{\xi}}_{\tilde{A}}^{2}
&0&0&0&0&{\tilde{\xi}}_{\tilde{A}}^{3}&3\,{\tilde{\xi}}_{\tilde{A}}\,{x}^{2}\\ \noalign{\medskip} {\tilde{\xi}}_{\tilde{A}}^{4}&0&{\tilde{\xi}}_{\tilde{A}}\,x&0&0&0&0\\ \noalign{\medskip}x&0&0&{\tilde{\xi}}_{\tilde{A}}^{3}&0&{\tilde{\xi}}_{\tilde{A}}&0\end {array} \right]  .
$$ If Theorem~\ref{katzmain} does not lead to the desired result, we then resort Lemma~\ref{katzlemma}. In the case of unperturbed system, it is also a matter of discussion whether this condition is necessary. M. Miyake recently \href{http://bcc.impan.pl/15AAGA/uploads/news/id67/poland-2.pdf}{claimed} to have such an example (where this condition is necessary).
\end{example}

\section{Examples}
\label{appe}
We treat with our algorithm examples from literature. For the  computation of the full the $\varepsilon$-exponential parts, we refer to our \textsc{Maple} package \textsc{ParamInt}~\cite{key500}.
\begin{example} [Continue Example~\ref{exmwasow1}]
\label{exmwasow2}
We resume the computation of the outer solution of system $[B_{\sigma_{B}}]$ given by:
$$[B_{\sigma_{B}}] \quad \quad \quad \quad \xi_{B}^2 x^5  \; \partial W^1 = B(x, \xi_B) W^1 = \{ \begin{bmatrix}  0 & 1\\ 0 & 0  \end{bmatrix}  +  \begin{bmatrix}  -1 & 0 \\ -1 & 0 \end{bmatrix} \xi_{B}  + \begin{bmatrix}  1 & -1 \\ 1 & -1 + x^4  \end{bmatrix} \xi_{B}^{2} + O(\xi_{B}^3) \; \}  \; W^1 ,$$
with $\sigma_B =-3$.  The leading matrix coefficient $ \begin{bmatrix}  0 & 1\\ 0 & 0  \end{bmatrix}$ is nilpotent and $\varepsilon$-irreducible. We compute the characteristic polynomial of $\frac{B(x,\xi_B)}{\xi_{B}^2 x^5}$:
$$ \lambda^2 - \frac{(\xi_B x^4-1)}{\xi_B x^5} \lambda + \frac{(\xi_B -1)(\xi^2_B x^4-1)}{\xi_B^3 x^{10}}. $$
Then the $\varepsilon$-formal exponential order is $\omega_{\varepsilon} = \frac{3}{2}$ and the $\varepsilon$-polynomial is $X^2 + \frac{1}{x}=0$.
Thus, by setting $\varepsilon = \tilde{\varepsilon}^2$, $\xi_B = \tilde{\xi}_B^2 x^3$ ($\tilde{\xi}_B = x^{-3} \tilde{\varepsilon}$) and applying $\tilde{\varepsilon}$-rank reduction of Section~\ref{moser} via $\operatorname{diag} (1, \tilde{\xi}_B)$, 
we get an $\tilde{\varepsilon}$-irreducible system whose leading matrix coefficient  is $\begin{bmatrix} 0 & 1 \\ - x^{-3} & 0  \end{bmatrix}$.
And so, by the turning point algorithm of Section~\ref{turnpt}, we apply the transformation $$\begin{bmatrix} 1 & 0 \\ 0& x^{-3/2} \end{bmatrix}$$  which yields (note that we can also express the former using a change of independent variable $x=t^2$) $$[\tilde{B}_{\sigma_{\tilde{B}}}] \quad \quad \quad \quad \tilde{\xi}_{\tilde{B}}^{3} x^{19/2} \partial U = \tilde{B}(x,\tilde{\xi}_{\tilde{B}}) U \quad \text{where} \quad \sigma_{\tilde{B}} = -3 \; \text{and} $$ 
$$ \tilde{B}(x,\tilde{\xi}_{\tilde{B}})= \begin{bmatrix} 0 & 1 \\ -1 &0 \end{bmatrix} + \begin{bmatrix}-x^{3/2} & 0 \\ 0 &0 \end{bmatrix} \tilde{\xi}_{\tilde{B}} + \begin{bmatrix}0 & 0 \\ x^{3} &0 \end{bmatrix}  \tilde{\xi}_{\tilde{B}}^2 + \begin{bmatrix}x^{9/2} & 0 \\ 0 & -x^{9/2} +x^{17/2} \end{bmatrix}  \tilde{\xi}_{\tilde{B}}^3 + O( \tilde{\xi}_{\tilde{B}}^4) .$$
The leading matrix has two distinct eigenvalues. Applying Splitting Lemma we get $$[\tilde{\tilde{B}}_{\sigma_{\tilde{\tilde{B}}}}] \quad \quad \quad \quad \tilde{\xi}_{\tilde{\tilde{B}}}^{3} x^{19/2} \partial R= \tilde{\tilde{B}}(x,\tilde{\xi}_{\tilde{\tilde{B}}}) R \quad \text{where} \quad \sigma_{\tilde{\tilde{B}}} = -3 \; \text{and}$$ 
\begin{eqnarray*} \tilde{\tilde{B}}(x,\tilde{\xi}_{\tilde{\tilde{B}}})= \begin{bmatrix} -i & 0 \\ 0 & i \end{bmatrix} + \begin{bmatrix} -\frac{1}{2}x^{3/2}  & 0 \\ 0 & -\frac{1}{2} x^{3/2}  \end{bmatrix} \tilde{\xi}_{\tilde{\tilde{B}}} + O( \tilde{\xi}_{\tilde{\tilde{B}}}^2) .\end{eqnarray*}
Thus this system can be decoupled to any desired precision. The solution follows by straightforward integration. Remark that, as expected, the system has the same $\varepsilon$-formal exponential order and $\varepsilon$-polynomials, as its equivalent scalar equation given in Example~\ref{exmwasow3}. 
The leading term of the $\varepsilon$-exponential part $Q$ is then given by: $\exp ( \begin{bmatrix} \frac{-2 i \sqrt{x}}{\varepsilon^{3/2}} & 0 \\ 0& \frac{2 i \sqrt{x}}{\varepsilon^{3/2}} \end{bmatrix})$. Moreover, since $\sigma_{\tilde{\tilde{B}}} = -3 $, we set $\rho_1= 1/3$. Further computations yield the full exponential part of outer solutions of the input system $[A_{\sigma_A}]$ of Example~\ref{exmwasow1}:
\begin{eqnarray*} & &\frac{x^2}{2\varepsilon^{2}} - \frac{1}{x\varepsilon} \quad \text{and} \quad \\ & & -\frac{1}{20 t^5 \varepsilon^{1/2}} -\frac{1}{2 t^2 \varepsilon} + \frac{2t}{\varepsilon^{3/2}} \quad \text{where} \quad x=-t^2. \end{eqnarray*} 
Moreover, since $\sigma_{\tilde{\tilde{B}}} = -3 $, we set $\rho_1= 1/3$.
For the initial system of Example~\ref{exmwasow1} given by: $$\varepsilon^2  \partial F=  \begin{bmatrix}  0 & 1& 0 \\ 0 & 0 & 1 \\  \varepsilon & 0 & x \end{bmatrix} F ,$$
we have obtained so far outer formal solutions. To compute inner (and probably intermediate) solutions, we set $\tau = x \varepsilon^{- \rho_A} =  x \varepsilon^{-1/3}$ which yields:
$$\varepsilon^{\frac{5}{3}}  \partial_{\tau} F =  \begin{bmatrix}  0 & 1& 0 \\ 0 & 0 & 1 \\  \varepsilon & 0 & \tau \varepsilon^{\frac{1}{3}} \end{bmatrix} F ,$$
which can be expressed equivalently as:
$$[E_{\sigma_E}] \quad \quad \quad \quad \tilde{\xi}_E^5 \partial_{\tau} F =  E(\tau, \tilde{\xi}_E) F = \begin{bmatrix}  0 & 1& 0 \\ 0 & 0 & 1 \\  \tilde{\xi}_E^3 & 0 & \tau \tilde{\xi}_E \end{bmatrix} F , \; \text{with} \; \sigma_E =0 \; \text{and} \; \xi_E = \varepsilon = \tilde{\varepsilon}^3= \tilde{\xi}_E^3.$$
By Algorithm \ref{algorithmkatzmain3}, one can compute the transformation $F = T G$ where 
$$T = \begin{bmatrix} 0 & 0 & 1\\ 0 & 1& 0 \\ 1 & 0&0\end{bmatrix}\; \begin{bmatrix} \tilde{\xi}_E^2 & 0 & 0 \\ 0 & \tilde{\xi}_E & 0 \\ 0 & 0&1\end{bmatrix}.$$
The resulting system is 
$$[\tilde{E}_{\sigma_{\tilde{E}}}] \quad \quad \quad \quad \tilde{\xi}_{\tilde{E}}^4 \partial_{\tau} G = \tilde{E}(x,\tilde{\xi}_{\tilde{E}}) G = \begin{bmatrix} \tau & 0&1 \\ 1 &0&0 \\ 0 & 1&0 \end{bmatrix} G ,\quad \text{with}\; \sigma_{\tilde{E}} = 0.$$
One can verify that $\tilde{E}_{0,0}$ has three distinct roots and so the system can be decoupled into three scalar equations. The leading term of the $\varepsilon$-exponential part in the inner subdomain is then given by: $$\exp ( \begin{bmatrix} \frac{\int^{\tau} dz + \dots}{\varepsilon^{4/3}} & 0 & 0 \\ 0& \frac{\int^{\tau} -1 + i/\sqrt{3}+ \dots dz}{2 \varepsilon^{4/3}} & 0 \\ 0 & 0 & \frac{\int^{\tau} -1 - i/\sqrt{3}+ \dots dz}{2 \varepsilon^{4/3}} \end{bmatrix}) .$$
And further computations show that the diagonal entries of the $\varepsilon$-exponential part of a fundamental matrix of formal solutions in the inner subdomain are given by:
\begin{eqnarray*} & &\frac{\tau + (1/6) \tau^2 + (1/27) \tau^3 + O(\tau^4)}{\varepsilon^{4/3}}\quad \text{and} \quad \\ & & \frac{\tau \operatorname{RootOf} (z^2 +z +1)+ (1/6) \tau^2 - (1/27) (\operatorname{RootOf} (z^2 +z +1) +1) \tau^3 + O(\tau^4)}{\varepsilon^{4/3}} .\end{eqnarray*} 
Since $\sigma_{\tilde{E}} = 0$, the reduction stops here.
\end{example}

\begin{example}[Iwano-Sibuya polygon]
Consider the following scalar equation whose $\sigma_a=0$:
 \begin{equation} 
\label{exm40}
[a_{\sigma_a}] \quad \quad \quad \quad \partial^4\;f  \;-\; \varepsilon^{-4} (2x+\varepsilon^3) \partial^2\;f\; \;+ \;\varepsilon^{-8} x^3  f \;=\; 0 ,
\end{equation}
One can verify that $\omega_{\varepsilon} =2$. The Iwano-Sibuya polygon  is given by the following set of points: $P_{\omega_{\varepsilon}(a)} = (2,-1)$, $P_{00} = (0,\frac{3}{4})$,  $P_{20} = (0,\frac{1}{2})$, and $P_{23} = (\frac{3}{2},0)$ leading to only one slope given by $\sigma = \frac{-3}{4}$. Then the behavior can be investigated with the help of $[\rho] \quad  0 < \rho_1 = -1/\sigma = {4/3}.$ Hence, there is only one stretching transformation to consider: $\tau = x \varepsilon^{-4/3}$.
We now wish to compute formal solutions and the sequence $[\rho]$ with the matricial representation using our algorithm: 
Let $F = [f, \varepsilon^2 \partial f, \varepsilon^4 \partial^2 f, \varepsilon^6 \partial^3 f]^T$ then~\eqref{exm40} can be expressed as the following first order differential system ($\sigma_A=0$)
\begin{equation} \label{exm41} [A_{\sigma_A}] \quad \quad \quad  \xi_A^2 \partial F = A(x, \xi_A) F = \begin{bmatrix} 0 & 1 & 0 &0 \\ 0 & 0 & 1 &0  \\ 0 & 0 & 0 &1  \\ -x^3 & 0 & \xi_A^3 + 2x &0  \end{bmatrix} F, \quad \text{where} \; \xi_A = \varepsilon. \end{equation}
\begin{itemize}
\item For the turning point treatment, we set $x=t^2$ and then apply the transformation $$T =  \left[ \begin {array}{cccc} 0&0&0&1\\ \noalign{\medskip}0&0&{t}^{2}&0
\\ \noalign{\medskip}{t}^{3}&0&0&0\\ \noalign{\medskip}0&{t}^{4}&0&0
\end {array} \right] $$ which yields the system
\begin{equation} \label{exm43} [\tilde{A}_{\sigma_{\tilde{A}}}] \quad \quad \quad  \quad \xi_{\tilde{A}}^2 t \partial_t G = \tilde{A}(t, \xi_{\tilde{A}}) G = \left[ \begin {array}{cccc} -3\,\xi_{\tilde{A}}^{2}&2&0&0\\ \noalign{\medskip}4
+2\,\xi_{\tilde{A}}^{3}{t}^{5/2}&-4\,\xi_{\tilde{A}}^{2}&0&-2\,t
\\ \noalign{\medskip}2&0&-2\,\xi_{\tilde{A}}^{2}&0\\ \noalign{\medskip}0&0&2\,t&0
\end {array} \right] 
 G.
\end{equation}
where $\sigma_{\tilde{A}} = -3/2$, i.e. $\xi_{\tilde{A}} = t^{-3/2} \varepsilon$. 

 Applying the Splitting lemma decouples~\eqref{exm43} into two subsystems (truncated at order $4$ in both $t$ and $\varepsilon$)  
 $$[S^1_{\sigma_{S^1}}]\quad \xi_{S^1}^2 t \partial_t U^1 = S^1(t,\xi_{S^1}) U^1 \quad \text{and} \quad [S^2_{\sigma_{S^2}}]\quad \xi_{S^2}^2 t \partial_t U^2 = S^2(t,\xi_{S^2}) U^2$$ 
 where $\sigma_{S^1} = \sigma_{S^2} = -3/2$, $\xi_{S^1} = \xi_{S^2}=t^{-3/2} \varepsilon$, and 
 $$ S^1(t,\xi_{S^1})=  \left[ \begin {array}{cc} -1/2\, \xi_{S^1}^{2} \left( \xi_{S^1}^{4}{t}^{2}+{t}
^{2}+6 \right) &-\xi_{S^1}^{4}{t}^{2}+2\\ \noalign{\medskip}2\,\xi_{S^1}^{3}{t
}^{5/2}-\xi_{S^1}^{4}{t}^{2}-{t}^{2}+4&-1/2\,\xi_{S^1}^{2} \left( 3\,\xi_{S^1}^{4}
{t}^{2}+{t}^{2}+8 \right) \end {array} \right] 
$$
whose leading coefficient has two constant nonzero eigenvalues (hence the dimension of the system /order of the equivalent equation can be reduced).  And,
\begin{equation} \label{S2} S^2(t,\xi_{S^2})= \begin{bmatrix} 1/2\,\xi_{S^2}^{2} ( 3\,\xi_{S^2}^{4}{t}^{2}+{
t}^{2}-4) &\gamma_1\\ 1/2\,t ( 2\,\xi_{S^2}^{4}{t}^{2}+{t}^{2}+4) & \gamma_2 \end{bmatrix} , \end{equation}
 where $$ \begin{cases} \gamma_1 = -1/8\,t ( 4\,\xi_{S^2}^{3}{t}^{5/2}-9\,\xi_{S^2}^{4}{t}^{2}-8\,\xi_{S^2}^{4}-8) \\ \gamma_2 = -1/16\,{t}^{2}\xi_{S^2}^{2} ( 8\,\xi_{S^2}^{3}{t}^{5/2}-105\,\xi_{S^2}^{4}{t}^{2}-40\,\xi_{S^2}^{4}-14\,{t}^{2}-16 . \end{cases} $$ Since $\varepsilon t^{-3/2} = \varepsilon x^{-3/4}$, we have  $\rho_1 = 4/3$.\\
If we wish to continue the reduction for the second subsystem, then, as usual, we first consider its leading coefficient:
$$S^{2}_{0}(t) =  \left[ \begin {array}{cc} 0&t\\ \noalign{\medskip}2t+1/2\,{t}^{3}&0
\end {array} \right] .$$
Hence, by the turning point resolution (here it suffices to factorize $t$), the system $[S^2_{\sigma_{S^2}}]$ can be rewritten as:
$$[\tilde{S}^2_{\sigma_{\tilde{S}^2}}] \quad \quad \quad \quad \xi_{\tilde{S}^2}^2 t^2 \partial_t R = \tilde{S}^2(t,\varepsilon) R$$ 
where $\sigma_{\tilde{S}^2}=-2$, $\xi_{\tilde{S}^2} = t^{-2} \varepsilon$, and
$$ \tilde{S}^2(t,\xi_{\tilde{S}^2}) =  \left[ \begin {array}{cc} 1/2\,\xi_{\tilde{S}^2}^{2} \left( 3\,{t}^{4}\xi_{\tilde{S}^2}^{4}+{
t}^{2}-4 \right) &{\tilde{\gamma}}_1\\ \noalign{\medskip}{t}^{4}\xi_{\tilde{S}^2}^{4}+1/2\,{t}^{
2}+2& {\tilde{\gamma}}_2\end {array}
 \right] 
$$
where $$ \begin{cases} {\tilde{\gamma}}_1 = -1/2\,\xi_{\tilde{S}^2}^{3}{t}^{4}+{\frac {9\,{t}^{4}\xi_{\tilde{S}^2}^{4}}{
8}}+\xi_{\tilde{S}^2}^{4}{t}^{2}+1 \\ {\tilde{\gamma}}_2 = -1/16\,{t}^{2}\xi_{\tilde{S}^2}^{2} \left( -105\,{t}^{4} \xi_{\tilde{S}^2}^{4}+8\,\xi_{\tilde{S}^2}^{3
}{t}^{4}-40\,\xi_{\tilde{S}^2}^{4}{t}^{2}-14\,{t}^{2}-16 \right). \end{cases}$$
The leading constant coefficient has two constant eigenvalues and hence can be decoupled.  
\item We consider again system $[A_{\sigma_A}]$ which is given by~\eqref{exm41}, and we apply the stretching $x = \varepsilon^{4/3} \tau$ which yields the system $[E_{\sigma_E}]$. One can verify that $\sigma_E=0$ and so with $\xi_E = \varepsilon = \tilde{\xi}_E^3 = \tilde{\varepsilon}^3$, we have:
\begin{equation} \label{exm47} [E_{\sigma_E}] \quad \quad \quad \quad \tilde{\xi}_E^{2} \partial_{\tau} H =  \begin{bmatrix} 0 & 1 & 0 &0 \\ 0 & 0 & 1 &0  \\ 0 & 0 & 0 &1  \\ -2 \tau^3 \tilde{\xi}_E^{12} & 0 & \tilde{\xi}_E^9 + 2 \tau \tilde{\xi}_E^{4} &0  \end{bmatrix} H. \end{equation}
Then the $\tilde{\varepsilon}$-rank reduction with $$T =  \left[ \begin {array}{cccc} 0&0&0&1\\ \noalign{\medskip}0&0& \tilde{\xi}_E^{2}&0\\ \noalign{\medskip}0& \tilde{\xi}_E^{4}&0&0\\ \noalign{\medskip}\tilde{\xi}_E^{6}&0&0&0\end {array} \right] 
$$ yields
$$ [\tilde{E}_{\sigma_{\tilde{E}}}] \quad \quad \quad \quad \partial_{\tau} V = E(\tau, \tilde{\xi}_{\tilde{E}}) V = \left[ \begin {array}{cccc} 0&\tilde{\xi}_{\tilde{E}}^{15}+2\,\tau&0&-2\,{\tau}
^{3}\tilde{\xi}_{\tilde{E}}^{12}\\ \noalign{\medskip}1&0&0&0\\ \noalign{\medskip}0&1
&0&0\\ \noalign{\medskip}0&0&1&0\end {array} \right] 
V, \quad {\sigma_{\tilde{E}}}=0, \tilde{\xi}_{\tilde{E}} = \tilde{\varepsilon}.$$
Evidently, The $\varepsilon$-exponential part of $[\tilde{E}_{\sigma_{\tilde{E}}}]$ is zero. A formal fundamental matrix of solutions can be constructed following Subsection~\ref{hzero}.
\end{itemize}
\end{example}
\begin{example}[Roo's equation~\cite{key7863}]
Consider the singularly-perturbed scalar differential equation with $\sigma_a=0$ ($\xi_a = \varepsilon$) : \begin{equation} \label{eq1} [a_{\sigma_a}] \quad \quad \quad \xi_a^4 \; \partial^2 f = (x^5 + \xi_a x^2 + \xi_a^2) f .\end{equation}
Iwano-Sibuya's polygon of $[a_{\sigma_a}]$ consists of three segments connecting four points~\cite[p. 607]{key7863}, $P_{00} = (0,5/2)$, $P_{01} = (1/2,1)$, $P_{02} = (1,0)$, and $P_{o_\varepsilon} = (2,-1)$, and thus having the three slopes, $-3, -2 $, and $-1$. Consequently, Iwano-Sibuya's sequence $[\rho]$ of positive rational numbers is:
$$[\rho] \quad \quad \quad 0 < \rho_1 = 1/3 < \rho_2 = 1/2 < \rho_3 = 1 .$$
We now wish to recompute the sequence $[\rho]$ with the matrix representation of $[a_{\sigma_a}]$ using our algorithm.  Let $F = [f, \xi_a^2  \partial f]^T$ then $[a_{\sigma_a}]$ can be rewritten as the following first order differential system ($\sigma_A=0$)
\begin{equation} \label{sys1} [A_{\sigma_A}] \quad \quad \quad \xi_A^2  \partial F = A(x, \xi_A) F = \begin{bmatrix} 0 & 1  \\ x^5 + \xi_A x^2 + \xi_A^2 & 0 \end{bmatrix} F. \end{equation}
\begin{itemize}
\item For the resolution of the turning point, let $T = \begin{bmatrix} 1 & 0 \\ 0& x^{5/2} \end{bmatrix}$ (resolution of turning point) then $F = T G$ yields
\begin{equation} \label{sys12} [\tilde{A}_{\sigma_{\tilde{A}}}] \quad \quad \quad  \xi_{\tilde{A}}^2 x^{8} \partial G = \tilde{A}(x, \xi_{\tilde{A}})= \begin{bmatrix} 0 & 1  \\ 1 + \xi_{\tilde{A}} + \xi_{\tilde{A}} x & -\frac{5}{2}\xi_{\tilde{A}}^2x^{5/2} \end{bmatrix} G \end{equation}
with $\sigma_{\tilde{A}} = -3$ and $\xi_{\tilde{A}} = x^{-3} \varepsilon$.  The system  $[\tilde{A}_{\sigma_{\tilde{A}}}]$ can be decoupled into two scalar equations using the Splitting lemma and we have $\rho_1 = 1/3$.
\item 
We consider again the system $[A_{\sigma_A}]$ and we perform the stretching $\tau = x \varepsilon^{-1/3}$ and ramification $\xi_E = \varepsilon = \tilde{\varepsilon}^3 = \tilde{\xi}_E^3$. This yields:
$$ [E_{\sigma_E}] \quad \quad \quad \tilde{\xi}_E^{5}  \partial_{\tau} U = E(\tau, \tilde{\xi}_E) U =  \begin{bmatrix} 0 & 1  \\ \tau^5 \tilde{\xi}_E^5+ \tilde{\xi}_E^5 \tau^2 + \tilde{\xi}_E^6 & 0 \end{bmatrix} U, \quad \text{with} \; \sigma_E=0 .$$
Let $T_1 = \begin{bmatrix} 1 & 0 \\ 0& \tilde{\xi}_E^3 \end{bmatrix}$ ($\tilde{\varepsilon}$-rank reduction) then $U = T_1 G$ yields
$$ [\tilde{E}_{\sigma_{\tilde{E}}}] \quad \quad \quad \quad \tilde{\xi}_{\tilde{E}}^3 \partial_{\tau} G = \tilde{E}(\tau, \tilde{\xi}_{\tilde{E}}) G =  \begin{bmatrix} 0 & \tilde{\xi}_{\tilde{E}}  \\ \tau^5 + \tau^2 + \tilde{\xi}_{\tilde{E}} & 0 \end{bmatrix} G, \quad \text{with} \;\sigma_{\tilde{E}}=0.$$
The leading coefficient matrix is nilpotent. The computation of the $\tilde{\varepsilon}$-exponential order suggests introducing an additional ramification in $\varepsilon$ of order $2$. Since we already have a ramification of order $3$, to avoid a proliferation of notations, we reset  $\varepsilon = {\tilde{\varepsilon}}^{6}$. This yields:
$$ [\tilde{\tilde{E}}_{\sigma_{\tilde{\tilde{E}}}}] \quad \quad \quad \quad \tilde{\xi}_{\tilde{\tilde{E}}}^{6}  \partial_{\tau} G = \tilde{\tilde{E}}(\tau, \tilde{\xi}_{\tilde{\tilde{E}}}) G = \begin{bmatrix} 0 & \tilde{\xi}_{\tilde{\tilde{E}}}^{2}  \\ \tau^5 +\tau^2 + \tilde{\xi}_{\tilde{\tilde{E}}}^{2} & 0 \end{bmatrix} G ,$$
with $\sigma_{\tilde{\tilde{E}}}=0$, $\varepsilon = \xi_E = \tilde{\xi}_{\tilde{\tilde{E}}}^6 = \tilde{\varepsilon}^6$. Let $T_2 = \begin{bmatrix} \tilde{\xi}_{\tilde{\tilde{E}}} & 0 \\ 0&1 \end{bmatrix}$ then $G = T_2 W$ yields the following ${\tilde{\varepsilon}}$-irreducible system:
$$[J_{\sigma_J}] \quad \quad \quad \quad \tilde{\xi}_{J}^{5}  \partial_{\tau} W =  \begin{bmatrix} 0 &1  \\  \tau^5 +  \tau^2 + \tilde{\xi}_{J}^{2} & 0 \end{bmatrix} W, \quad \text{with} \; \sigma_J =0,\; \text{and}\;  \tilde{\xi}_{J}^6 = \tilde{\varepsilon}^6= \varepsilon.$$
Let $T_3 = \begin{bmatrix} 1 & 0 \\ 0& \tau \end{bmatrix}$ (resolution of turning point) then $W = T_3 V$ yields
$$[\tilde{J}_{\sigma_{\tilde{J}}}] \quad \quad \quad \quad \tilde{\xi}_{\tilde{J}}^{5} \tau^{4} \partial_{\tau} V = \tilde{B}(x, \tilde{\xi}_{\tilde{J}}) V = \begin{bmatrix} 0 &1 \\  \tau^3 +  1 + \tilde{\xi}_{\tilde{J}}^{2} &  -\tilde{\xi}_{\tilde{J}}^{5} \tau^{3}  \end{bmatrix} V $$
with $\sigma_{\tilde{J}} = -1$, $\tilde{\xi}_{\tilde{J}} = \tau^{-1} \tilde{\varepsilon} = \tau^{-1} \varepsilon^{1/6}$. Clearly, system  $[\tilde{J}_{\sigma_{\tilde{J}}}]$  can be decoupled into two  scalar equations. To find the restraining index, we rewrite $\tau$ and $\tilde{\varepsilon}$ in terms of $x$ and $\varepsilon$: 
$$\tau^{-1} \tilde{\varepsilon}  = \tau^{-1} \varepsilon^{1/6} = (x^{-1} \varepsilon^{1/3}) \varepsilon^{1/6}  = x^{-1} \varepsilon^{1/2} = (x^{-2} \varepsilon)^{1/2}, $$
and so $\rho_2 = 1/2$. 
\item 
We consider again system $ [A_{\sigma_A}]$ and perform the stretching $\tau = x \varepsilon^{-1/2}$, and the ramification $\varepsilon = \tilde{\varepsilon}^2$. This yields:
$$  [C_{\sigma_C}]  \quad \quad \quad \quad \tilde{\xi}_C^{3}  \partial_{\tau} Z = C(\tau, \tilde{\xi}_C) Z =  \begin{bmatrix} 0 & 1  \\ \tau^5 \tilde{\xi}_C^{5} + \tau^2 \tilde{\xi}_C^{4}+ \varepsilon^{2} & 0 \end{bmatrix} Z ,$$
with $\sigma_C =0$, $\xi_C = \varepsilon = \tilde{\varepsilon}^2 = \tilde{\xi}_C^2$. Let $T_1 = \begin{bmatrix} 1 & 0 \\ 0& \tilde{\xi}_C^2 \end{bmatrix}$ ($\tilde{\varepsilon}$-rank reduction) then $Z = T_1 G$ yields
$$  [\tilde{C}_{\sigma_{\tilde{C}}}]  \quad \quad \quad \quad \varepsilon^{1/2}  \partial_{\tau} G =  \begin{bmatrix} 0 & 1  \\ \tau^5 \varepsilon^{1/2} + \tau^2  + 1 & 0 \end{bmatrix} G, \quad \text{with} \;\sigma_{\tilde{C}}=0 \; \text{and} \; \tilde{\xi}_{\tilde{C}} = \tilde{\xi}_C,$$
which can be decoupled into two scalar equations by Splitting lemma. Since $\sigma_{\tilde{C}}=0$, the restraining index is infinity. 
\end{itemize}
And so, by our techniques, we have obtained the two slopes $\rho_1 = 1/3$ and  $\rho_2 = 1/2$.
\end{example}

\begin{example}[(\cite{key218}, Example 1, Section 9.5, p. 446)]
\label{exmbender1}
Consider the following singularly-perturbed scalar differential equation $$\varepsilon \partial^3 f - \partial f + x f = 0 $$
which can be rewritten as the following differential first order system
$$[A_{\sigma_A}] \quad \quad \quad \xi_A \partial F = A(x, \xi_A) F = \begin{bmatrix} 0 & \xi_A & 0 \\ 0 & 0 & \xi_A \\ -x & 1 & 0 \end{bmatrix} F$$
where $F = [f, \partial f , \partial^2 f]^T$, $\sigma_A=0$, and $\xi_A=\varepsilon$. Since
$A_0(x)$ is nilpotent, Splitting lemma cannot be applied but no treatment of turning points is required. 
The transformation $F = \left[ \begin {array}{ccc} 0&1&0\\ \noalign{\medskip}1&0&0
\\ \noalign{\medskip}0&0&1\end {array} \right] \; \left[ \begin {array}{ccc} 1&x&0\\ \noalign{\medskip}0&1&0
\\ \noalign{\medskip}0&0&1\end {array} \right]  \; G$ yields the following $\varepsilon$-irreducible equivalent system:
$$[\tilde{A}_{\sigma_{\tilde{A}}}] \quad \quad \quad \xi_{\tilde{A}} \partial G = \tilde{A}(x,\xi_{\tilde{A}}) G =  \left[ \begin {array}{ccc} -\xi_{\tilde{A}} \,x&-{x}^{2}\xi_{\tilde{A}}-\xi_{\tilde{A}}&
\xi_{\tilde{A}}\\ \noalign{\medskip}\xi_{\tilde{A}}&\xi_{\tilde{A}}\,x&0
\\ \noalign{\medskip}1&0&0\end {array} \right]G ,
$$
with $\sigma_{\tilde{A}} = 0$ and $\theta_{\tilde{A}}(\lambda) = -\lambda - x$.
Since the leading matrix coefficient is still nilpotent, we have to compute the $\varepsilon$-exponential oder. The characteristic polynomial of $\frac{\tilde{A}(x, \xi_{\tilde{A}})}{\xi_{\tilde{A}}}$ is given by:
$$\lambda^3 + \frac{\xi_{\tilde{A}} - 1}{\xi_{\tilde{A}}} \lambda + \frac{x}{\xi_{\tilde{A}}} .$$
Consequently, the $\varepsilon$-exponential oder is given by:
$\omega_\varepsilon (\tilde{A}) = \; max \; \{ \frac{1}{2} , \frac{1}{3} \}  = \frac{1}{2} .$
Let $\varepsilon = \tilde{\varepsilon}^2$ and $\xi = \tilde{\xi}^2$ then we have:
$$[\tilde{\tilde{A}}_{\sigma_{\tilde{\tilde{A}}}}] \quad \quad \quad   \tilde{\xi}^2 \partial G = \tilde{\tilde{A}}(x, \tilde{\xi}) G =  \left[ \begin {array}{ccc} - \tilde{\xi}^2\,x&-{x}^{2} \tilde{\xi}^2- \tilde{\xi}^2&
 \tilde{\xi}^2\\ \noalign{\medskip} \tilde{\xi}^2& \tilde{\xi}^2\,x&0
\\ \noalign{\medskip}1&0&0\end {array} \right] G, \quad \text{with}\; \sigma_{\tilde{\tilde{A}}} = 0.
$$
The transformation $G = \begin{bmatrix}  \tilde{\xi} & 0 & 0 \\ 0 & \tilde{\xi} & 0 \\ 0 & 0 & 1 \end{bmatrix} \left[ \begin {array}{ccc} 0&1/2&1/2\\ \noalign{\medskip}1&0&0
\\ \noalign{\medskip}0&-1/2&-1/2\end {array} \right] 
 U$ yields the following $\tilde{\varepsilon}$-irreducible system:
  $$[B_{\sigma_B}] \quad \quad \quad  \tilde{\xi}_B \partial U = B(x, \tilde{\xi}_B) U = \left[ \begin {array}{ccc} \tilde{\xi}_B\,x&\tilde{\xi}_B/2&\tilde{\xi}_B/2
\\ \noalign{\medskip} \left( -{x}^{2}-1 \right) \tilde{\xi}_B&-1/2\,
\tilde{\xi}_B\,x-1&-1/2\,\tilde{\xi}_B\,x\\ \noalign{\medskip} \left( -{x}^{2}-1
 \right) \tilde{\xi}_B&-1/2\,\tilde{\xi}_B\,x&-1/2\,\tilde{\xi}_B\,x+1\end {array}
 \right] 
 U,
$$
with $\sigma_B = 0, \tilde{\xi}_B = \tilde{\varepsilon}$.
The leading matrix coefficient $B_{0,0}$ has three distinct roots and consequently, the system can be decoupled into three first order scalar differential equations by applying the Splitting lemma. A fundamental matrix of formal solutions is then given by:
$$\Phi (x, {\varepsilon}^{1/2}) \exp (\begin{bmatrix} 0  & 0 & 0 \\ 0 & \frac{-x}{{\varepsilon}^{1/2}} & 0 \\ 0 & 0& \frac{x}{{\varepsilon}^{1/2}} \end{bmatrix} ) .$$ 
where $\Phi(x, {\varepsilon}^{1/2}) $ is the product of all the transformations performed including that of the Splitting lemma. Since $\sigma_B = 0$, the origin is not a turning point for this system and we do not need stretchings.
\end{example}

\begin{example}[\cite{key221}]
\label{exmiwano}
We consider the scalar differential equation
$$[a_{\sigma_a}] \quad \quad \quad \varepsilon^3 \partial^2 f  - (x^3 + \varepsilon) f = 0, $$
with $\sigma_a=0$. Setting $F = [f, \partial f]^T$, we get the following differential first order system
$$[A_{\sigma_A}] \quad \quad \quad  \xi_A^3 \partial F = A(x, \xi_A) F = \begin{bmatrix} 0 & \xi_A^3  \\ x^3 + \xi_A & 0 \end{bmatrix} F$$
with $\sigma_A=0$  and $\xi_A = \varepsilon$. 
\begin{itemize}
\item Since $A_{0,0}$ is nilpotent and $A_0(x)$ is not, we first start with the treatment of the turning point at $x=0$. It suffices however to factorize $x^3$ from $A_0(x)$ which results in $\sigma_A =-3$ and $\xi_A= x^{-3} \varepsilon$. We then apply the transformation $F = \begin{bmatrix}  \xi_A & 0 \\ 0  & 1 \end{bmatrix}  G$ to get the $\varepsilon$-irreducible system: 
$$ [\tilde{A}_{\sigma_{\tilde{A}}}] \quad \quad \quad \xi_{\tilde{A}}^2 x^6 \partial G = \tilde{A}(x, \xi_{\tilde{A}})  G = \begin{bmatrix} 3 x^5 \xi_{\tilde{A}}^2 & \xi_{\tilde{A}} x^6 \\ 1 + \xi_{\tilde{A}} &  0 \end{bmatrix} G$$
with $\sigma_{\tilde{A}} = -3$ and $\xi_{\tilde{A}}= x^{-3} \varepsilon$. 
Since $\tilde{A}_0(x)$ is still nilpotent, we proceed to computing the $\varepsilon$-exponential order from the characteristic polynomial of $\frac{\tilde{A}(x, \varepsilon)}{\xi_{\tilde{A}}^2 x^6}$ which is given by:
$${\lambda}^{2}-\frac{3}{x} \lambda - \frac{1+\xi_{\tilde{A}}}{ \xi_{\tilde{A}}^3 x^6}.
$$
Hence, $\omega_{\varepsilon}(\tilde{A}) = \frac{3}{2}$ and the $\varepsilon$-polynomial is given by: $X^2 -x^3 =0$. Let $\varepsilon= \tilde{\varepsilon}^2$ and $\tilde{\xi}_{\tilde{A}} = x^{-3} \tilde{\varepsilon}$, then with $\xi_{\tilde{A}} = \tilde{\xi}_{\tilde{A}}^2 x^3$ we have:
$$  [{\tilde{A}}_{\sigma_{\tilde{A}}}] \quad \quad \quad x^{12} \tilde{\xi}_{\tilde{A}}^4 \partial G = \tilde{A}(x,\tilde{\xi}_{\tilde{A}})  G = \begin{bmatrix} 3x^{11} \tilde{\xi}_{\tilde{A}}^4 & \tilde{\xi}_{\tilde{A}}^2 x^9 \\ 1 + x^3\tilde{\xi}_{\tilde{A}}^2 &  0 \end{bmatrix} G$$
with $\sigma_{\tilde{A}} = -3$.  
We then perform $G = \begin{bmatrix}  \tilde{\xi}_{\tilde{A}} & 0 \\ 0  & 1 \end{bmatrix}  U$ to get the $\tilde{\varepsilon}$-irreducible system: 
$$  [B_{\sigma_{B}}] \quad \quad \quad  x^{12}  {\tilde{\xi}}_B^3 \partial U = B(x, {\tilde{\xi}}_B)  U = \begin{bmatrix} 3 x^{11} {\tilde{\xi}}_B^2 ( {\tilde{\xi}}_B +1) & x^{9}  \\ 1 + x^3 {\tilde{\xi}}_B^2 &  0 \end{bmatrix} U, $$
with $\sigma_B =-3$ and ${\tilde{\xi}}_B = x^{-3} \tilde{\varepsilon}$.

Now that $B_{0}(x)$ has two distinct eigenvalues and $B_{0,0}$ is nilpotent,
we perform $U= \begin{bmatrix}  x^{9/2} & 0 \\ 0  & 1 \end{bmatrix}  W$ to treat the turning point. Then we have: 
$$ [\tilde{B}_{\sigma_{\tilde{B}}}] \quad \quad \quad  x^{15/2} {\tilde{\xi}}_{\tilde{B}}^3 \partial W = \tilde{B}(x, {\tilde{\xi}}_{\tilde{B}})  W = \begin{bmatrix} 3 x^{13/2} {\tilde{\xi}}_{\tilde{B}}^2 (1 - \frac{1}{2}{\tilde{\xi}}_{\tilde{B}})  & 1 \\ 1 + x^3 {\tilde{\xi}}_{\tilde{B}}^2 &  0 \end{bmatrix} W ,$$
with $\sigma_{\tilde{B}} =-3$ and ${\tilde{\xi}}_{\tilde{B}} = x^{-3} \tilde{\varepsilon}$.

Since $\tilde{B}_{00}$ has two distinct eigenvalues, $\pm1$, one can proceed to the Splitting lemma which decouples the system thoroughly. From the the eigenvalues of $\tilde{B}_{00}$ or the $\varepsilon$-polynomial, one can read the $\varepsilon$-exponential part of an outer solution  of this system:
$$\exp (\begin{bmatrix} \frac{\int^x z^{3/2} dz}{\varepsilon^{3/2}} & 0 \\ 0 & -\frac{\int^x z^{3/2} dz}{\varepsilon^{3/2}}  \end{bmatrix}) ,$$ and deduce that $\rho_1= 1/3$.

\item We consider again system $[A_{\sigma_A}]$ and perform the stretching $\tau = x \varepsilon^{-1/3}$. After introducing the ramification $\varepsilon = \tilde{\varepsilon}^3$, we have:
$$[S_{\sigma_S}] \quad \quad \quad \tilde{\xi}_S^{5} \partial_{\tau} U = S(\tau, \tilde{\xi}_S)= \begin{bmatrix} 0 & \tilde{\xi}_S^6 \\ \tau^3 + 1& 0 \end{bmatrix}  U $$
with $\sigma_S = 0$, $ \tilde{\xi}_S = \tilde{\varepsilon}$, $\xi_S = \tilde{\xi}_S^3$. The transformation $U = \begin{bmatrix}  \tilde{\xi}_S^3 & 0 \\ 0  & 1 \end{bmatrix}  V$ yields 
$$[\tilde{S}_{\sigma_{\tilde{S}}}] \quad \quad \quad \tilde{\xi}_{\tilde{S}}^{2}  \partial_{\tau}V =  \tilde{S} (\tau, \tilde{\xi}_{\tilde{S}}) V = \begin{bmatrix} 0 & 1 \\ \tau^3 + 1& 0 \end{bmatrix}  V ,$$
with $\sigma_{\tilde{S}} = 0$, $ \tilde{\xi}_{\tilde{S}} = \tilde{\varepsilon}$.
Due to the nature of the eigenvalues of $\tilde{S}_{0,0}$, the system $[\tilde{S}_{\sigma_{\tilde{S}}}]$ can be decoupled into two scalar equations. $\sigma_{\tilde{S}} = 0$ suggests that there are no more stretchings to perform.
\item However, if we consider again system $[A_{\sigma_A}]$ and experiment with the stretching $\tau = x \varepsilon^{-1}$, we get:
$$[S_{\sigma_S}] \quad \quad \quad \xi_S \partial_{\tau} U = S(\tau, \xi_S) U=  \begin{bmatrix} 0 & \xi_S^2 \\ \tau^3 \xi_S^2 + 1 & 0 \end{bmatrix} U $$
with $\sigma_S = 0$ and $\xi_S = \varepsilon$. By the transformation $U = \begin{bmatrix}  \varepsilon & 0 \\ 0  & 1 \end{bmatrix}  V$ and some simplification we get
$$[\tilde{S}_{\sigma_{\tilde{S}}}] \quad \quad \quad \partial_{\tau} G=  \tilde{S} (\tau, \xi_{\tilde{S}}) V = \begin{bmatrix} 0 & 1 \\ \tau^3 \xi_{\tilde{S}}^2 + 1& 0 \end{bmatrix} V ,$$
with $\sigma_{\tilde{S}} = 0$ and $\xi_{\tilde{S}} = \varepsilon$.
Due to the nature of the eigenvalues of $\tilde{S}_{0,0}$, the system $[\tilde{S}_{\sigma_{\tilde{S}}}]$ can be decoupled into two scalar equations. Its $\varepsilon$-exponential part is zero.
\end{itemize}
So far, this matricial approach can determine the stretching $\tau = x \varepsilon^{-1/3}$ but not the stretching $\tau = x \varepsilon^{-1}$. However, both stretchings can be obtained from Iwano-Sibuya's polygon by the treatment of this system as a scalar equation~\cite[Last section]{key221}.
\end{example}
\begin{example}[\cite{key218}, Example 1, Section 9.5, p. 453]
\label{exmbender2}
We consider the scalar equation $$ \varepsilon^3 x \partial^2 f + x^2 \partial f - (x^3 + \varepsilon) f= 0 ,$$
which can be rewritten as the following first order differential system
$$[A_{\sigma_A}] \quad \quad \quad x \xi_A^3 \partial F = A(x, \xi_A) F = \begin{bmatrix} 0 & \xi_A^3 x \\ x^3 + \xi_A &  - x^2 \end{bmatrix} F$$
where $\sigma_A=0$, $\xi_A = \varepsilon$, and $F = [f, \partial f]^T$. 
We first start with the treatment of the turning point at $x=0$ since $A_{0,0}$ is nilpotent while $A_0(x)$ is not. It suffices here to factorize $x^2$ which yields 
$$[\tilde{A}_{\sigma_{\tilde{A}}}] \quad \quad \quad x^{5} \xi_{\tilde{A}}^3 \partial F = \tilde{A}(x, \xi_{\tilde{A}} )  F = \begin{bmatrix} 0 &  \xi_{\tilde{A}}^3 x^5 \\ x +  \xi_{\tilde{A}}  &  - 1 \end{bmatrix} F$$
where $ \sigma_{\tilde{A}} = -2$ and $\xi_{\tilde{A}} = x^{-2} \varepsilon$. Now that $\tilde{A}_{0,0}$ has two distinct eigenvalues, we can proceed to the Splitting lemma. The transformation $F = \begin{bmatrix} 0 & 1 \\ 1 & x\end{bmatrix} G$ block-diagonalizes $\tilde{A}_0(x)$ which yields
$$[\tilde{\tilde{A}}_{\sigma_{\tilde{\tilde{A}}}}] \quad \quad \quad x^{5} \xi_{\tilde{\tilde{A}}}^3 \partial G = \tilde{\tilde{A}}(x, \xi_{\tilde{\tilde{A}}})  G = \begin{bmatrix} - \xi_{\tilde{\tilde{A}}}^3 x^6 -1 & -\xi_{\tilde{\tilde{A}}}^3 x^7 - \xi_{\tilde{\tilde{A}}}^3 x^5 + \xi_{\tilde{\tilde{A}}} \\  \xi_{\tilde{\tilde{A}}}^3 x^5  &   \xi_{\tilde{\tilde{A}}}^3 x^6  \end{bmatrix} G ,$$
where $ \sigma_{\tilde{\tilde{A}}} = -2$ and $\xi_{\tilde{\tilde{A}}} = x^{-2} \varepsilon$ and from which we can deduce that the $\varepsilon$-exponential part is given by:
$$\exp (\begin{bmatrix}  \frac{- x^2}{ 2 {\varepsilon}^{3}} + O(ln(x) \varepsilon^{-2}) & 0  \\ 0 & O( ln(x) \varepsilon^{-2}) \end{bmatrix} ) .$$ 
Since the restraining index is nonzero, we can deduce that $\rho_1 = 1/2$ and apply the stretching $\tau = x \varepsilon^{-1/2}$ in the input system $[A_{\sigma_A}]$ and repeat the formal reduction procedure. 
\end{example}

\begin{example}
\label{introexm2}
In this example, we use our algorithm to explain some of the computations in the introductory example in the light of the eigenvalues of the matrix of the system. Thus, we consider again the system $[A_{\sigma_A}]$ given by~\eqref{introexmsys}:
$$[A_{\sigma_A}] \quad \quad \quad \varepsilon \; \partial F = A(x, \varepsilon) F= \begin{bmatrix} 0 & 1 \\ x^3 - \varepsilon & 0 \end{bmatrix} F. $$

We remark that the eigenvalues of $A(x, \varepsilon)$ are given by $\lambda_{1,2} = \pm (x^3 - \varepsilon)^{1/2}$. We thus have the two possible expansions which correspond to the two subdomains above:
\begin{itemize}
\item $\lambda_{1,2} = \pm (x^3 - \varepsilon)^{1/2} = \pm i  \varepsilon^{1/2} (1 - \frac{1}{2}\varepsilon^{-1} x ^{3} - \frac{1}{8} \varepsilon^{-2} x ^{6} + \dots )$, valid for $|\varepsilon^{-1} x^{3}| < 1$.
\item 
$\lambda_{1,2} = \pm (x^3 - \varepsilon)^{1/2} = \pm x^{3/2} (1 - \frac{1}{2}\varepsilon x ^{-3} - \frac{1}{8} \varepsilon^2 x^{-6} + \dots )$ , valid for $|\varepsilon x^{-3}| < 1$. 
\end{itemize}
\end{example}

\end{document}